\documentclass[11pt, a4paper]{article}
\usepackage[comma, authoryear, numbers]{natbib}
\usepackage{amsmath,amsfonts,amsthm,bm}
\usepackage[parfill]{parskip}
\usepackage{hyperref}
\usepackage{graphicx}
\usepackage[a4paper, margin=25mm]{geometry}
\usepackage{bm}
\usepackage{bbm}
\usepackage{breqn}
\usepackage[toc,page]{appendix}
\usepackage{tabularx}
\usepackage{amsthm}
\usepackage[dvipsnames]{xcolor}
\usepackage{multirow}
\usepackage{longtable}
\usepackage{booktabs}
\usepackage{bibentry}
\usepackage[multiple]{footmisc}
\usepackage{mathtools}
\usepackage{algpseudocode}
\usepackage{algorithmicx}
\usepackage{algorithm}
\usepackage{subcaption}
\usepackage{etoolbox}
\usepackage{dsfont}
\usepackage{breqn}
\usepackage{tabulary}
\usepackage{tocbibind}
\usepackage{comment}

\DeclareMathOperator*{\argmax}{argmax} 

\newcommand{\fa}[0]{\ \forall \ }

\newcommand{\E}[0]{\mathbb{E}}
\newcommand{\R}{\mathbb{R}}
\renewcommand{\P}[0]{\mathbb{P}}
\newcommand{\m}[1]{\mathcal{#1}}

\renewcommand{\l}{\left}
\renewcommand{\r}{\right}

\usepackage{wrapfig}
\mathtoolsset{showonlyrefs}

\allowdisplaybreaks

\title{ 
    Robust Markov decision processes under parametric transition distributions
}
\author{
    Ben Black\footnote{STOR-i Centre for Doctoral Training, Lancaster
    University, United Kingdom. Email:
    \href{mailto:b.black1@lancaster.ac.uk}{b.black1@lancaster.ac.uk}}
    \footnote{Corresponding author.}\ , 
    Trivikram Dokka\footnote{Advanced Analytics Group, Air Products Plc, United Kingdom. Email:
    \href{mailto:Trivikram.Dokka@yahoo.co.uk}{Trivikram.Dokka@yahoo.co.uk}}\ , 
    Christopher Kirkbride\footnote{Department of Management Science, Lancaster
    University Management School, United Kingdom. Email:
    \href{mailto:c.kirkbride@lancaster.ac.uk}{c.kirkbride@lancaster.ac.uk}.}
}
\date{}
\graphicspath{{images^tables/}}

\usepackage{amsthm}
\usepackage[dvipsnames]{xcolor}

\begin{document}
\setcitestyle{nosort}
\pagenumbering{roman}
\newpage
\maketitle
\pagenumbering{arabic}

\begin{abstract}
    This paper considers robust Markov decision processes under parametric transition distributions. We assume that the true transition distribution is uniquely specified by some parametric distribution, and explicitly enforce that the worst-case distribution from the model is uniquely specified by a distribution in the same parametric family. After formulating the parametric robust model, we focus on developing algorithms for carrying out the robust Bellman updates required to complete robust value iteration. We first formulate the update as a linear program by discretising the ambiguity set. Since this model scales poorly with problem size and requires large amounts of pre-computation, we develop two additional algorithms for solving the robust Bellman update. Firstly, we present a cutting surface algorithm for solving this linear program in a shorter time. This algorithm requires the same pre-computation, but only ever solves the linear program over small subsets of the ambiguity set. Secondly, we present a novel projection-based bisection search algorithm that completely eliminates the need for discretisation and does not require any pre-computation. We test our algorithms extensively on a dynamic multi-period newsvendor problem under binomial and Poisson demands. In addition, we compare our methods with the non-parametric phi-divergence based methods from the literature. We show that our projection-based algorithm completes robust value iteration significantly faster than our other two parametric algorithms, and also faster than its non-parametric equivalent. 
\end{abstract}

\normalsize 

\section{Introduction}

Markov decision processes~\citep{puterman} (MDPs) are a mathematical framework for modelling dynamic decision making problems under uncertainty. Under the MDP framework, at each decision epoch in a finite or infinite time horizon, a decision maker utilises information about the current \textit{state} of a system in order to select an \textit{action} that yields them a reward. The action taken can affect the next state of the system, which is stochastically governed by a set of \textit{transition probabilities}. The goal of the decision maker is to make decisions at each epoch in order to maximise the total (discounted) expected reward that they receive over the entire horizon. A solution of an MDP is understood as a \textit{policy}, which provides an action or a distribution over the set of actions to be taken in each state of the system. The policy is found prior to any decisions being made, and in practice the decision maker can instantaneously generate their action from the policy at any given epoch. Policies are usually found from algorithms based on dynamic programming and Bellman's optimality equations~\citep{bellman1966dynamic}, which are based around the concept of \textit{value functions}. Value functions give the expected total future reward from starting in each state and following an optimal policy thereafter.

In classical MDPs, it is assumed that all parameters of the model (rewards, transition probabilities, etc.) are known exactly. However, in practice it can be difficult to determine these parameters exactly and they must often be replaced with estimates. However, it has been found that replacing true parameters with estimates thereof can lead to policies that fail drastically when implemented, due to errors in estimation~\citep{le2007robust, Wiesemann13} and that the resulting value function estimates can have large variance and bias~\citep{Mannor07}. Due to these issues, \textit{robust MDPs}~\citep{Wiesemann13} (RMDPs) have been proposed to explicitly represent uncertainty in model parameters. RMDPs do not assume that all parameters are known, but that they are known to lie in some pre-determined set. The decision maker then aims to find a policy with the best worst-case total reward over all parameters in the set. This limits the potential hazards of poor estimation.

In the case where only the transition probabilities are not known, we refer to this set as an \textit{ambiguity set}. Ambiguity sets are designed so that the decision maker can be confident that the true transition distribution lies within the set. There are many ways of constructing an ambiguity set. Early sets placed bounds on each transition probability~\citep{Satia73, GIVAN2000}. In more recent papers, it has become more common to bound the distance between any distribution in the set and some nominal distribution. For example, one can use the Kullback-Leibler divergence, modified $\chi^2$-distance or $L_1$-norms~\citep{Iyengar05}, or more general classes of distance measures such as $\phi$-divergence functions~\citep{ho2022robust}. The choice of ambiguity set strongly affects the tractability of the resulting RMDP model. For general ambiguity sets, it is known that RMDPs are NP-hard. However, this paper considers a special type of ambiguity set called $s$-rectangular ambiguity sets~\citep{le2007robust}. Such ambiguity sets allow the transition distributions for each state to be chosen independently of one another. The resulting RMDP is solvable in polynomial time via \textit{robust value iteration}~\citep{Wiesemann13}.

Robust value iteration starts with some initial estimate of the value functions, then iteratively updates these estimates until Bellman's optimality equations are satisfied. We refer to the process of finding the next value function estimates as solving a \textit{robust Bellman update}. Much of the recent RMDP literature has focused on developing fast algorithms for solving the robust Bellman update in $s$-rectangular RMDPs. Due to the fact that only the value estimates themselves (and not the optimal policies responsible) are required to complete robust value iteration, many of these algorithms employ simple methods like bisection search~\citep{Grand_Clement_Kroer_2021, ho2022robust} to solve the update. In this paper, we will focus on developing algorithms for carrying out robust Bellman updates, with one key difference from the existing literature. In particular, we will focus on RMDPs where the true state transition distribution either lies in some parametric family or is specified by some external random variable that lies in some parametric family (e.g. demand, service times, failure rates). In the existing RMDP literature, transition distributions are assumed to be non-parametric, and ambiguity sets typically contain non-parametric distributions. However, in the case where the transition distribution is indeed parametric, non-parametric ambiguity sets necessarily contain distributions that cannot be equal to the true distribution. Our models use \textit{parametric ambiguity sets} to enforce that every potential distribution in the set lies in the same parametric family as the true distribution. 

Constructing an RMDP in this fashion has a number of benefits. Firstly, it means that we only need to find the worst-case parameters and not the entire worst-case distribution. The worst-case parameter is typically of much smaller dimension than the worst-case distribution, meaning that finding it can be much less cumbersome. As such, instead of ambiguity sets for the true distribution, we use ambiguity sets for the true parameters. Secondly, explicitly using ambiguity sets for the true parameters and using the corresponding parametric distributions in the model means that every worst-case distribution generated by the model will lie in the correct parametric family. In addition, we can make use of maximum likelihood estimation to build confidence sets for the true parameters and use these as ambiguity sets in our models. Finally, parametric distributions are natural models for random variables affecting MDP state transitions in a number of problems. An example of such a problem is a dynamic multi-period newsvendor problem~\citep{arrow1958}. In newsvendor models, demand is often considered as a parametric random variable. More specifically, newsvendor demand has been modelled as normal~\citep{Nahmias1994739}, negative binomial~\citep{Agrawal96}, lognormal and gamma~\citep{GALLEGO2007}, and exponential~\citep{Siegel21}. In addition, for such problems it has been shown that assuming that parameter estimates are truth can lead to poor cost estimation~\citep{ROSSI2014674, Siegel21}. Hence, a parametric ambiguity set provides a way to hedge against parameter uncertainty while ensuring the worst-case distribution is also parametric.

This paper extends the concept of parametric ambiguity sets from~\cite{black}, who studied a static multi-period resource planning problem under binomial demand, into the RMDP literature. We formulate $s$-rectangular parametric ambiguity sets and solve the resulting RMDP via robust value iteration. Under such ambiguity sets, the robust Bellman update is a parametric distributionally robust optimisation problem~\citep{black}. As a benchmark, we reformulate the robust Bellman update as a linear program (LP) by discretising the ambiguity set. Since this LP can become very slow for large problems, we develop two additional algorithms for carrying out the update. The first is a \textit{cutting surface}~\citep{mehrotra2014cutting} (CS) algorithm that iteratively solves the LP over increasing subsets of the ambiguity set. The second algorithm is a parametric projection-based bisection search algorithm. This algorithm does not rely on any discretisation of the ambiguity set, and we will show that this means that it solves the robust Bellman updates orders of magnitude faster than both CS and LP.

In summary, the contributions of this paper are as follows:
\begin{enumerate}
    \item We extend the concept of parametric ambiguity sets from~\cite{black} into the RMDP literature. Such ambiguity sets have only been used for static distributionally robust optimisation (DRO) problems in the past. Since the DRO model used in the robust Bellman update must be solved multiple times in an iterative fashion, scalability of algorithms and computation is even more of a challenge in RMDPs.
    \item We develop a fast projection-based bisection search algorithm for solving a robust Bellman-update, that does not rely on any discretisation of the parametric ambiguity set.
    \item We apply our methods to a dynamic multi-period newsvendor model under binomial and Poisson demands. The results show that the parametric robust value iteration is tractable and can be solved faster than its non-parametric equivalent.
\end{enumerate}


\section{Literature review}

\subsection{Robust Markov decision processes}

RMDPs are a framework for modelling MDPs with unknown parameters, which has become common in recent years due to the fact that MDPs are extremely sensitive to small changes in their parameters~\citep{Mannor07}. RMDPs have been studied in the literature since the 1970s, where the first example of an RMDP used ambiguity sets based on assigning upper and lower bounds to each transition probability~\citep{Satia73}. Such ambiguity sets were common in the early MDP literature. \cite{GIVAN2000} also studied bounded parameter RMDPs, which were solved by solving a collection of exact MDPs. Later, \cite{bagnell2001solving} generalised the concepts of RMDPs to a variety of other ambiguity sets. Their only assumption was that the sets were convex and compact, meaning that the class they considered covered interval ambiguity sets as a special case. In addition, finite horizon RMDPs were also studied, bringing forth a robust version of dynamic programming (DP)~\citep{Nilim05}. These authors further developed the ambiguity sets used to encompass distance-based sets, such as those based on the Kullback-Leibler divergence. Using such ambiguity sets allowed the Bellman optimality equations to be reformulated using dualisation and hence solved exactly or via bisection. \cite{Iyengar05} formalised these concepts further, studying finite and infinite horizon RMDPs with a variety of ambiguity sets. For example, they studied ambiguity sets built using the Kullback-Leibler divergence, modified $\chi^2$-distance and $L_1$ norm.

Since these early papers, $s$-rectangular ambiguity sets have become very common in RMDPs. An $s$-rectangular ambiguity set~\citep{le2007robust} is one arising from the situation in which the state transitions for each state are independent of one another. Hence, the worst-case distributions for each state can be extracted independently of one another. Solving an RMDP with an $s$-rectangular ambiguity set is equivalent to finding a fixed point of the robust Bellman operator~\citep{Wiesemann13}, hence allowing a robust value iteration algorithm to solve the infinite horizon case. Many recent papers have developed fast algorithms for solving the robust Bellman updates required by robust value iteration. For example, \cite{Ho21} studied $s$-rectangular ambiguity sets defined by the $L_\infty$ norm. They developed a homotopy method that was implemented within a bisection search algorithm for solving the robust Bellman update. \cite{ho2021partial} applied a similar concept to weighted $L_1$ norm ambiguity sets, although they used a partial policy iteration algorithm instead of value iteration. For ellipsoidal and Kullback-Leibler ambiguity sets, \cite{Grand_Clement_Kroer_2021} proposed a first order method that is embedded in robust value iteration. Their algorithm is based on the observation that solving the robust Bellman update is equivalent to solving $S$ bilinear saddle point problems. \cite{ho2022robust} studied $\phi$-divergence ambiguity sets, and showed that solving the update in this case corresponds to solving a set of highly structured simplex projection problems. They used dualisation to represent each projection problem as a univariate convex optimisation problem. Different to these algorithms, \cite{derman2021twice} showed that solving an $s$-rectangular RMDP with reward uncertainty is equivalent to solving a regularised MDP.

In general, $s$-rectangular ambiguity sets are common due to the tractability of the resulting RMDP. However, since the state transitions for different states are not always independent, more general ambiguity sets have also been presented in the literature. \cite{NEURIPS2018_7ec0dbee} state that $s$-rectangular ambiguity sets can lead to conservative policies, and does not facilitate knowledge transfer between states or across different decision processes. They instead use non-rectangular ambiguity sets that bound the moments of state-action features, which are taken over entire MDP trajectories and not just those for one state. These RMDPs are solved by finding the optimal policy for a mixture of non-robust MDPs. Following a similar argument with regards to the conservativeness of $s$-rectangular ambiguity set, \cite{goyal2022robust} develop a new class of non-rectangular ambiguity sets: factor matrix ambiguity sets. Each distribution in such an ambiguity set is a convex combination of a set of common feature vectors. This ambiguity set allows for the modelling of dependence across states and for the RMDP to be efficiently solved by a hybrid value iteration-policy improvement algorithm. 

This paper studies RMDPs with $s$-rectangular ambiguity sets, but with one key difference from those discussed here. We study transition distributions that are parametric, and our ambiguity sets contain only distributions that lie in the same parametric family as the true transition distribution. This represents, for example, problems in which the state transitions are defined by some external random variables such as demand, and that these random variables take parametric distributions. Despite the fact that the underlying transition distributions may be parametric, in RMDPs, ambiguity sets always contain non-parametric distributions. However, any distribution in the set that is not part of the same family as the true transition distribution cannot be equal to the true distribution. As a result of this, we formulate \textit{parametric RMDPs}, where the ambiguity sets used contain potential parameters of the transition distribution, not potential distributions. This allows us to ensure that the worst-case distribution lies in the correct parametric family, and in addition we only need to find the worst-case parameter, not the entire distribution. To solve the resulting RMDP, we present three algorithms that are used inside a robust value iteration. The first two are based on discretising the ambiguity set of parameters and formulating the update as a linear program with one constraint for each parameter. The second is a fast bisection search algorithm that solves simplex projection problems to compute the update, similar to the approach of~\cite{ho2022robust}. 

\subsection{Newsvendor models}

The model that we will use to illustrate our methods is the newsvendor model~\citep{Newsboy1}. The newsvendor model is a classical model in inventory and operations management that describes a retailer deciding on how much stock to purchase in order to meet uncertain future demand as closely as possible. The newsvendor model has the distinguishing feature that failing to meet demand in any way is penalised. If too much stock is purchased, the newsvendor pays a holding cost in order to keep that stock for future customers. If demand is not met by the stock purchased, the newsvendor pays a backorder cost in order to meet the unmet demand. Due to this, demand uncertainty and correctly modelling said uncertainty plays a strong role in maximising profits. Since the initial model of~\cite{Newsboy1}, the newsvendor model has been extended in many ways. The extension most relevant to this paper is the multi-period newsvendor model~\citep{arrow1958}. This is the natural extension of the problem to the case where the newsvendor needs to meet demand in multiple time periods, and is able to make separate orders for each.  

Although some papers study static newsvendor models~\citep{MATSUYAMA2006170, CHEN2017, Mehran19}, where the newsvendor must commit to their order quantities prior to the selling period, it is more common in the literature to consider dynamic newsvendor models. In a dynamic newsvendor model, at the start of each period in the horizon, the newsvendor selects their order quantity for that day. This way, they have exact knowledge of the amount of inventory remaining at the time of ordering, as opposed to static models where future inventory levels must be estimated beforehand. Early dynamic models were finite horizon discrete DP models where base-stock policies were optimal. An example of this comes from~\cite{Bouakiz92}, who considered the case where the demand random variables are independent and identically distributed with a known distribution.  A continuous time version of the dynamic model was later solved by~\cite{KOGAN2003448}, who showed it to be equivalent to solving a set of discrete-time problems. Soon after its introduction, papers on the dynamic multi-period model considered more complex situations with respect to demand behaviour and knowledge about its distribution. \cite{Retsef07} developed policies based on only samples from the true demand distribution, with no assumptions being made about the distribution itself. Other extensions include models with partially observable demand~\citep{Benoussan07}, non-stationary demand~\citep{KIM2015139} and service-dependent demand~\citep{Deng14}. These papers illustrate the importance of accurate demand modelling in dynamic newsvendor models, and highlight that it is very common in such problems for demand information to be incomplete. 

Another important extension of the newsvendor problem is the \textit{distribution free} (DF) newsvendor model~\citep{oneperiodrobustsol}. This model represents situations where the true demand distribution is not known exactly, but some of its moments are known exactly. The model then maximises the worst-case profit over the ambiguity set containing all distributions with said moments. Since the work of~\cite{oneperiodrobustsol} for the single-period single-product DF model, the DF concept has received significant attention in the newsvendor literature. Early extensions include models with multiple products and random yield~\citep{gallego_extensions}, balking~\citep{DF_balking}, uncertainty in cost parameters~\citep{OUYANG20021}, shortage penalty costs and budget constraints~\citep{ALFARES2005465}. Later models included additional complexities such as advertising and the costs thereof~\citep{Lee11Advertising}, risk- and ambiguity-aversion~\citep{AverseNV} and carbon emissions~\citep{LiuCarbon}. Due to the fact that these models are not dynamic, they can usually be solved by either KKT conditions or Cauchy-Schwarz bounds on the worst-case cost. Although much less common, the DF concept has also been applied to the multi-period model. For example, \cite{AHMED2007226} studied a DF model arising from using coherent risk measures in the objective function. The model was solved as a finite horizon DP, and it was shown that a base-stock policy was optimal. \cite{LEVINA2010516} considered a model where the only distributional information came from aggregating the opinions of multiple experts. This work was later extended to the case with shortage penalty costs by~\cite{Zhang17}. These authors framed the problem as online learning with expert advice, as opposed to an MDP model.  \cite{Mehran19} found optimal policies for static multi-period distribution free models with moment-based ambiguity sets. As far as we are aware, \cite{AHMED2007226} is the only example of an MDP-based DF newsvendor model.

It is clear from the literature on the DF model that newsvendor models commonly lack distributional information, but many of the MDP models for multi-period newsvendor problems do not account for this. With the recent advancements in RMDPs combined with the fact that many multi-period newsvendor models are formulated as MDPs, this problem is a very appropriate application of our methods. Our research differs from the existing newsvendor literature in two key ways. Where the majority of the DF newsvendor literature focuses on the case where some moments of the demand distribution are known, we do not make any such assumption. Using DF methods usually entails estimating the moments that are assumed to be known, but studies have found that this can lead to various complications. For example, it has been found that this can lead to overly conservative solutions~\citep{Likelihood_RNVP}, suboptimal solutions~\citep{DRNV_Wasserstein} and poor estimates of the true cost function~\citep{ROSSI2014674}. As such, our approach is closer to the more recent papers in RMDPs~\citep{Grand_Clement_Kroer_2021, ho2022robust}, where ambiguity sets contain all distributions that can be considered to be close to some nominal distribution. In addition, unlike these two papers, we consider parametric ambiguity sets. This allows us to model cases where the newsvendor demand is parametric, and enforce that the worst-case distribution lies in the same parametric family as the true demand distribution. As discussed earlier, it is very common to assume that demand distributions are parametric~\citep{Nahmias1994739, Agrawal96, GALLEGO2007, Siegel21}, but DF models do not incorporate this. Our methodology allows parametric distributions to be used, but without the pitfalls of assuming that parameter estimates are truth.

\section{Modelling and algorithms}

In this section, we define our model and present the algorithms used to solve it. The general robust formulation is presented in Section~\ref{sec:general_robust}. The robust value iteration algorithm is presented in Section~\ref{sec:RVI}. Following this, Section~\ref{sec:phi} presents $\phi$-divergence based non-parametric ambiguity sets and Section~\ref{sec:NP_proj} details how the resulting RMDP is solved. We detail these methods since they will act as benchmarks for our parametric methods. In Section~\ref{sec:param_sets} we formulate our parametric ambiguity sets, and in Section~\ref{sec:solve_param} we detail our solution algorithms. 

\subsection{General robust model}\label{sec:general_robust}

The RMDP we consider is formulated as follows. The state and action spaces are defined as $\m{S} = \{1,\dots,S\}$ and $\m{A} = \{1,\dots,A\}$, respectively. Decisions are made at each epoch $t \in \m{T} = \mathbb{N}$. The state at time $t$ is a random variable, denoted by $S_t$. Similarly, we denote by $a_t$ the action taken at time $t$. The reward for selecting action $a \in \m{A}$ when in state $s \in \m{S}$ and transitioning to state $s' \in \m{S}$ is given by $r_{s, a, s'} \in \R_+$. We denote by $\Delta_n$ the probability simplex in $\R^n$: $\Delta_n = \{\bm{P} \in \R^n_+: \sum_{i=1}^n P_i = 1\}$. The distribution of the initial state $S_0$, i.e.\ the state at time $t=0$, is denoted by $\bm{Q} \in \Delta_S$.  The distribution of $S_{t+1}$ given that action $a$ is taken in state $s$ at time $t$ is given by  the unknown distribution $\bm{P}^0_{s, a} = (P^0_{s, a, 1},\dots,P^0_{s, a, S}) \in \Delta_S$. Here, $P^0_{s, a, s'} = \P(S_{t+1} = s'|S_t = s, a_t = a)$ for any $t \ge 0$. Similarly, we write $\bm{P}^0_s$ to denote a matrix where the element on the $a$th row and $s'$th column is $P^0_{s,a,s'}$. Denote by $\Pi = (\Delta_A)^S$ the set of all stationary, randomised policies. A policy $\bm{\pi}$ is a matrix $\bm{\pi} = (\pi_{s, a})_{s \in \m{S}, a \in \m{A}} \in \Pi$ such that $\pi_{s, a}$ gives the probability of taking action $a$ when in state $s$ for each $(s, a) \in \m{S} \times \m{A}$ under policy $\bm{\pi}$. Denote by $\m{P} \subseteq (\Delta_S)^{S \times A}$ an \textit{ambiguity set} for $\bm{P}^0$. Each $\bm{P} \in \m{P}$ and $\bm{\pi} \in \Pi$ induce a stochastic process $\l\{(s_t, a_t\r)\}_{t=0}^{\infty}$ on the space $(\m{S} \times \m{A})^{\infty}$ of sample paths, and $\E_{\bm{P}, \bm{\pi}}$ refers to the expectation w.r.t.\ this process. Then, the robust MDP problem is given by:
\begin{equation}\label{eq:RMDP}
    \max_{\bm{\pi} \in \Pi}\min_{\bm{P} \in \m{P}} \E_{\bm{P}, \bm{\pi}}\l[\sum_{t=0}^{\infty}\gamma^t r_{s_t, a_t, s_{t+1}}\bigg| S_0 \sim \bm{Q} \r],
\end{equation}
where $\gamma \in (0,1)$ is a \textit{discount factor}. We consider $s$-rectangular ambiguity sets, which are of the form:
\begin{equation}
    \m{P} = \m{P}_1 \times \dots \times \m{P}_s, \quad \m{P}_s \subseteq (\Delta_{S})^A \fa s \in \m{S}.
\end{equation}

\subsection{Statewise Bellman equations and robust value iteration}\label{sec:RVI}

Given an initial estimate $v^0_s \fa s \in \m{S}$, robust value iteration is performed by iteratively updating the estimates using the robust Bellman equation~\eqref{eq:robust_bell} for $n = 0,1,\dots$:
\begin{equation}\label{eq:robust_bell}
    v^{n+1}_s = \max_{\bm{\pi}_s \in \Delta_A}\min_{\bm{P}_s \in \m{P}_s} \sum_{a \in \m{A}} \pi_{s, a} \sum_{s' \in \m{S}} P_{s, a, s'} (r_{s, a, s'} + \gamma v^n_{s'}) \fa s \in \m{S}. 
\end{equation}
Adapting the pseudocode by~\cite{powell2007approximate}, this leads to the following robust value iteration algorithm:
\begin{enumerate}
    \item Initialise $n = 0$, $\Delta = 0$, $\bm{v}^0 = \bm{0}$, and select $\varepsilon$.
    \item While $\Delta \ge \frac{\varepsilon \gamma}{1-2\gamma}$:
    \begin{enumerate}
        \item For each $s \in \m{S}$, solve~\eqref{eq:robust_bell} to find the value of $v^{n+1}_s$.
        \item Set $\Delta = ||\bm{v}^{n+1} - \bm{v}^n||$ where $||\bm{v}|| = \max_{s \in \m{S}} v_s$.
        \item Set $n=n+1$.
    \end{enumerate}
    \item Set $\bm{v}^* = \bm{v}^n$ and let the policy that solves~\eqref{eq:robust_bell} under $\bm{v}^n = \bm{v}^*$ be $\bm{\pi}^*$. 
    \item Return $\bm{\pi}^*$ and compute the optimal total reward under $\bm{\pi}^*$ as $\sum_{s\in S} Q_s v^*_s$.
\end{enumerate}
Step 2(a) is referred to as solving a robust Bellman update. 

\subsection{\texorpdfstring{$\phi$-divergence ambiguity sets}{XD}}\label{sec:phi}

The most common ambiguity sets in RMDPs are non-parametric, i.e.\ they do not make use of any information about the family of distributions in which the true distribution lies. Common non-parametric ambiguity sets are distance-based~\citep{Grand_Clement_Kroer_2021, ho2022robust}. Such ambiguity sets contain only distributions that lie within a pre-prescribed maximum distance from a nominal or estimated distribution $\hat{\bm{P}}_s$. In other words, a non-parametric distance-based ambiguity set is of the form given in~\eqref{eq:NP_DB}.
\begin{equation}\label{eq:NP_DB}
    \m{P}_s = \l\{\bm{P}_s \in (\Delta_{A})^S: \sum_{a \in \m{A}} d_a(\bm{P}_{s, a}, \hat{\bm{P}}_{s, a}) \le \kappa\r\} \fa s \in \m{S}.
\end{equation}
Here, $d_a: \Delta_S \times \Delta_S \to \R_+$ is a \textit{distance measure}. We will consider cases where $d_a$ is a \textit{$\phi$-divergence}, i.e.\ it satisfies:
\begin{equation}
    d_a(\bm{P}_{s, a}, \hat{\bm{P}}_{s, a}) =  \sum_{s'=1}^S \hat{P}_{s, a, s'}\phi\l(\frac{P_{s, a, s'}}{\hat{P}_{s, a, s'}}\r),
\end{equation}
where $\phi:\R_+ \to \R_+$ is a \textit{$\phi$-divergence function}. With different choices of $\phi$, the class of $\phi$-divergences encompasses many distances measures, such as the Kullback-Leibler divergence (KLD), $\chi^2$ distance, and Burg entropy. As described by~\cite{DenHertog13}, one benefit of such ambiguity sets is that we can choose $\kappa$ such that $\m{P}_s$ is an approximate confidence set for the true
distribution. Suppose that the true distribution for state $s$, $\bm{P}^0_s$, lies in a parameterised
set $\{\bm{P}^{\bm{\theta}}_s\ | \ \bm{\theta}_s \in \Theta_s\}$, and let the true parameter be $\bm{\theta}^0_s$. We will assume that only $\bm{\theta}_s$ is required to compute $\bm{P}^{\bm{\theta}}_s$ and that $\bm{P}^{\bm{\theta}} = (\bm{P}^{\bm{\theta}}_1,\dots,\bm{P}^{\bm{\theta}}_S)$ is parameterised by $\bm{\theta} = (\bm{\theta}_1,\dots,\bm{\theta}_S)$. Also suppose that the distributions $\bm{P}^0_{s, a}$ are independent. Then, for each $(s, a) \in \m{S} \times \m{A}$, $\bm{P}^0_{s, a}$ is a distribution parameterised by $\bm{\theta}^0_{s, a}$. Suppose that we take $N$ sample transitions from each   
$\bm{P}^0_{s, a}$ and use these to create a maximum likelihood estimate (MLE) $\hat{\bm{\theta}}_{s, a}$ of $\bm{\theta}^0_{s, a}$. Then, if we choose $\kappa$ according to~\eqref{eq:rho}, the set $\m{P}_s$ is an approximate $100(1-\alpha)\%$
confidence set for $\bm{P}^0_{s}$ around $\hat{\bm{P}}_{s} = \bm{P}^{\hat{\bm{\theta}}}_{s}$.
\begin{equation}\label{eq:rho}
    \kappa = \frac{\phi''(1)}{2N} \chi^2_{oA, 1-\alpha}.
\end{equation}
In~\eqref{eq:rho}, $o$ is the dimension of $\Theta_s$ and $\chi^2_{o, 1-\alpha}$
is the $100(1-\alpha)^{\text{th}}$ percentile of the $\chi^2$ distribution with
$o$ degrees of freedom. Note that, while~\cite{DenHertog13} use $o$ degrees of freedom, we use $oA$ in~\eqref{eq:rho} since $\sum_{a \in \m{A}} d_a(\bm{P}_{s, a}, \hat{\bm{P}}_{s, a})$ is the sum of $A$ independent $\chi^2_{o}$ random variables; i.e.\ it is a $\chi^2_{oA}$ random variable.

\subsection{Solving the robust Bellman update}\label{sec:NP_proj}

In this section, we detail the algorithms that we will use to solve the robust Bellman update under non-parametric ambiguity sets, that will act as benchmarks for our methods. In Section~\ref{sec:QP_reform}, we describe how to reformulate the update using the conjugate of the $\phi$-divergence function. In Section~\ref{sec:NP_projection}, we describe the projection-based bisection search algorithm of~\cite{ho2022robust}.

\subsubsection{Solution via reformulation} \label{sec:QP_reform}

The robust Bellman update problem under $\phi$-divergence ambiguity sets can be reformulated using the \textit{convex conjugate} of a $\phi$-divergence function:
\begin{equation}
    \phi^*(z) = \sup_{\tau\ge 0}\{z \tau  - \phi(\tau)\}.
\end{equation}
Using this definition, following the steps given by~\cite{DenHertog13}, we dualise the inner problem of~\eqref{eq:robust_bell} to arrive at the following reformulation (with $\bm{v} = \bm{v}^n$):
\begin{align}\label{eq:bellman_reformulated}
    \max_{\bm{\pi}_s \in \Delta_A, \bm{\nu}, \eta}\l\{\nu - \eta \kappa - \sum_{a \in \m{A}}\sum_{s' \in \m{S}} \eta \hat{P}_{s,a, s'} \phi^*\l(\frac{\nu_a - \pi_{s, a}(r_{s, a, s'} + \gamma v_{s'})}{\eta}\r): \bm{\nu} \in \R^A, \eta \in \R_+\r\}
\end{align}
where $\nu_a \in \R$ is the Lagrange multiplier for the constraint $\sum_{s' \in \m{S}} P_{s, a, s'} = 1$ for each $a \in \m{A}$, $\nu = \sum_{a \in \m{A}} \nu_a$ and $\eta \in \R_+$ is the Lagrange multiplier for the constraint $\sum_{a \in \m{A}} d_a(\bm{P}_{s, a}, \hat{\bm{P}}_{s, a}) \le \kappa$.  For a derivation of this reformulation, see Appendix~\ref{sec:gen_reform}. 

The model requires different approaches for different $\phi$ functions, due to the different forms $\phi^*$ can take. As an example, for the modified $\chi^2$ divergence, this model can be reformulated as the following conic quadratic program:
\begin{equation}\label{eq:QP}
\begin{aligned}
    \max_{\bm{\pi}_s}&\l\{\nu + \eta (A-\kappa) - \frac{1}{4}\sum_{a \in \m{A}}\sum_{s' \in \m{S}} \hat{P}_{s, a, s'} u_{s, a, s'}\r\}\\
    \text{s.t. } & \sqrt{4\zeta_{s,a,s'}^2 + (\eta - u_{s,a,s'})^2} \le (\eta + u_{s, a, s'}) \fa a \in \m{A} \fa s' \in \m{S}\\
    &\zeta_{s, a, s'} \ge 2\eta + \nu_a - \pi_{s, a}(r_{s, a, s'} + \gamma v_{s'}) \fa a \in \m{A} \fa s' \in \m{S}\\
    &\zeta_{s, a, s'} \ge 0 \fa a \in \m{A} \fa s' \in \m{S}\\
    &\sum_{a \in \m{A}} \pi_{s, a} = 1 \\
    &\pi_{s, a} \ge 0 \fa a \in \m{A} \\
    &\eta \ge 0\\
    &\bm{\nu} \in \R^A.
\end{aligned}
\end{equation}
For more details on the derivation of this reformulation, see Appendix~\ref{sec:mchisq_reform}.

\subsubsection{Projection-based bisection search algorithms}\label{sec:NP_projection}

Model~\eqref{eq:bellman_reformulated} can become large when $A$ and $S$ are large, and so it is not always reasonable to solve it in every step of the value iteration algorithm. Hence,  \cite{ho2022robust} presented a fast projection-based algorithm for solving the corresponding robust Bellman update. We define a simplex projection problem as follows:
\begin{equation}\label{eq:projection}
    \mathfrak{P}(\hat{\bm{P}}_{s,a}; \bm{b}, \beta) = \l[\begin{array}{rl}
        \min_{\bm{P}_{s, a}} &d_a(\bm{P}_{s, a}, \hat{\bm{P}}_{s, a}) \\
        \text{s.t. } & \sum_{s' \in \m{S}} b_{s'}{P}_{s, a, s'} \le \beta \\
        & \bm{P}_{s, a} \in \Delta_S
    \end{array}\r].
\end{equation}
Then, the outline of the algorithm presented by~\cite{ho2022robust} is as follows. In each iteration $n$ of the value iteration algorithm, for each $s \in \m{S}$, the Bellman update is solved via bisection search on the value of $v^{n+1}_s$. This is done via the following algorithm, which we will call non-parametric bisection search (NBS):
\begin{enumerate}
    \item Initialise $\epsilon$ and define $\overline{v}^0_s = \bar{R}_s(\bm{v}^n) = \frac{\max_{(a, s') \in \m{A} \times \m{S}} r_{s, a, s'}}{1-\gamma}, \underline{v}^0_s = \max_{a \in \m{A}} \min_{s' \in \m{S}}\l\{r_{s, a, s'} + \gamma v^n_{s'}\r\}$ and $\delta =  \frac{\epsilon \kappa}{2 A + \bar{R}_s(\bm{v}^n) + A \epsilon}$.
    \item For each $i = 0,\dots$:
    \begin{enumerate}
        \item Set $\beta = \frac{\overline{v}^i_s + \underline{v}^i_s}{2}$.
        \item For each $a \in \m{A}$,
        \begin{enumerate}
            \item If $\mathfrak{P}(\hat{\bm{P}}_{s, a}; \bm{r}_{s, a} + \gamma \bm{v}^n, \beta)$ is infeasible, i.e.\ $\min_{s'}\l\{r_{s, a, s'} + \gamma v^n_{s'}\r\} > \beta$, then set $\underline{d}_a$ and $\overline{d}_a$ equal to  $\kappa + 1$. Go to step 2(c).
            \item Otherwise, solve the projection problem  to $\delta$-optimality to obtain parameter action-wise upper and lower bounds $\underline{d}_a, \overline{d}_a$ on its objective value. 
        \end{enumerate}
        \item Use these bounds to update $\underline{v}^i_s$ and $\overline{v}^i_s$:
        \begin{equation}
            (\underline{v}^{i+1}_s, \overline{v}^{i+1}_s) = \begin{cases}
                (\underline{v}^i_s, \beta) &\text{ if } \sum_{a\in\m{A}} \overline{d}_a \le \kappa,\\
                (\beta, \overline{v}^i_s) &\text{ if } \sum_{a \in \m{A}} \underline{d}_a > \kappa
            \end{cases}]
        \end{equation}
        \item $\overline{v}^{i+1}_s - \underline{v}^{i+1}_s < \epsilon$ or $\kappa \in [\sum_{a \in \m{A}} \underline{d}_a, \sum_{a \in \m{A}} \overline{d}_a)$ then go to step 3.
    \end{enumerate}
    \item Return $\beta = \frac{\overline{v}^{i+1}_s + \underline{v}^{i+1}_s}{2}$.
\end{enumerate}
This generates the updated value estimates $v^{n+1}_s$ for $s \in \m{S}$. Now, in step 2(b), the projection problem is also usually solved via bisection search. The logic behind step 2(b)i is that, if the problem is infeasible then there is no $\bm{P}_{s} \in \m{P}_s$ that achieves an objective value of $\beta$. Hence, this scenario should be treated the same as when the problem is feasible and gives $\sum_{a \in \m{A}} \underline{d}_a > \kappa$. Since the actual value of the objective function does not matter as long as this inequality holds, we set it to $\kappa + 1$. For ambiguity sets defined by $\phi$-divergences such as the Kullback-Leibler divergence and $\chi^2$-distance, \cite{ho2022robust} showed how to solve the projection problem efficiently. For the modified $\chi^2$-distance, their method involves first dividing the projection problem into $S+1$ subproblems, and then reformulating each one as a univariate optimisation problem with at most 3 potential solutions that can be found analytically. Solving the subproblem then corresponds to evaluating each of these potential solutions, and choosing the best of those that are feasible. Then, the subproblems' solutions are compared and the best one is selected. Details of this algorithm can be found in Appendix~\ref{sec:sort_alg}.  Following the completion of the value iteration algorithm, a policy must be retrieved. Since the algorithm of~\cite{ho2022robust} does not return a policy, it must be extracted from solving~\eqref{eq:bellman_reformulated}, using $\bm{v} = \bm{v}^*$.

\subsection{Parametric ambiguity sets}\label{sec:param_sets}

We now present our formulation for the RMDP under parametric transition distributions. Suppose that the true transition distribution $\bm{P}^0$ is uniquely defined by the probability mass function (PMF) and/or cumulative distribution function (CDF) of a parametric probability distribution. In this section, we detail how our model allows us to enforce that the worst-case distribution maintains this structure.

\subsubsection{Formulation}\label{sec:P_formulation}

Suppose that, for each $(s, a) \in \m{S} \times \m{A}$, $\bm{P}^0_{s,a}$ is uniquely defined by the distribution of some exogenous random variable $X_{s, a}$ with support $\m{X}_{s, a}$. Let $f_{X_{s, a}}$ and $F_{X_{s, a}}$ be the PMF and CDF of $X_{s, a}$, which are parameterised by the parameter $\bm{\theta}^0_{s, a} = ({\theta}^0_{s, a, 1}, \dots, \theta^0_{s, a, o})$. Assume that the current state $S_t = s$ and action $a_t = a$ are given. We assume that the next state $S_{t+1}$ is specified by some simple, known function $g$ of the exogenous random variable $X_{s, a}$:
$$S_{t+1} = g(X_{s, a}|s, a).$$
In other words, for a given realisation $x$ of $X_{s, a}$, we can compute the next state as $s_{t+1} = g(x|s, a)$. We define the set of all realisations of $X_{s, a}$ that lead to $S_{t+1} = s'$ as:
\begin{equation}
    \m{X}_{s, a}(s') = \l\{x \in \m{X}_{s, a}: g(x|s, a) = s'\r\}.
\end{equation}
Then, the transition matrix corresponding to the parameter $\bm{\theta}^0$ is given by:
\begin{align}
    P^0_{s, a, s'} &= \P(S_{t+1} = s'|S_t = s, a_t = a)\\
    &= \P(g(X_{s, a}|s, a) = s') \\
    &= \sum_{x \in \m{X}_{s, a}(s')} f_{X_{s, a}}(x|\bm{\theta}^0_{s, a}) \fa s' \in \m{S}.
\end{align}
Since $g$ is known, in this case the value of $\bm{P}^0$ is uniquely specified by $\bm{\theta}^0$. Therefore, the only unknown element required to find the true distribution is $\bm{\theta}^0$. Hence, given that the worst-case distribution should maintain the structure of $\bm{P}^0$, we can simply construct ambiguity sets for $\bm{\theta}^0$. More specifically, we consider ambiguity sets of the form:
\begin{equation}
    \Theta_s \subseteq \R^o \fa s \in \m{S}, \quad \Theta =  \Theta_1 \times \dots \times \Theta_S.
\end{equation}
We can then reformulate the RMDP as:
\begin{equation}\label{eq:PRMDP}
    \max_{\bm{\pi} \in \Pi}\min_{\bm{\theta} \in \Theta}
    \E_{\bm{\theta}, \bm{\pi}}\l[\sum_{t=0}^{\infty}\gamma^t r_{s_t, a_t, s_{t+1}}\bigg| S_0 \sim \bm{Q} \r].
\end{equation}
Let $\bm{P}^{\bm{\theta}}$ represent the transition probabilities corresponding to $\bm{\theta}$. Similarly, for any $s \in \m{S}$ and $\bm{\theta}_s \in \Theta_s$, write $\bm{P}^{\bm{\theta}}_s = (P^{\bm{\theta}}_{s, a, s'})_{a \in \m{A}, s' \in \m{S}}$. Note that, although the superscript for $\bm{P}^{\bm{\theta}}_s$ is $\bm{\theta}$, only $\bm{\theta}_s$ is required to compute it and by rectangularity we can obtain $\bm{P}^{\bm{\theta}}$ simply by obtaining $\bm{P}^{\bm{\theta}}_s$ for all $s \in \m{S}$. Similarly, only $\bm{\theta}_{s, a}$ is required to compute $\bm{P}^{\bm{\theta}}_{s, a} =  (P^{\bm{\theta}}_{s, a, s'})_{s' \in \m{S}}$. Now, using the information about $\bm{P}^0$'s structure, we compute $\bm{P}^{\bm{\theta}}_s$ according to:
\begin{equation}
    P^{\bm{\theta}}_{s, a, s'} = \sum_{x \in \m{X}_{s, a}(s')} f_{X_{s, a}}(x|\bm{\theta}_{s, a}) \quad \fa (a, s') \in \m{A} \times \m{S}.
\end{equation}
The robust state-wise Bellman equation can then be written as:
\begin{equation}\label{eq:robust_bell_param}
    v^{n+1}_s = \max_{\bm{\pi}_s \in \Delta_A}\min_{\bm{\theta}_s \in {\Theta}_s} \sum_{a \in \m{A}} \pi_{s, a} \sum_{s' \in \m{S}} P^{\bm{\theta}}_{s, a, s'} (r_{s, a, s'} + \gamma v^n_{s'}) \fa s \in \m{S}.
\end{equation}
As discussed by~\cite{black}, the non-linearities of the PMFs as functions of the parameters mean that above model is not tractable as a mathematical program if the parameters are treated as decision variables. One way to find $v^{n+1}_s$ approximately is to use a discretisation $\Theta'_s$ of the ambiguity set $\Theta_s$. This allows us to reformulate the problem in~\eqref{eq:robust_bell_param} as:
\begin{equation}\label{eq:P_RBU}
    v^{n+1}_s = \max_{\bm{\pi}_s \in \Delta_A} \l\{\vartheta: \vartheta \le \sum_{a \in \m{A}} \pi_{s, a} \sum_{s' \in \m{S}} P^{\bm{\theta}}_{s, a, s'} (r_{s, a, s'} + \gamma v^n_{s'})  \fa \bm{\theta}_s \in \Theta'_s\r\} \fa s \in \m{S}.
\end{equation}
This problem can be solved as an LP with $\lvert\Theta'_s\rvert + 1$ constraints. Due to this, if a fine discretisation of $\Theta'_s$ is used, this model can be very slow to solve. While this is the approach used in parametric DRO prior to this paper, in robust value iteration we only need to compute $v^{n+1}_s$ and not the optimal policy. Hence, in certain cases, no mathematical programming formulation is necessary. We will discuss this in more detail in Section~\ref{sec:param_proj}. However, please note that solving~\eqref{eq:P_RBU} is currently the only way to extract the optimal policy and worst-case probabilities, to the best of our knowledge.

\subsubsection{Confidence sets for the true parameter}\label{sec:confidence_sets}

We assume that we have access to $N$ samples from the true distribution of $\bm{X}$, i.e.\ the distribution that characterises $\bm{P}^0$. This allows us to create an MLE $\hat{\bm{\theta}}$ of the true parameter $\bm{\theta}^0$. In addition, by standard results in maximum likelihood theory~\citep{MillarRussellB2011MLEa} we have:
\begin{equation}
    \l(\hat{\bm{\theta}}_{s, a} - \bm{\theta}^0_{s, a}\r)^T I_{\E}\l(\bm{\theta}^0_{s, a}\r)\l(\hat{\bm{\theta}}_{s, a} - \bm{\theta}^0_{s, a}\r) \sim \chi^2_o
\end{equation}
approximately, for large $N$. Here, $I_{\E}\l(\bm{\theta}^0_{s, a}\r)$ is the expected \textit{Fisher information matrix}, which is defined by~\eqref{eq:fisher}. In~\eqref{eq:fisher}, $\ell$ is the log-likelihood function for the observed data.
\begin{equation}\label{eq:fisher}
    I_{\E}\l(\bm{\theta}_{s, a}\r) = \l(- \E_{X_{s, a}}\l[\frac{\partial}{\partial \theta_{s,a,i}\partial \theta_{s, a, j}}\ell(\bm{\theta_{s, a}})\r]\r)_{i,j=1,\dots,o}.
\end{equation}
By independence of the random variables $X_{s, a}$ for $a \in \m{A}$, we have that:
\begin{equation}
    \sum_{a \in \m{A}}\l(\hat{\bm{\theta}}_{s, a} - \bm{\theta}^0_{s, a}\r)^T I_{\E}\l(\bm{\theta}^0_{s, a}\r)\l(\hat{\bm{\theta}}_{s, a} - \bm{\theta}^0_{s, a}\r) \sim \chi^2_{oA}.
\end{equation}
Since the two are asymptotically equivalent, we can replace $I_{\E}\l(\bm{\theta}^0_{s, a}\r)$ with $I_{\E}\l(\hat{\bm{\theta}}_{s, a}\r)$. Therefore, an approximate $100(1-\alpha)\%$ confidence set for $\theta_s$ is given by:
\begin{equation}
    \Theta^\alpha_s = \l\{\bm{\theta}_s \in \R^A \times \R^o: \sum_{a \in \m{A}}\l(\hat{\bm{\theta}}_{s, a} - \bm{\theta}_{s, a}\r)^T I_{\E}\l(\hat{\bm{\theta}}_{s, a}\r)\l(\hat{\bm{\theta}}_{s, a} - \bm{\theta}_{s, a}\r) \le \chi^2_{oA, 1-\alpha}\r\}.
\end{equation}
In our experiments, we will use $\Theta^\alpha_s$ as an ambiguity set for our parametric model, for each $s \in \m{S}$. We will refer to a discretisation of this set as $(\Theta^{\alpha}_s)'$.

\subsection{Solving the parametric robust Bellman update}\label{sec:solve_param}

It is often cited (e.g., by~\cite{ho2022robust}) that solving an infinite-horizon RMDP efficiently boils down to being able to solve the robust Bellman update efficiently. In our parametric formulation, if we use the LP approximation, then the model that we solve in each iteration for each state $s \in \m{S}$ is the LP~\eqref{eq:P_RBU}, which has $|(\Theta^{\alpha}_s)'| + 1$ constraints. However, depending on the fineness of the discretisation  used to construct $(\Theta^{\alpha}_s)'$, this set can impose thousands of constraints on the model. Hence, \eqref{eq:P_RBU} can be slow to solve. For this reason, we develop two algorithms for solving the robust Bellman update with parametric transition distributions.

\subsubsection{A cutting surface algorithm}\label{sec:CS}

In our previous paper on DRO~\citep{black}, cutting surface (CS) algorithms have performed very well at solving parametric DRO problems that are formulated using discrete ambiguity sets in the same way as~\eqref{eq:P_RBU}. Hence, we now describe the CS algorithm that we will use for the RMDP.
The idea behind the CS algorithm is as follows. Suppose we are at iteration  $n$ of the value iteration algorithm and currently solving for state $s \in \m{S}$. Start with some initial singleton subset $\Theta_s^1 = \{\bm{\theta}_s^{\text{init}}\}$. Solve~\eqref{eq:P_RBU} using $\Theta_s = \Theta_s^1$ to generate a policy $\bm{\pi}_s^1$. Next, solve the \textit{distribution separation problem}~\eqref{eq:DS_VI} with $k=1$ to find the worst-case parameter $\bm{\theta}_s^1$ for the policy $\bm{\pi}_s^1$. Set $\Theta_s^2 = \Theta_s^1 \cup \{\bm{\theta}_s^1\}$ and repeat until stopping criteria are met. 
\begin{equation}\label{eq:DS_VI}
    \min_{\bm{\theta}_s \in (\Theta^{\alpha}_s)'} \sum_{a \in \m{A}} \pi^k_{s, a} \sum_{s' \in \m{S}} P^{\bm{\theta}}_{s, a, s'} (r_{s, a, s'} + \gamma v^n_{s'})
\end{equation}
The appeal of this algorithm is that it only ever solves the approximate robust Bellman update~\eqref{eq:P_RBU} over some small subset $\Theta_s^k$ of $(\Theta^{\alpha}_s)'$, meaning the LP concerned only has $|\Theta_s^k| \ + 1 = k + 1$ constraints at iteration $k$. Typically, in our previous research, we found that this algorithm typically never runs for more than $k=5$ iterations. A formal description of the algorithm for iteration $n$ of the value iteration algorithm for state $s$ is given below.
\begin{enumerate}
    \item Initialise $\Theta_s^1 = \{\bm{\theta}_s^{\text{init}}\}$ for some $\bm{\theta}_s^{\text{init}} \in (\Theta^{\alpha}_s)'$, set $k=1$.
    \item While $k \le k^{\max}$:
    \begin{enumerate}
        \item Solve the LP~\eqref{eq:P_RBU} using $\Theta_s = \Theta_s^k$ to obtain policy $\bm{\pi}_s^k$, which has a worst-case reward of $\tilde{R}^k$ over $\Theta_s^k$.
        \item Evaluate the worst-case rewards:
        \begin{equation}
            R(\bm{\pi}_s^k|\bm{\theta}_s) = \sum_{a \in \m{A}} \pi^k_{s, a} \sum_{s' \in \m{S}} P^{\bm{\theta}}_{s, a, s'} (r_{s, a, s'} + \gamma v^n_{s'}) \fa \bm{\theta}_s \in (\Theta^{\alpha}_s)',
        \end{equation}
        and find $\bm{\theta}_s^k = \argmax_{\bm{\theta}_s \in(\Theta^{\alpha}_s)'} R(\bm{\pi}_s^k|\bm{\theta}_s)$. Set $R^k = R(\bm{\pi}_s^k|\bm{\theta}_s^k)$.
        \item If $\tilde{R}^k \le R^k + \frac{\varepsilon}{2}$ or $\bm{\theta}_s^k \in \Theta_s^k$ then set $k = k^{\max} + 1$.
    \end{enumerate}
    \item Return $\bm{\pi}_s^k$ with worst-case parameter $\bm{\theta}_s^k$ and worst-case reward $R^k$.
\end{enumerate}
For this paper, this algorithm will serve as a method for solving the approximate robust Bellman update~\eqref{eq:P_RBU}. It will therefore be embedded into step 2(a) of the robust value iteration algorithm in Section~\ref{sec:RVI}. 

\subsubsection{A projection-based algorithm for single parameter distributions}\label{sec:param_proj} 

The main algorithms of~\cite{ho2022robust} are based around solving the robust Bellman update using bisection search, Within each iteration of the bisection algorithm, a set of $|\m{A}|$ simplex projection problems~\eqref{eq:projection} are solved to generate the next upper and lower bounds on the value function. The benefit of this is that the projection problem, in the non-parametric case with $\phi$-divergence ambiguity sets, can often be reformulated as a univariate convex optimisation problem. Solving the projection problem $\mathfrak{P}(\hat{\bm{P}}_{s,a}; \bm{b}, \beta)$ corresponds to finding the closest distribution to $\hat{\bm{P}}_{s, a}$ that yields an objective value that is no larger than $\beta$, when action $a$ is taken in state $s$. In the case of distributions where $\bm{P}^0_{s, a}$ is parametrised by only one parameter (such as when $X_{s, a}$ is binomial with a fixed number of trials, or Poisson), the parametric equivalent of this problem can be stated as: 
\begin{equation}\label{eq:projection_param}
    \tilde{\mathfrak{P}}(\hat{\theta}_{s,a}; \bm{b}, \beta) = \l[\begin{array}{rl}
        \min_{\theta_{s, a}} &\l(\hat{\theta}_{s, a} - \theta_{s, a}\r)^2 I_{\E}\l(\hat{\theta}_{s, a}\r)\\
        \text{s.t. } & \sum_{s' \in \m{S}} b_{s'}{P}^{\bm{\theta}}_{s, a, s'} \le \beta \\
        & \theta_{s, a} \in [\theta^{\min}_{s, a}, \theta^{\max}_{s, a}]
    \end{array}\r],
\end{equation}
If $\sum_{s' \in \m{S}} b_{s'}{P}^{\bm{\theta}}_{s, a, s'}\le \beta$ then the model is trivially solved by $\theta_{s, a} = \hat{\theta}_{s, a}$ with an objective value of 0. Therefore, suppose that $\sum_{s' \in \m{S}} b_{s'}{P}^{\bm{\theta}}_{s, a, s'}> \beta$. Without any type of reformulation, the model is a univariate optimisation problem that can be solved via bisection. The only complication in solving this problem via bisection is the constraint $\sum_{s' \in \m{S}} b_{s'}{P}^{\bm{\theta}}_{s, a, s'} \le \beta$. Note that, since $I^{-1}_{\E}(\hat{\theta}_{s, a})$ is the asymptotic variance of the MLE $\hat{\theta}_{s, a}$, we have that $I_{\E}(\hat{\theta}_{s, a}) \ge 0$. Hence, since $I_{\E}(\hat{\theta}_{s, a})$ is constant in $\theta_{s, a}$, the objective of~\eqref{eq:projection_param} is equivalent to:
$$\min_{\theta_{s, a}} |\hat{\theta}_{s,a} - \theta_{s, a}|.$$
Therefore, it is clear that the optimal solution to~\eqref{eq:projection_param} is the closest $\theta_{s, a}$ to $\hat{\theta}_{s, a}$ in terms of absolute value that satisfies $\sum_{s' \in \m{S}} b_{s'}{P}^{\bm{\theta}}_{s, a, s'} \le \beta$. Since $\sum_{s' \in \m{S}} b_{s'}{P}^{\bm{\theta}}_{s, a, s'}> \beta$, the optimal solution must satisfy $\sum_{s' \in \m{S}} b_{s'}{P}^{\bm{\theta}}_{s, a, s'} = \beta$. To see this, observe that any feasible solution with $\sum_{s' \in \m{S}} b_{s'}{P}^{\bm{\theta}}_{s, a, s'} < \beta$ must be further left or right of $\hat{\theta}_{s, a}$ than a solution with $\sum_{s' \in \m{S}} b_{s'}{P}^{\bm{\theta}}_{s, a, s'} = \beta$.
Suppose that the problem is feasible and let $\theta^{\min}_{s, a}$ and $\theta^{\max}_{s, a}$ be global lower and upper bounds on $\theta_{s, a}$. Then, there must be at least one $\theta_{s, a} \in [\theta^{\min}_{s, a}, \theta^{\max}_{s, a}]$ such that  $\sum_{s' \in \m{S}} b_{s'}{P}^{\bm{\theta}}_{s, a, s'} = \beta$. Based on this, we have three potential scenarios as discussed below:
\begin{enumerate}
    \item There exists a root of $\sum_{s' \in \m{S}} b_{s'}{P}^{\bm{\theta}}_{s, a, s'} = \beta$ in $[\theta^{\min}_{s, a}, \hat{\theta}_{s, a}]$. Let $\theta^l_{s, a}$ be the closest root of $\sum_{s' \in \m{S}} b_{s'}{P}^{\bm{\theta}}_{s, a, s'} = \beta$ to $\hat{\theta}_{s, a}$ in the interval $[\theta^{\min}_{s, a}, \hat{\theta}_{s, a}]$. 
    
    \item There exists a root of $\sum_{s' \in \m{S}} b_{s'}{P}^{\bm{\theta}}_{s, a, s'} = \beta$ in $[\hat{\theta}_{s, a}, \theta^{\max}_{s, a}]$. Let $\theta^u_{s, a}$ be the closest root of $\sum_{s' \in \m{S}} b_{s'}{P}^{\bm{\theta}}_{s, a, s'} = \beta$ to $\hat{\theta}_{s, a}$ in the interval $[\hat{\theta}_{s, a}, \theta^{\max}_{s, a}]$.
    
    \item $\theta^l_{s, a}$ and $\theta^u_{s, a}$ both exist as defined above. 
\end{enumerate}
Solving the projection problem then amounts to finding $\theta^l_{s,a}$ and $\theta^u_{s, a}$, and checking which is closest to $\hat{\theta}_{s, a}$. Given this, we solve our projection problem to $\delta$-optimality for a given $s, a$ using the following algorithm:
\begin{enumerate}
    \item Initialise a gap $\tilde{\epsilon}$, the set of root containing intervals as $\rho = \emptyset$, and upper and lower bounds on $\theta_{s, a}$ as $\theta^{\min}_{s, a}, \theta^{\max}_{s, a}$.
    \item Find interval containing closest left root:
    \begin{enumerate}
        \item Initialise $\theta_{s, a} = \hat{\theta}_{s, a}$, $E =  \beta$.
        \item While $E \ge \beta$ and $\theta_{s, a} \neq \theta^{\min}_{s, a}$:
        \begin{enumerate}
            \item Set $\theta_{s,a} = \max\{\theta_{s,a} - \tilde{\epsilon}, \theta^{\min}_{s, a}\}$.
            \item Compute $\bm{P}^{\theta}_{s, a}$ and set $E = \sum_{s' \in \m{S}} b_{s'}{P}^{\bm{\theta}}_{s, a, s'}$.
        \end{enumerate}
        \item If $E \le \beta$ then set $\rho = \rho \cup \{[\theta_{s, a}, \theta_{s, a} + \tilde{\epsilon}]\}$.
    \end{enumerate}
    \item Find interval containing closest right root:
    \begin{enumerate}
        \item Initialise $\theta_{s, a} = \hat{\theta}_{s, a}$, $E =  \beta$.
        \item While $E \ge \beta$ and $\theta_{s, a} \neq \theta^{\max}_{s, a}$:
        \begin{enumerate}
            \item Set $\theta_{s,a} = \min\l\{\theta_{s,a} + \tilde{\epsilon}, \theta^{\max}_{s, a}\r\}$.
            \item Compute $\bm{P}^{\theta}_{s, a}$ and set $E = \sum_{s' \in \m{S}} b_{s'}{P}^{\bm{\theta}}_{s, a, s'}$.
        \end{enumerate}
        \item If $E \le \beta$ then set $\rho = \rho \cup \{[\theta_{s, a} - \tilde{\epsilon}, \theta_{s, a}]\}$.
    \end{enumerate}
    \item Carry out a bisection search in interval in $\rho$ to find the roots $\theta^l_{s, a}$ and $\theta^u_{s, a}$, stopping once the difference between the upper and lower bounds on the objective function in the bisection interval is no larger than $\delta$. Store the intervals $\l[\underline{\theta}^x_{s, a}, \overline{\theta}^x_{s, a}\r]$ for $x \in \{l, u\}$.
    \item Return the interval $\l[\underline{\theta}^*_{s, a}, \overline{\theta}^*_{s, a}\r]$ whose midpoint is closest to $\hat{\theta}_{s, a}$ in terms of absolute value.
\end{enumerate}
We use an iterative procedure starting from $\hat{\theta}_{s, a}$ in steps 2 and 3 in order to reduce the number of times we need to compute $\bm{P}^{\theta}_{s, a}$. Since we are only interested in the closest roots to $\hat{\theta}_{s, a}$, there is no need to enumerate all intervals of width $\tilde{\epsilon}$. Note that, in some cases, $\theta_{s,a}$ may not have both a global lower and upper bound. For example, if $\theta_{s, a}$ is a Poisson parameter, then it has no upper bound. However, if the solution to the projection problem does not lie in the ambiguity set, then it is treated the same as if the problem is infeasible. Hence, we are only interested in roots inside the ambiguity set and so in such cases we can use the bounds from the ambiguity set. If $\theta_{s, a}$ is a binomial parameter then we can pick either $\theta^{\min}_{s, a}, \theta^{\max}_{s, a}$ to be $0, 1$ or $\min((\Theta^{\alpha}_s)'), \max((\Theta^{\alpha}_s)')$. Since we will typically split the interval $[\theta^{\min}_{s, a}, \theta^{\max}_{s, a}]$ into an equal number of sub-intervals and hence each choice results in the same amount of computation, which upper and lower bounds we pick are not of particular importance.

Given the above, we adapt the non-parametric bisection search algorithm from Section~\ref{sec:NP_proj} into the following parametric algorithm, which we call parametric bisection search (PBS):
\begin{enumerate}
    \item Initialise $\epsilon$ and define $\overline{v}^0_s = \bar{R}_s(\bm{v}^n) = \frac{\max_{(a, s') \in \m{A} \times \m{S}} r_{s, a, s'}}{1-\gamma}, \underline{v}^0_s = \max_{a \in \m{A}} \min_{s' \in \m{S}}\l\{r_{s, a, s'} + \gamma v^n_{s'}\r\}$ and $\delta =  \frac{\epsilon \kappa}{2 A + \bar{R}_s(\bm{v}^n) + A \epsilon}$.
    \item For each $i = 0,\dots$:
    \begin{enumerate}
        \item Set $\beta = \frac{\overline{v}^i_s + \underline{v}^i_s}{2}$.
        \item For each $a \in \m{A}$, 
        \begin{enumerate}
            \item If $\bar{\mathfrak{P}}(\hat{\theta}_{s, a}; \bm{r}_{s, a} + \gamma \bm{v}^n, \beta)$ is infeasible, i.e.\ $\min_{s'}\l\{r_{s, a, s'} + \gamma v^n_{s'}\r\} > \beta$, then set $\underline{c}_a$ and $\overline{c}_a$ equal to  $\chi^2_{oA, 1-\alpha} + 1$. Go to step 2(c).
            \item Otherwise, solve the projection problem  to $\delta$-optimality to obtain parameter action-wise upper and lower bounds $\underline{c}_a, \overline{c}_a$ on its objective value. If projection algorithm returns no solutions, set both to $\chi^2_{oA, 1-\alpha} + 1$.
        \end{enumerate}
        \item Use these bounds to update $\underline{v}^i_s$ and $\overline{v}^i_s$:
        \begin{equation}
            (\underline{v}^{i+1}_s, \overline{v}^{i+1}_s) = \begin{cases}
                (\underline{v}^i_s, \beta) &\text{ if } \sum_{a\in\m{A}} \overline{c}_a \le \chi^2_{oA, 1-\alpha},\\
                (\beta, \overline{v}^i_s) &\text{ if }\sum_{a\in\m{A}} \underline{c}_a > \chi^2_{oA, 1-\alpha}
            \end{cases}
        \end{equation}
        \item $\overline{v}^{i+1}_s - \underline{v}^{i+1}_s < \epsilon$ or $\chi^2_{oA, 1-\alpha} \in [\sum_{a\in\m{A}} \underline{c}_a, \sum_{a\in\m{A}} \overline{c}_a)$ then go to step 3.
    \end{enumerate}
    \item Return $\beta = \frac{\overline{v}^{i}_s + \underline{v}^{i}_s}{2}$. 
\end{enumerate}
Given this algorithm, we can efficiently carry out value iteration without ever needing a solver. However, after this is complete, the optimal policy will need to be retrieved by solving the approximate MIP reformulation~\eqref{eq:P_RBU} of the robust Bellman update. This can be done using the cutting surface algorithm of Section~\ref{sec:CS}.



\section{A capacitated dynamic multi-period newsvendor problem}

As an example problem, we consider a dynamic multi-period newsvendor problem. This version of the problem has discrete demands and actions, and a capacity limiting the amount that can be held in inventory for any given period. In Section~\ref{sec:NV_model}, we describe the model in detail. Then, in Section~\ref{sec:binom}, we formulate the model under binomial demands and perform computational experiments to test our algorithms in this case. Finally, in Section~\ref{sec:poisson}, we formulate the model and test our algorithms under Poisson demands.

\subsection{Model}\label{sec:NV_model}

Suppose that $S_t$ represents the amount of inventory in a system of some product affected by uncertain demand. Let, $a_t$ be the amount of this product to order at the start of period $t$, to be sold during period $t$. Products are delivered immediately. We assume that there is a capacity $C$ for holding stock in inventory, so that $\m{S} = \{0,\dots,C\}$.  Given that action $a$ is taken in state $s$, if $s + a > C$ then any excess product is lost as it cannot be stored. Although the newsvendor could technically order infinite stock, they have no reason to. Hence, $\m{A} = \{0, \dots, C\}$, and so $S = \lvert\m{S}\rvert = C + 1$ and $A = \lvert\m{A}\rvert = C+ 1$.  We assume that every unit of stock that must be held for a period incurs a holding cost of $h$, and if the newsvendor runs out of stock then they pay a stockout cost of $b'$. Furthermore, assume that one unit of stock sells for $c$ and is purchased for $w < c$. Let the demand for the product, $X_{s, a}$, be a random variable whose distribution is parameterised by the unknown parameter $\theta^0_{s, a}$ for each $(s, a) \in \m{S} \times \m{A}$. Then, given $S_t = s$ and $a_t = a$, we have:
\begin{equation}
    S_{t+1} = \max\{0, \min\{s + a, C\} - X_{s, a}\}.
\end{equation}
For shorthand, let $\bar{s} = \min\{s + a, C\}$ be the post-action pre-demand state. 
Then, we have that $g(x|s, a) = \max\{0, \bar{s} - x\}$, and therefore:
\begin{equation}
    \m{X}_{s, a}(s') = \begin{cases}
        \l\{\bar{s}, \bar{s}  + 1,\dots, C-1, C\r\} &\text{ if } s' = 0\\
        \bar{s} - s' &\text{ if } s' > 0.
    \end{cases}
\end{equation}
Therefore, the transition distribution satisfies:
\begin{align}
    P^0_{s, a, s'} &= \begin{cases}
        \sum_{x = \bar{s}}^{C} f_{X_{s, a}}(x | \bm{\theta}^0_{s, a}) & \text{ if } s' = 0,\\
        f_{X_{s, a}}(\bar{s} - s' | \bm{\theta}^0_{s, a})& \text{ if } s' > 0.
    \end{cases}\\
    &= \begin{cases}
        1 - \sum_{x = 0}^{\bar{s} - 1} f_{X_{s, a}}(x| \bm{\theta}^0_{s, a}) & \text{ if } s' = 0,\\
        f_{X_{s, a}}(\bar{s} - s'| \bm{\theta}^0_{s, a})& \text{ if } s' > 0.
    \end{cases}
\end{align}
The reward for taking action $a$ in state $s$ and moving to state $s'$ is given by the following. Define the event of a stockout as $\mathds{1}\{s'=0\}$. Then, the rewards are given by:
\begin{equation}\label{eq:NV_rewards}
    r_{s, a, s'} = c\max\{\bar{s} - s', 0\} - wa - h(\bar{s} - \max\{\bar{s} - s', 0\}) - b' \mathds{1}\{s'=0\}.
\end{equation}
The term $b'\mathds{1}\{s' = 0\}$ will charge the newsvendor a flat cost of $b'$ whenever they miss demand. It is more common in the newsvendor literature to incur a backorder cost for every unit of missed demand, representing the newsvendor paying an additional cost to meet this demand after initially not meeting it. This would involve adding cost of $b'\max\{X_{s, a} - \bar{s}, 0\}$ instead of $b'\mathds{1}\{s' = 0\}$. Hence, the rewards would depend on $X_{s, a}$ and we need to formulate the robust Bellman update~\eqref{eq:robust_bell} in a different fashion. The main change would be that the distribution of $X_{s, a}$ would be required to calculate the expected one-stage rewards as opposed to simply the transition matrix. Hence, we would replace the inner minimisation over $\bm{P}_s$ with a minimisation over candidates $\bm{P}'_s \in \m{P}'_s$ for the true distribution of demand $\bm{X}_s$. We would then replace the inner expected value with respect to the next state with an expectation with respect to $X_{s, a}$. Since each $\bm{P}'_{s, a}$ has dimension $|\m{X}_{s, a}|$, this would remain tractable in the non-parametric case for finite support demand random variables. However, since we do not know any moments of the distribution of $X_{s, a}$, it would result in an infinite number of decision variables for infinite support demand random variables. This would not affect the parametric model, however, which would still find the worst-case parameter directly. For more details on the reformulations in the case of a backorder cost, see Appendix~\ref{sec:X_reformulation}.  

Considering a stockout cost instead of a backorder cost means that the robust Bellman update can be computed via an expectation over the finite set $\m{S}$, regardless of whether or not $\m{X}_{s, a}$ is finite. The downside of this formulation is that it can penalise the newsvendor for meeting demand exactly. However, if this is a concern then one can set $w + h > b'$ to ensure that the newsvendor would still prefer to meet demand exactly and pay a stockout cost than to purchase too much stock and hold one item for the following period. Also, it is important to note that shortage costs are implicitly represented in this model via missed profits and the newsvendor can see how much demand was lost after the period is complete.

\subsection{Numerical experiments with binomial demands}\label{sec:binom}

To examine the efficacy of the algorithms described in this paper, we now carry out numerical experiments on the dynamic newsvendor problem. Firstly, in Section~\ref{sec:binom_AS}, we describe the binomial ambiguity sets used. Then, in Section~\ref{sec:exp_des}, we describe the parameters used. Following this, in Section~\ref{sec:times}, we discuss the times taken by each algorithm to finish value iteration and compute the optimal policy. Finally, in Section~\ref{sec:P_vs_NP}, we compare the parametric and non-parametric value functions and resulting policies.

\subsubsection{Ambiguity sets}\label{sec:binom_AS}

Suppose that $X_{s, a} \sim \text{Bin}(C, p^0_{s, a})$ for $(s, a) \in \m{S} \times \m{A}$, and hence:
\begin{equation}
    f_{X_{s, a}}(x|p^0_{s, a}) = \binom{C}{x}(p^0_{s, a})^x(1-p^0_{s, a})^{C-x} \quad (x \in \{0,\dots,C\}).
\end{equation}
We set the number of trials as $C$ for the following reasons. Since a binomial random variable is bounded above by the number of trials, binomial demands might correspond to a scenario where a restriction is placed on demand by the newsvendor. In this case, the number of trials represents the maximum demand allowed by the newsvendor before no more orders are allowed. The number of trials provides a way for the newsvendor to limit the amount of unmet demand that is possible. Since any demand above $C$ is guaranteed to be unmet regardless of the current stock levels, it is not reasonable for the number of trials to be set above $C$. This would not have any benefit for the newsvendor or the customers.  Another logical choice for the number of trials may be $\min\{s+a, C\}$. However, this would imply that the newsvendor would need to update the upper bound on demand after every order, and this information would need to be conveyed to customers. In addition, it suggests that the newsvendor is always able to meet demand exactly, which is not a realistic modelling assumption. Also, the newsvendor would be unable to observe how much demand was lost or if the demand met the capacity, which is inconvenient for improving their decision making and capacity levels. When the maximum demand is $C$, the newsvendor can infer whether or not more capacity is required from how often a demand of $C$ occurs. Similarly, they can decide if they have too much capacity if, for example, the demand is always less than the capacity.

Since the number of trials is fixed, the distribution of $\bm{X}$ is uniquely parameterised by $\bm{p}^0 = (p^0_{s, a})_{(s, a) \in \m{S} \times \m{A}}$. In the notation of Section~\ref{sec:confidence_sets}, this means that $o=1$. Suppose that we take a sample $\bm{x}_{s, a} = (x^1_{s, a}, \dots, x^N_{s, a})$ from the distribution of $X_{s, a}$ for each $(s, a) \in \m{S} \times \m{A}$. Then, the MLE $\hat{p}_{s, a}$ of $p^0_{s, a}$ is given by:
\begin{equation}
    \hat{p}_{s, a} = \frac{\sum_{j=1}^N x^j_{s, a}}{NC} \fa (s, a) \in \m{S} \times \m{A}.
\end{equation}
In addition, the Fisher information is given by:
\begin{equation}
    I_{\E}(\hat{p}_{s,a}) = \frac{NC}{\hat{p}_{s, a}(1-\hat{p}_{s, a})}.
\end{equation}
Therefore, our approximate $100(1-\alpha)\%$ confidence set for $\bm{p}^0_s$ is given by:
\begin{equation}\label{eq:binom_conf}
    \Theta^{\alpha}_s = \l\{\bm{p}_s \in [0, 1]^A: \sum_{a \in \m{A}} \frac{NC(p_{s, a} - \hat{p}_{s, a})^2}{\hat{p}_{s, a}(1-\hat{p}_{s, a})} \le \chi^2_{A, 1-\alpha} \r\}
\end{equation}
As discussed in Section~\ref{sec:P_formulation}, in order for the parametric robust Bellman update~\eqref{eq:robust_bell_param} to be tractable, we consider discrete ambiguity sets. Since~\eqref{eq:binom_conf} is a multivariate set, it is difficult to discretise directly. Therefore, we will construct a set $\Theta^{\text{base}}_s$ such that $\Theta^{\alpha}_s \subseteq \Theta^{\text{base}}_s$ and discretise $\Theta^{\text{base}}_s$ instead. Then, we construct a discretisation of $\Theta^{\alpha}_s$ by extracting all elements of $\Theta^{\text{base}}_s$ that also lie in $\Theta^{\alpha}_s$. Observe that the definition of $\Theta^{\alpha}_s$ implies that every $\bm{p}_s \in \Theta^{\alpha}_s$ satisfies:
\begin{equation}
    p_{s, a} \in p^{\text{I}}_{s, a} = \l[\max\l\{0, \hat{p}_{s, a} - \sqrt{\frac{\chi^2_{A, 1 - \alpha} \hat{p}_{s, a}(1-\hat{p}_{s, a})}{NC}}\r\}, \min\l\{1, \hat{p}_{s, a} + \sqrt{\frac{\chi^2_{A, 1 - \alpha} \hat{p}_{s, a}(1-\hat{p}_{s, a})}{NC}}\r\}\r]
\end{equation}
for all $a \in \m{A}$. Therefore, defining:
\begin{equation}
    \Theta^{\text{base}}_s = p^{\text{I}}_{s, 1} \times \dots \times p^{\text{I}}_{s, A}, 
\end{equation}
we have $\Theta^{\alpha}_s \subseteq \Theta^{\text{base}}_s$. Furthermore, define $p^l_{s, a}$ and $p^u_{s, a}$ as the lower and upper bounds of $p^{\text{I}}_{s, a}$ for each $(s, a) \in \m{S} \times \m{A}$. We can then find discretisations of each $p^{\text{I}}_{s, a}$ containing $M$ points as follows:
\begin{equation}
    \tilde{p}^{\text{I}}_{s, a} = \l\{p^l_{s, a} + m \frac{p^u_{s, a} - p^l_{s, a}}{M - 1}\r\}.
\end{equation}
Then, a discretisation of $\Theta^{\text{base}}_s$ is given by $(\Theta^{\text{base}}_s)' = \tilde{p}^{\text{I}}_{s, 1} \times \dots \times \tilde{p}^{\text{I}}_{s, A}$. Finally, a discretisation of $\Theta^{\alpha}_s$ is given by $(\Theta^{\alpha}_s)' = (\Theta^{\text{base}}_s)' \cap \Theta^{\alpha}_s$. 

\subsubsection{Experimental design}\label{sec:exp_des}

We now detail the experiments used to test our algorithms on the dynamic newsvendor problem. The parameters used were as follows. We considered $w, h, b', c \in \{1, 5, 10\}$ such that $w > c$. The capacities we considered we $C \in \{1, 2, 3, 7, 9, 14\}$. This leads to $\lvert \m{S} \rvert =  \lvert \m{A} \rvert \in \{2, 3, 4, 8, 10, 15\}$. We used a discount parameter of $\gamma = 0.5$ in all cases. For each algorithm, the value iteration algorithm was run for a maximum of $n^{\max} = 1000$ iterations. With regard to ambiguity sets, we always used $\alpha = 0.05$, the discretisation parameter was $M \in \{3, 5, 10\}$ and we took $N \in \{10, 50\}$ samples to create the MLEs. Each algorithm was given a maximum time of 4 hours to complete value iteration and find the optimal policy after value iteration ended. In addition, the parametric algorithms were given a maximum of 4 hours to complete their precomputation, i.e.\ computing the discrete ambiguity set and corresponding transition probabilities. Note that this is not required for solving value iteration with PBS, but it is required to compute the optimal policy after value iteration ends. If an algorithm ran for 4 hours and the model was not solved, then the algorithm is said to have  \textit{timed out} for this instance. Both the parametric and non-parametric models used $100(1-\alpha)\%$ confidence sets as ambiguity sets. The parametric model used~\eqref{eq:binom_conf} or a discretisation thereof, and the non-parametric model used~\eqref{eq:NP_DB} where $\kappa$ is defined by~\eqref{eq:rho}. In addition, we used a value iteration tolerance of $\varepsilon = 10^{-6}$ and we initialised the value functions as $\bm{v}^0 = \bm{0}$. In PBS, we used a gap of $\tilde{\epsilon} = 0.01$ with $(\theta^{\min}_{s, a}, \theta^{\max}_{s, a}) = (0, 1)$ for all $(s, a) \in \m{S} \times \m{A}$. Finally, the bisection search tolerance used for PBS and NBS was $\epsilon = 10^{-7}$.

The above inputs generated 810 instances. We ran value iteration on each instance using 5 different algorithms, where each one is defined by how it solves each robust Bellman update. The algorithms and how they solve the update are as follows:
\begin{enumerate}
    \item PBS: the parametric projection-based bisection search algorithm of Section~\ref{sec:param_proj}.
    \item CS: the cutting surface algorithm of Section~\ref{sec:CS}.
    \item LP: using Gurobi to solve the approximate LP reformulation~\eqref{eq:P_RBU} of the parametric update.
    \item QP: using Gurobi to solve the CQP reformulation~\eqref{eq:QP} of the non-parametric update~\eqref{eq:NP_DB}.
    \item NBS: the non-parametric projection-based bisection search algorithm of Section~\ref{sec:NP_projection}.
\end{enumerate}
Value iteration was run until either $n^{\max}$ iterations had been completed, 4 hours of run time had been used, or the algorithm converged. After value iteration ended, for LP, CS and QP the policy was returned. For PBS, the policy was extracted using CS. For NBS, the policy was extracted using QP. 

\subsubsection{Times taken}\label{sec:times}

In this section, we summarise the times taken by the algorithms. We first present the number of times that each algorithm timed out. Firstly, LP and CS timed out while running value iteration 56 and 54 times respectively. No other algorithm timed out while running value iteration. Secondly, although PBS never timed out while running value iteration, it timed out twice while computing the optimal policy. As we will show, PBS is a fast algorithm in itself, and these two timeouts are a result of the slowness of CS in instances with large ambiguity sets.  

Due to the above result, we present the times taken to run value iteration separately from the times taken to compute the policy. Table~\ref{tab:VI_times} summarises the amount of time that each algorithm spent running value iteration. This table shows that PBS took 31 seconds on average to finish value iteration, while CS took 17 minutes 30 seconds and LP took 26 minutes. It is therefore clear that PBS results in greatly reduced times to complete value iteration compared to these solver-based algorithms. CS also saves approximately 12 minutes per iteration compared with LP on average. Note that LP and CS's average times per iteration are large because, when they timed out, they usually timed out after only one iteration. In addition, NBS took 43 seconds on average to complete value iteration, which is 33\% slower than PBS. On average, NBS is much faster than its solver-based equivalent QP, which took over 6 minutes on average to finish value iteration. However, it is important to note that QP led to convergence issues in our experiments. While all other algorithms always finished value iteration in around 31 iterations, when using QP value iteration failed to converge in 378 instances. This was likely due to Gurobi being unable to provide precise enough optimal objective values. In addition, QP was also the fastest algorithm per iteration, and was only slow overall due to value iteration's failure to converge when using this algorithm.  
\begin{table}[htbp!]
    \centering
    \begin{tabular}{llll}
\toprule
Algorithm &   Mean Time &    Max Time & Mean Time Per Iteration \\
\midrule
      PBS &  0:00:31.09 &  0:04:29.37 &              0:00:01.07 \\
       CS &  0:17:30.90 &     4:00:00 &              0:04:25.77 \\
       LP &  0:26:00.73 &     4:00:00 &              0:16:20.66 \\
       QP &  0:06:16.24 &  0:55:59.19 &              0:00:00.41 \\
      NBS &  0:00:43.52 &  0:09:21.02 &              0:00:01.48 \\
\bottomrule
\end{tabular}

    \caption{Summary of times taken to run value iteration (binomial)}
    \label{tab:VI_times}
\end{table}

\begin{figure}[htbp!]
    \begin{subfigure}[t]{0.45\textwidth}
        \centering
        \includegraphics[width=0.9\textwidth]{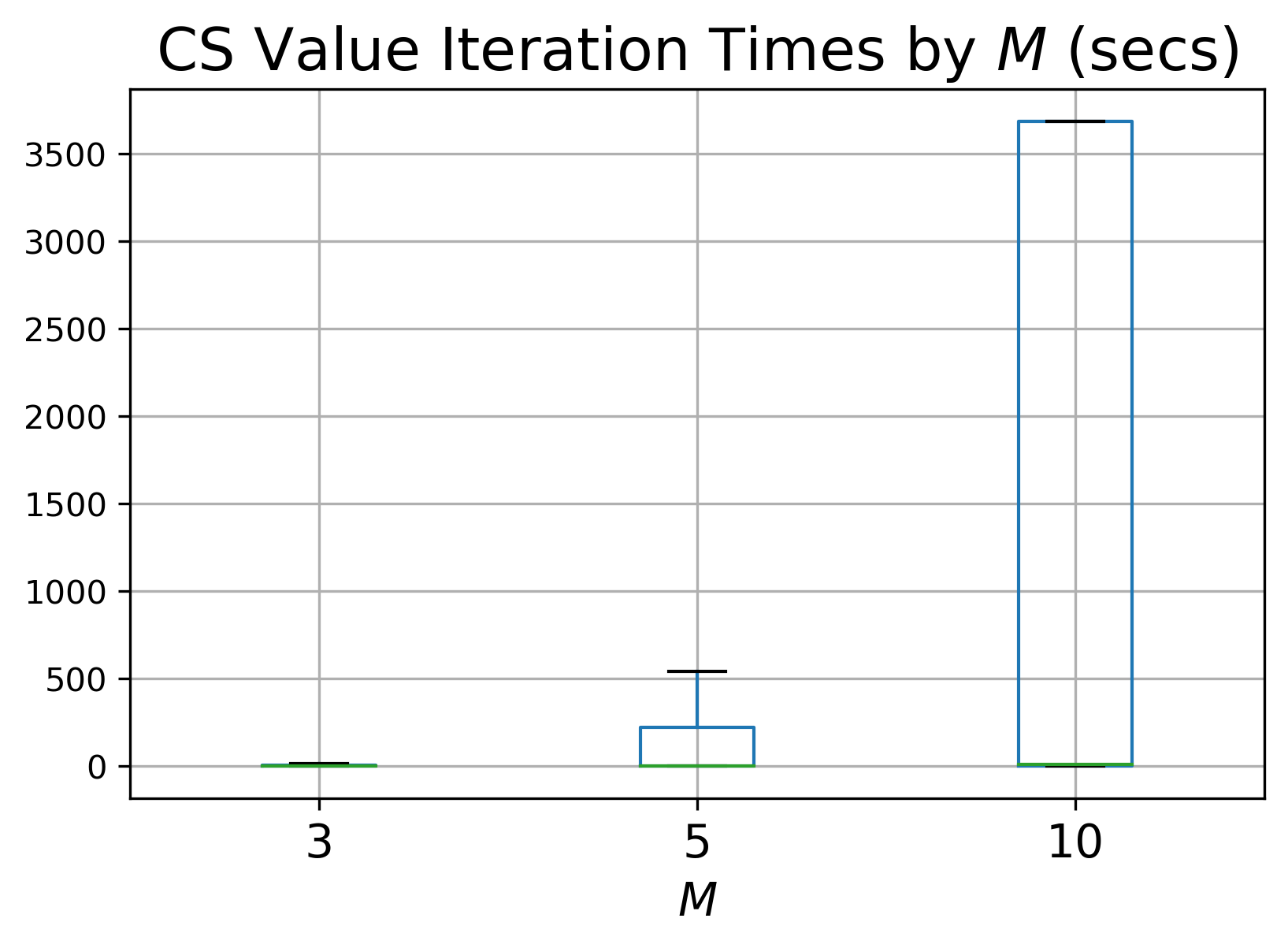}
        \caption{}
        \label{fig:CS_by_M}
    \end{subfigure}
    \begin{subfigure}[t]{0.45\textwidth}
        \centering
        \includegraphics[width=0.90\textwidth]{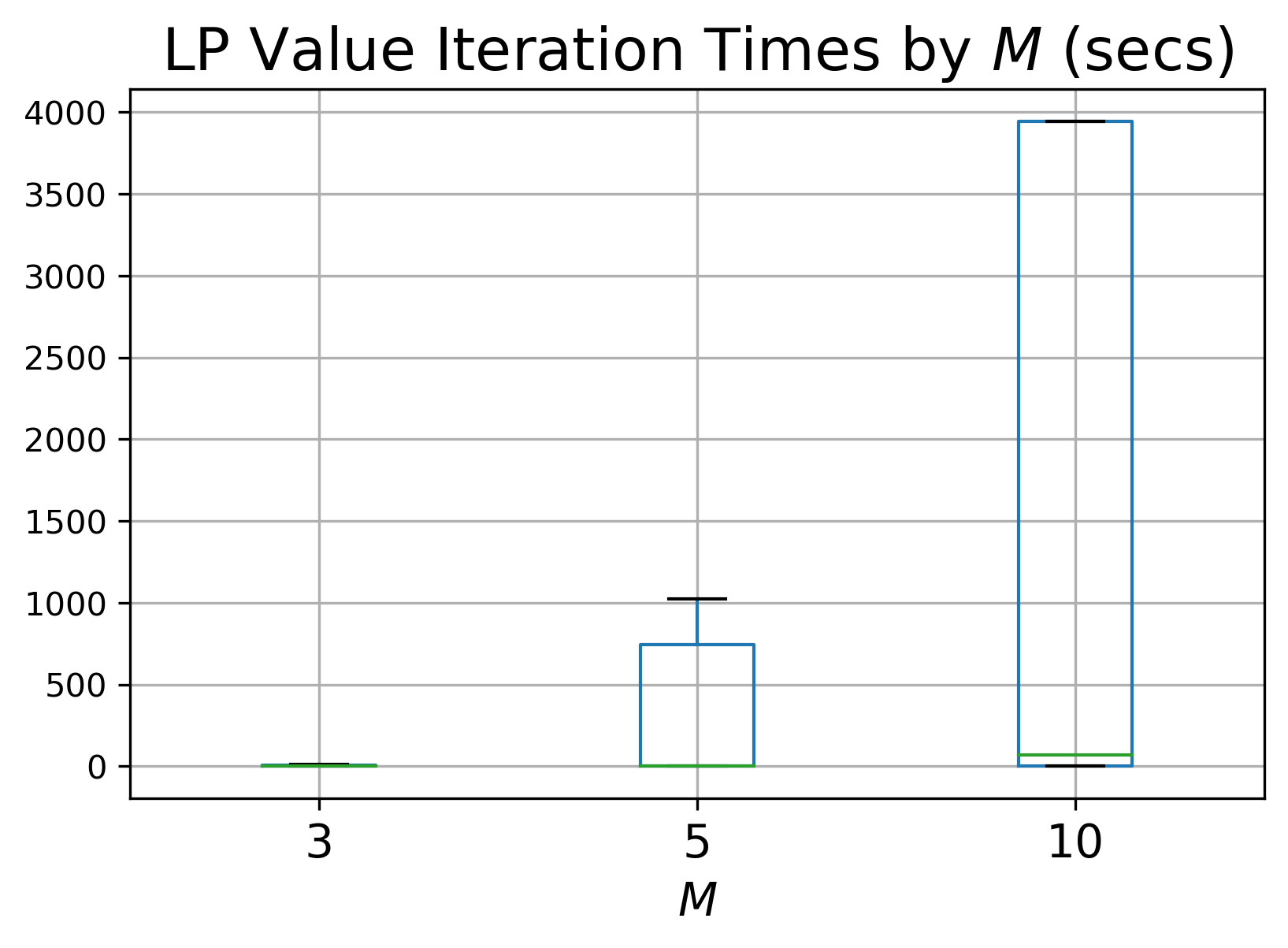}
        \caption{}
        \label{fig:LP_by_M}
    \end{subfigure}
    \caption{Boxplots of value iteration run times of (a) CS and  (b) LP, by $M$.}
    \label{fig:CS_LP_M}
\end{figure}

It is clear from this table that CS and LP can both become very slow. The main reason for this is $M$, the parameter defining the fineness of the discretisation of $\Theta^{\alpha}_s$ used by the parametric solver-based algorithms. We confirm this with Figure~\ref{fig:CS_LP_M}, which shows boxplots of CS and LP's value iteration run times by $M$. Figures~\ref{fig:CS_by_M} and~\ref{fig:LP_by_M} show that both CS and LP scale poorly with $M$ in terms of value iteration run times. However, the effect of $M$ is not particularly noticeable until $M=10$, as was the case for PBS's policy times. CS scales better than LP, but it still becomes slow for instances with large $M$ or large $C$. However, since PBS does not return a policy, an algorithm like CS is required to generate the optimal policy. Please note that, unlike the CS algorithm of~\cite{black}, this CS algorithm is the optimal version which finds the worst-case parameter over the \textit{entire} ambiguity set in every iteration. This explains why it does not offer the same level of time savings when compared with LP as the CS algorithm of~\cite{black}.

\begin{figure}[htbp!]
    \begin{subfigure}[t]{0.45\textwidth}
        \centering
        \includegraphics[width=0.9\textwidth]{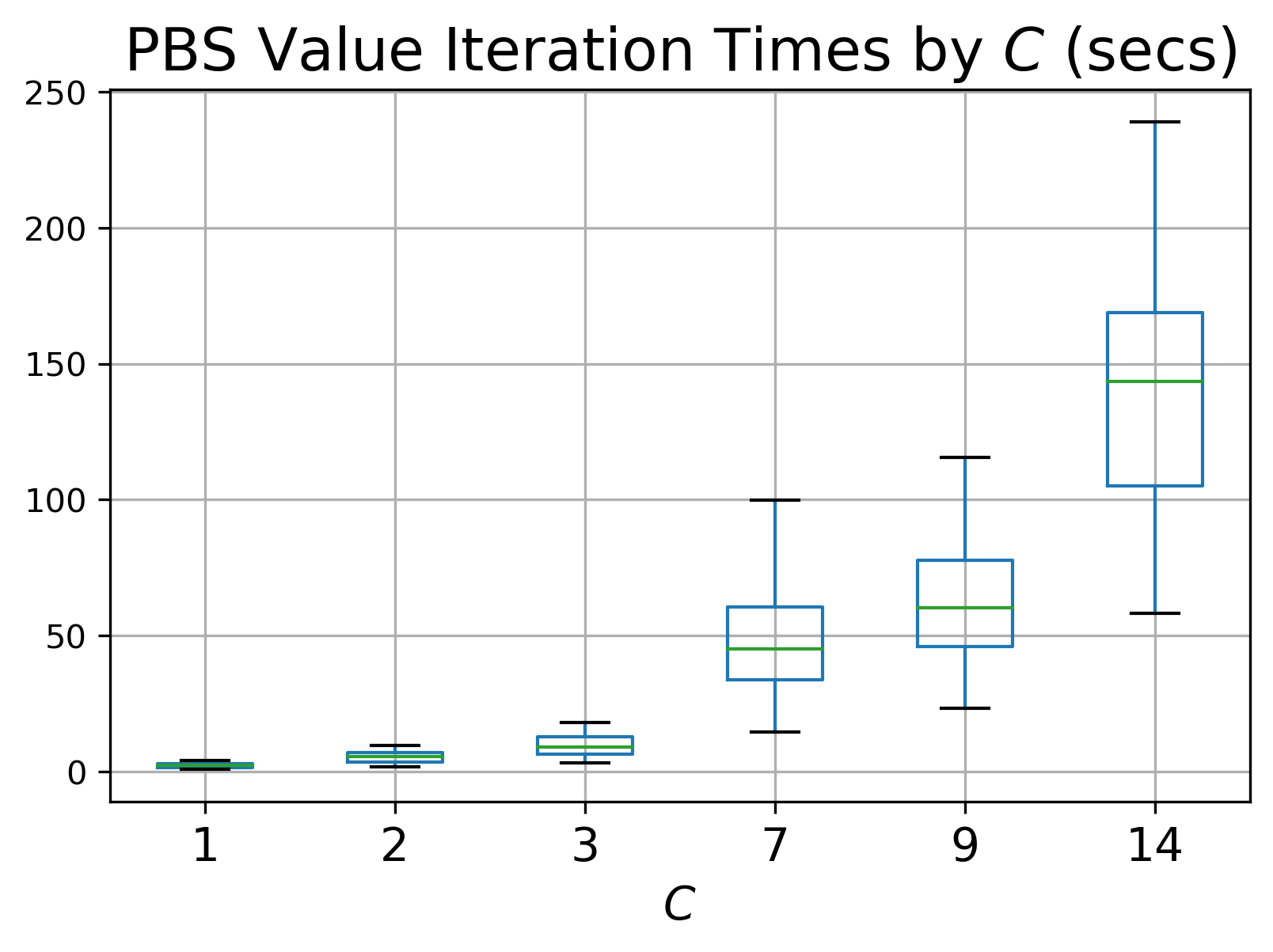}
        \caption{}
        \label{fig:PBS_by_C}
    \end{subfigure}
    \begin{subfigure}[t]{0.45\textwidth}
        \centering
        \includegraphics[width=0.90\textwidth]{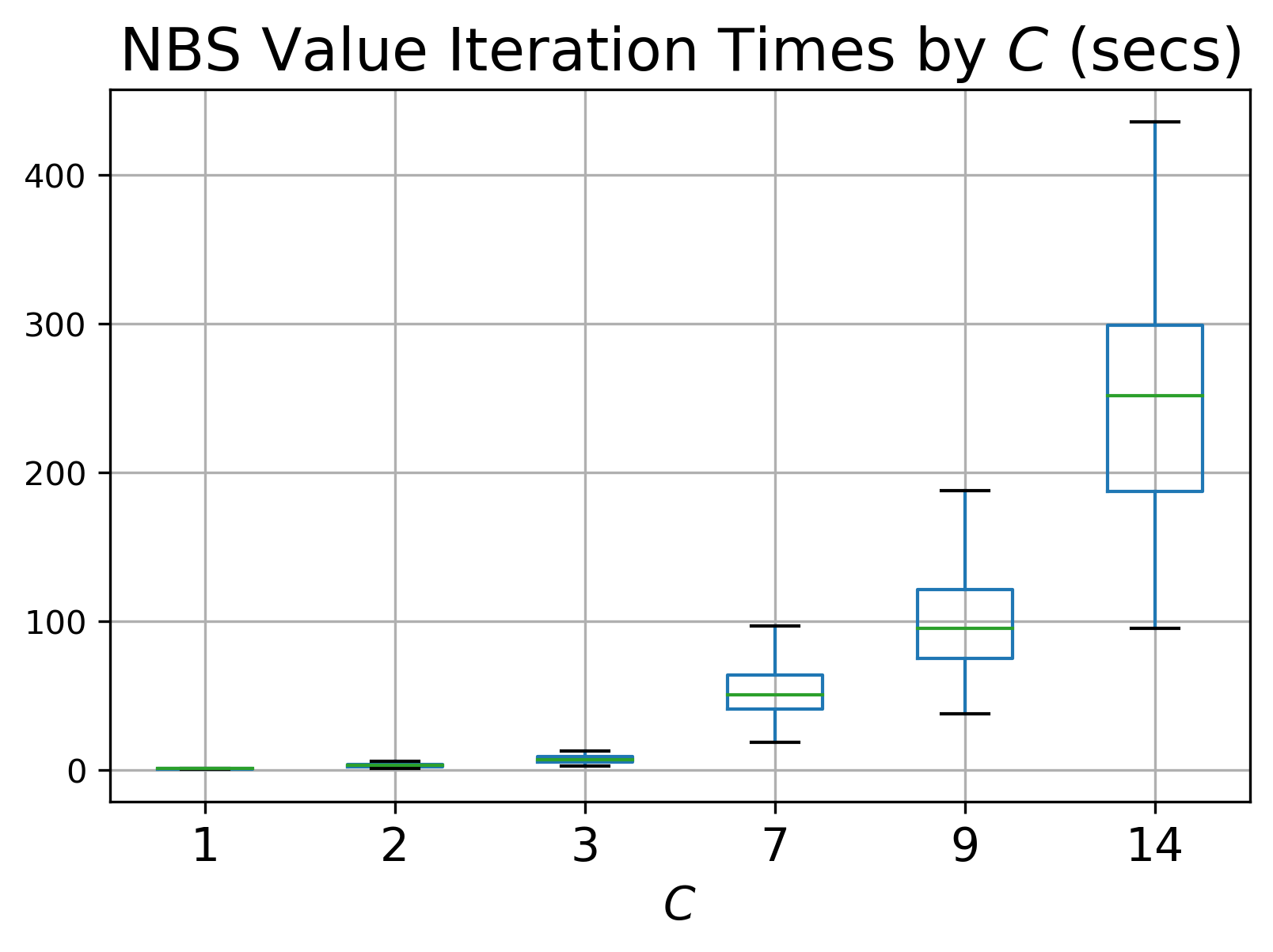}
        \caption{}
        \label{fig:NBS_by_C}
    \end{subfigure}
    \caption{Boxplots of value iteration run times of (a) PBS and  (b) NBS, by $C$ (binomial).}
    \label{fig:times_C}
\end{figure}

Since PBS and NBS do not rely on a discrete ambiguity set, their value iteration times are not affected by $M$. Therefore, the main parameter affecting their value iteration times is $C$. We present boxplots of PBS and NBS's value iteration times by $C$ in Figure~\ref{fig:times_C}. Figures~\ref{fig:NBS_by_C} and~\ref{fig:NBS_by_C} show similar increases in times as $C$ increases, but it is clear that PBS generally scaled better with $C$ than NBS. For $C=14$, PBS typically took no longer than 250 seconds to complete value iteration, while NBS typically took no longer than 450 seconds. On the other hand, NBS was slightly faster than PBS for small $C$. The reason for the difference in scaling is likely because NBS solves $S + 1 = C + 2$ sub-problems in order to solve a projection problem, whereas PBS always carries out a 3-step procedure to solve its projection problems. Hence, the number of steps involved in solving a projection problem increases with $C$ for NBS, but not for PBS.

\begin{figure}[htbp!]
    \begin{subfigure}[t]{0.45\textwidth}
        \centering
        \includegraphics[width=0.9\textwidth]{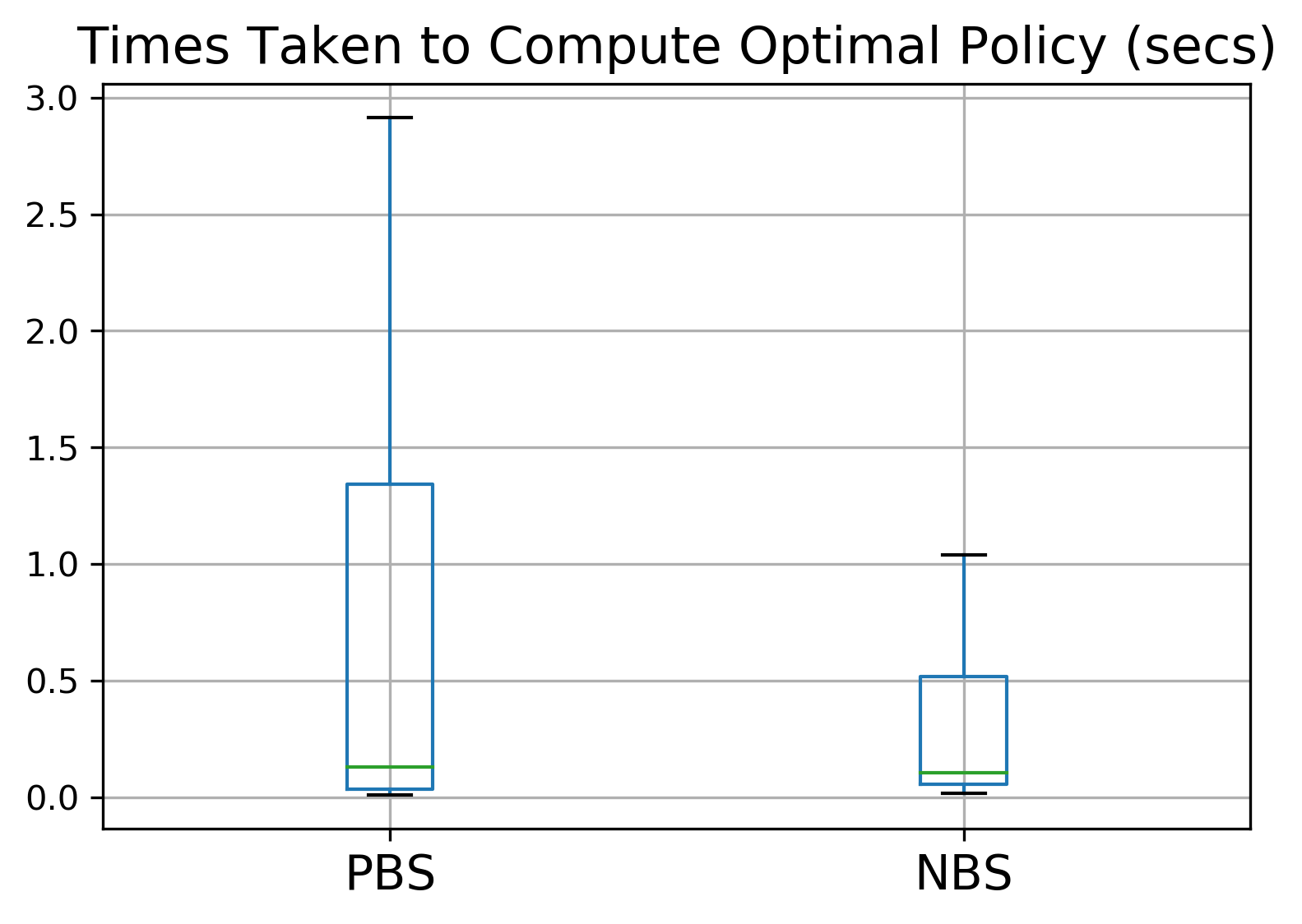}
        \caption{}
        \label{fig:policy_times}
    \end{subfigure}
    \begin{subfigure}[t]{0.45\textwidth}
        \centering
        \includegraphics[width=0.90\textwidth]{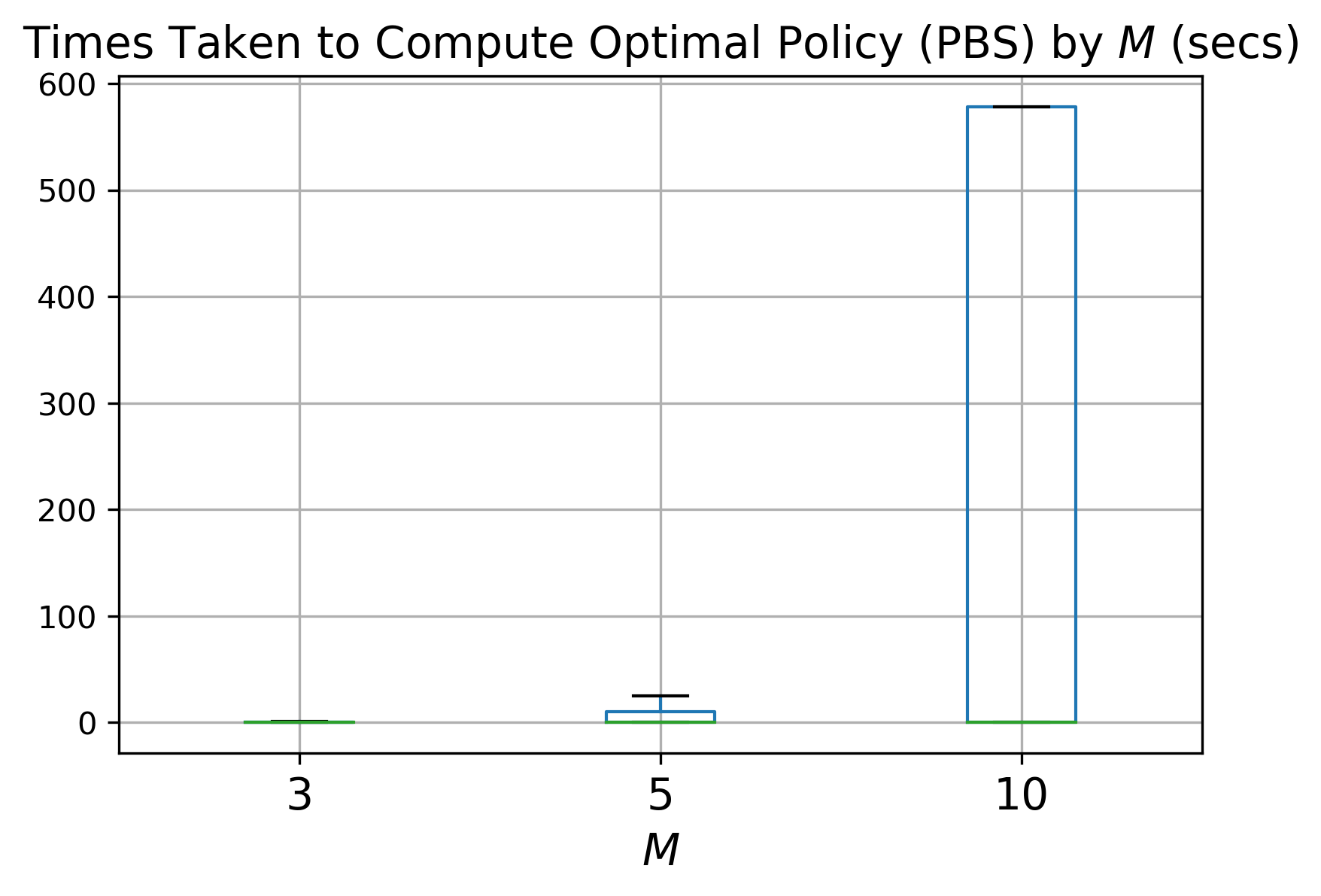}
        \caption{}
        \label{fig:policy_times_M}
    \end{subfigure}
    \caption{Boxplots of times taken to compute optimal policy for (a) both BS algorithms and (b) PBS by $M$ (binomial).}
    \label{fig:policy_times_bin}
\end{figure}

Although PBS was faster than NBS in running value iteration, it generally took longer for the parametric optimal policy to be computed in the parametric case than the non-parametric. This is because computing the policy is not part of PBS, and so CS had to be used for this. On average, it took 8 minutes 27 seconds for CS to compute the policy for PBS's value function, and only 0.43 for QP to compute NBS's. However, the parametric average is greatly affected by a small selection of very slow instances. Figure~\ref{fig:policy_times} shows a boxplot of the times taken to compute the policy for PBS and NBS's value functions after value iteration ended. Note that this boxplot does not show outliers, which are defined as any data that are further than 1.5 times the interquartile range above the 75th percentile or below the 25th percentile.  Figure~\ref{fig:policy_times_bin} shows that while the average time to compute the policy for PBS's values was 8 minutes 27 seconds, the median time was only 0.128 seconds. The 25th and 75th percentiles of the time taken to compute the parametric policy are 0.0325 and 1.34 seconds. 

Figure~\ref{fig:policy_times_M} explains why the average time to compute the policy for PBS was large. In particular, it shows that this starts to take a very long time when $M=10$. This parameter defines the fineness of the discretisation of $\Theta^{\alpha}_s$ used by CS when computing the policy. Since PBS uses $\Theta^{\alpha}_s$ directly in value iteration, this parameter does not affect its value iteration run times. Hence, the slow times to compute a policy for PBS are a reflection on CS's scaling with respect to $M$, not PBS's. This can be confirmed by comparing Figure~\ref{fig:policy_times_M} with Figure~\ref{fig:CS_by_M}, and observing that the exact same pattern is present in both. Due to this, it is a downside that PBS does not provide an optimal policy. The same applies to NBS, but since QP is faster than CS in generating a policy, the effect is not so severe. Comparing the times taken to compute optimal policies is therefore not comparing the run times of PBS and NBS, but comparing the run times of CS and QP.

\subsubsection{Comparison of value functions, distributions and policies}\label{sec:P_vs_NP}

In this section, we compare the values, policies and worst-case transition distributions from the 5 algorithms tested. We first discuss the effect of the discretisation of $\Theta^{\alpha}_s$ on the value functions from LP and CS. Following this, we compare the outputs from PBS and NBS in order to assess the benefits of incorporating additional distributional information into the model. 

\begin{figure}[htbp!]
    \begin{subfigure}[t]{0.45\textwidth}
        \centering
        \includegraphics[width=0.9\textwidth]{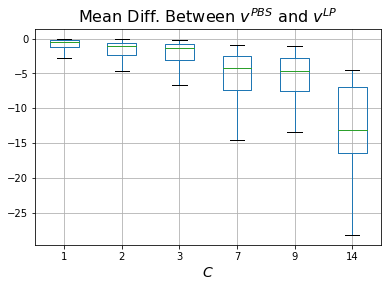}
        \caption{}
        \label{fig:LP_PBS_C}
    \end{subfigure}
    \begin{subfigure}[t]{0.45\textwidth}
        \centering
        \includegraphics[width=0.90\textwidth]{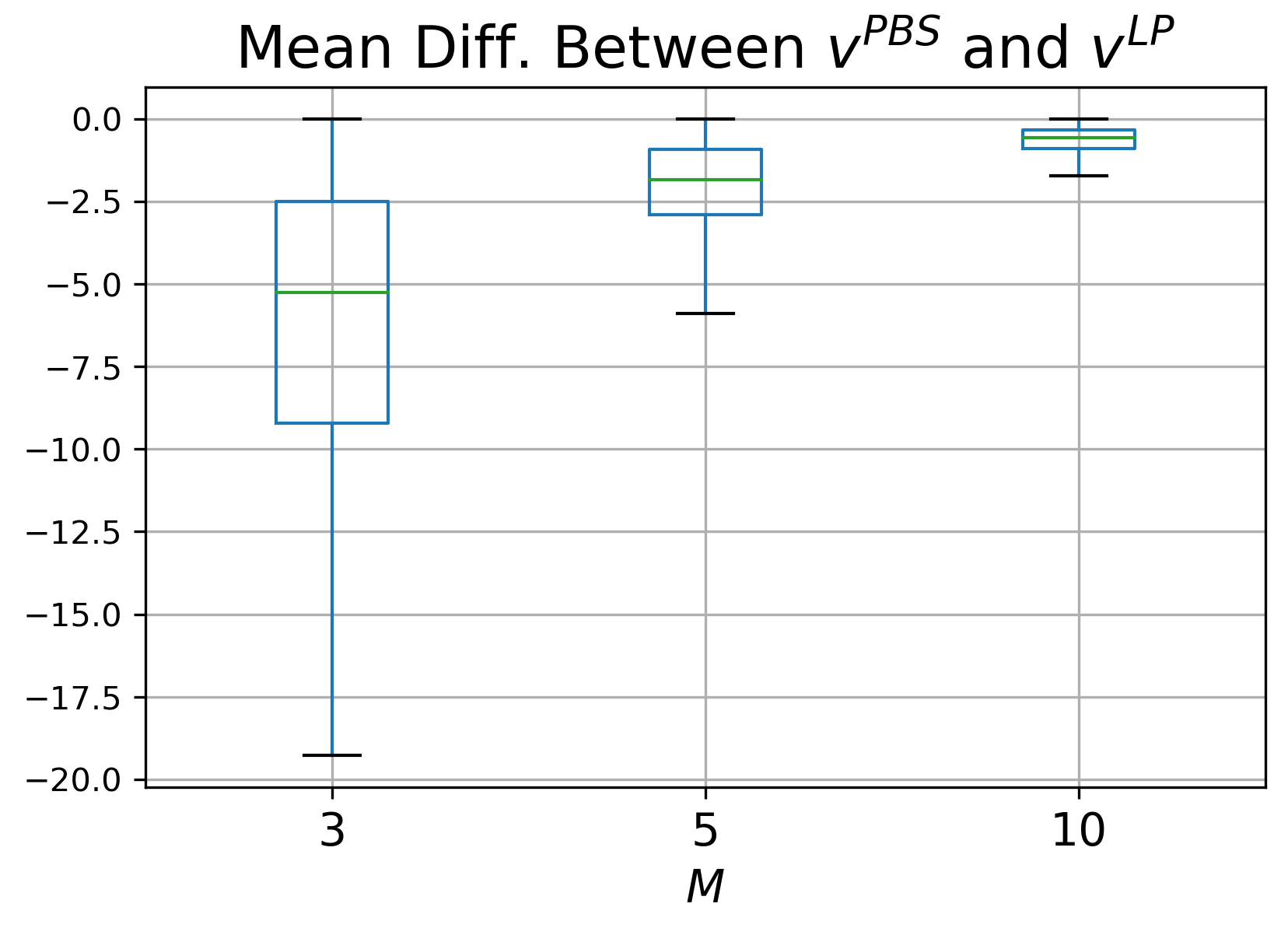}
        \caption{}
        \label{fig:LP_PBS_mean}
    \end{subfigure}
    \caption{Boxplots of mean difference between $\bm{v}^{\text{LP}}$ and $\bm{v}^{\text{PBS}}$ by (a) $C$ and (b) $M$ (binomial).}
    \label{fig:v_param}
\end{figure}

Let $\bm{v}^y$ be the value function generated by running value iteration with algorithm $y \in \m{Y} =\{\text{PBS}, \text{CS}, \text{LP}, \text{QP}, \text{NBS}\}$. Similarly, define $\bm{\pi}^y$ as the policy and $\bm{P}^y$ as the worst-case transition distribution from using algorithm $y$ in the value iteration algorithm. Then, we can summarise the differences between LP's approximate value functions and PBS's optimal value functions via Figure~\ref{fig:v_param}. Figures~\ref{fig:LP_PBS_C} and~\ref{fig:LP_PBS_mean} show boxplots of the mean difference between $\bm{v}^{\text{LP}}$ and $\bm{v}^{\text{PBS}}$ over all instances where LP did not time out, by $C$ and $M$, respectively. The mean differences are calculated as:
\begin{equation}\label{eq:v_diffs}
     \frac{1}{\lvert S \rvert}\sum_{s \in \m{S}}\l(v^{\text{PBS}}_s - v^{\text{LP}}_s\r).
\end{equation}
We see from Figure~\ref{fig:LP_PBS_C} that the average difference between $\bm{v}^{\text{PBS}}$ and $\bm{v}^{\text{LP}}$ was always negative, with the magnitude of the difference growing larger as $C$ increases. This means that LP's value function estimates are higher than PBS's optimal values. This is intuitive, since LP uses a discrete subset of $\Theta^{\alpha}_s$ and therefore cannot find the true worst-case parameter for any given $\bm{\pi}$, in general. Therefore, LP's value function estimates overestimate the worst-case reward for a given policy. As is reflected in Figure~\ref{fig:LP_PBS_mean}, the two value functions get closer as $M$ increases, with average differences that are less than 2.5 in absolute value for $M = 10$. These plots suggest two results. Firstly, the effect of the discretisation increases as $C$ increases. In other words, for larger $C$, the value function estimates resulting from the discrete approximations are less accurate. Secondly, the value function estimates resulting from the discretised ambiguity set appear to converge to their optimal values over the complete (not discretised) ambiguity set.

\begin{figure}[htbp!]
        \centering
        \includegraphics[width=0.5\textwidth]{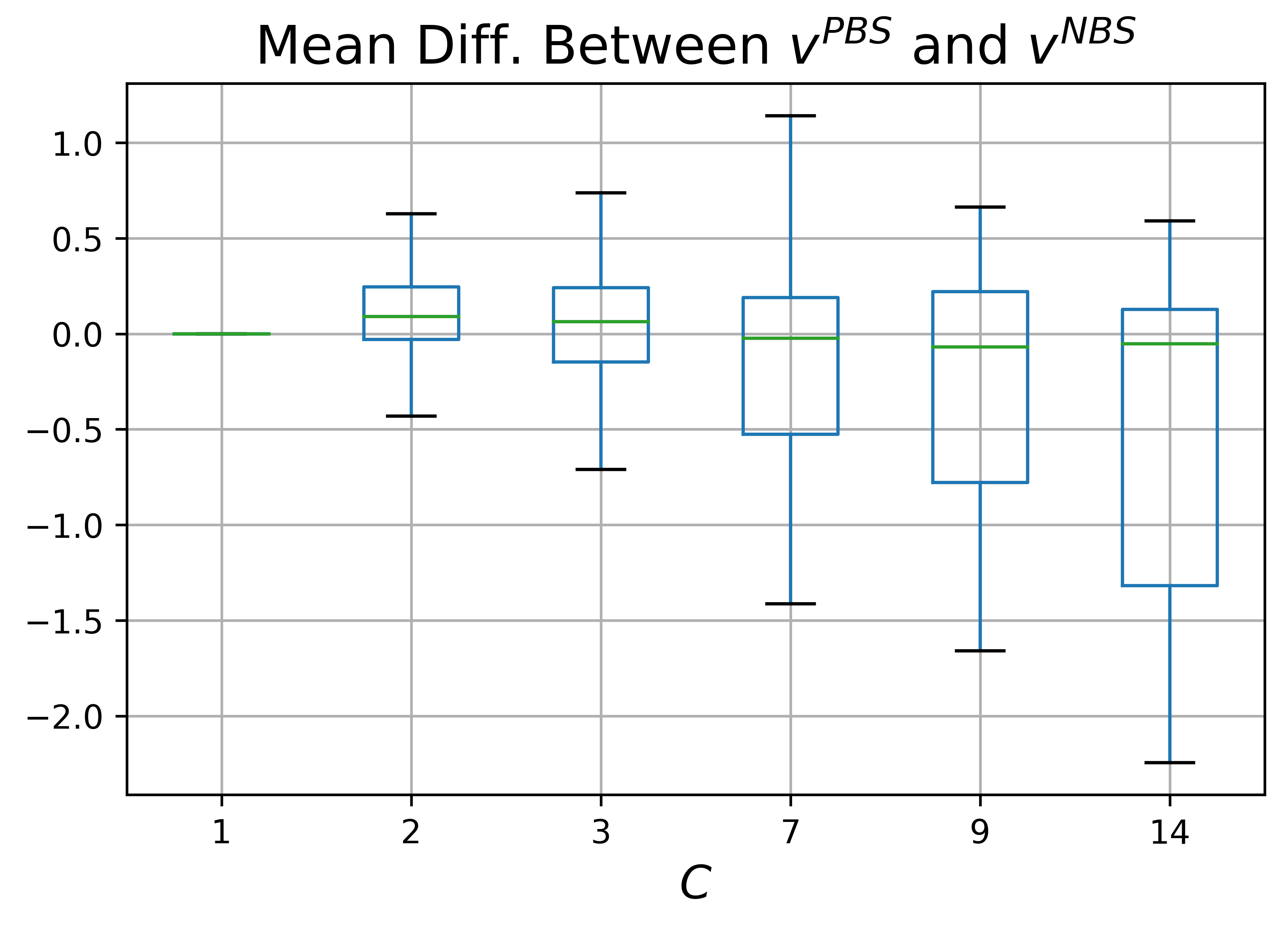}
    \caption{Boxplot of mean difference between $\bm{v}^{\text{NBS}}$ and $\bm{v}^{\text{PBS}}$ by $C$ (binomial).}
    \label{fig:v_diffs_NBS}
\end{figure}

\begin{figure}[htbp!]
    \begin{subfigure}[t]{0.45\textwidth}
        \centering
        \includegraphics[width=0.9\textwidth]{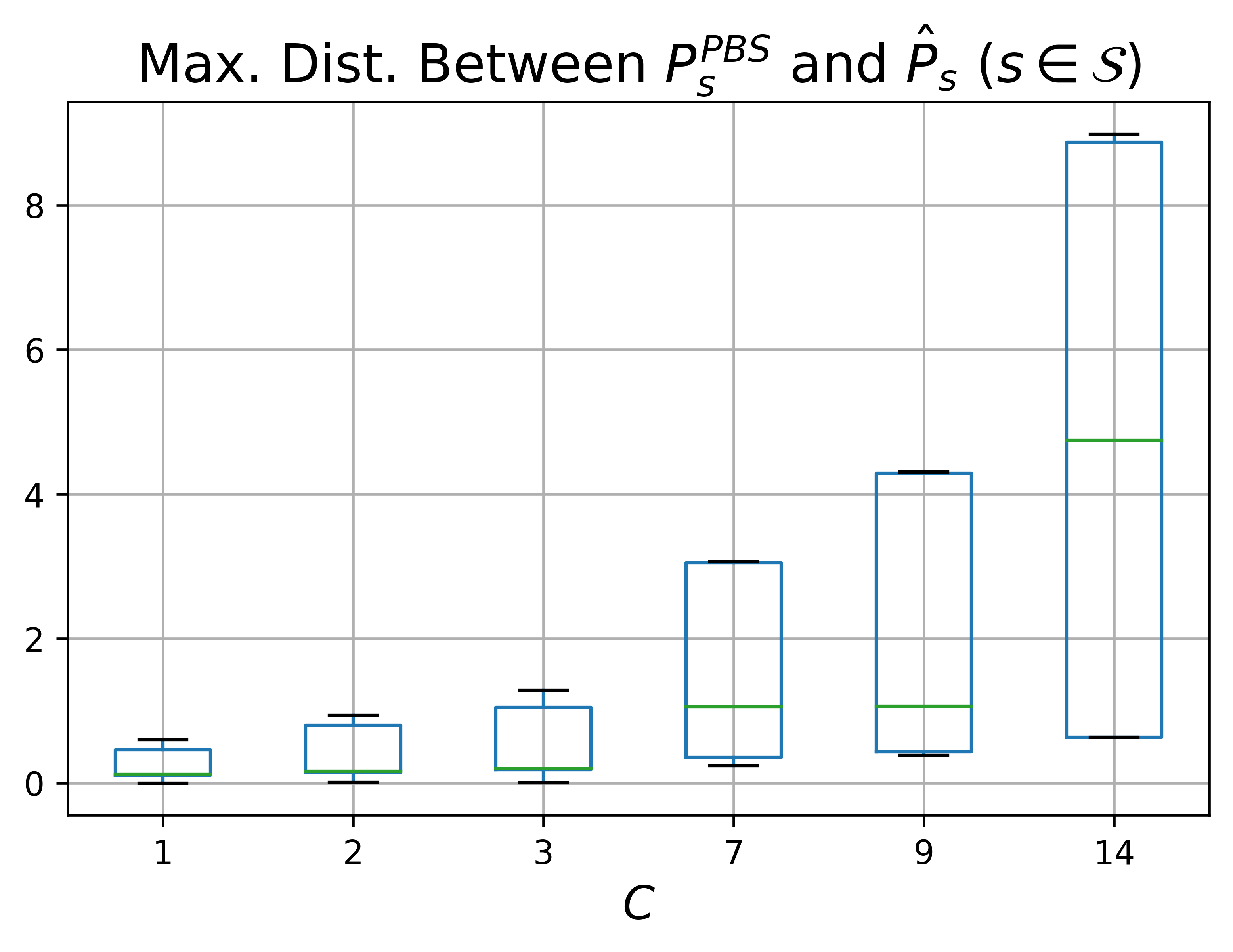}
        \caption{}
        \label{fig:P_vs_Phat}
    \end{subfigure}
    \begin{subfigure}[t]{0.45\textwidth}
        \centering
        \includegraphics[width=0.90\textwidth]{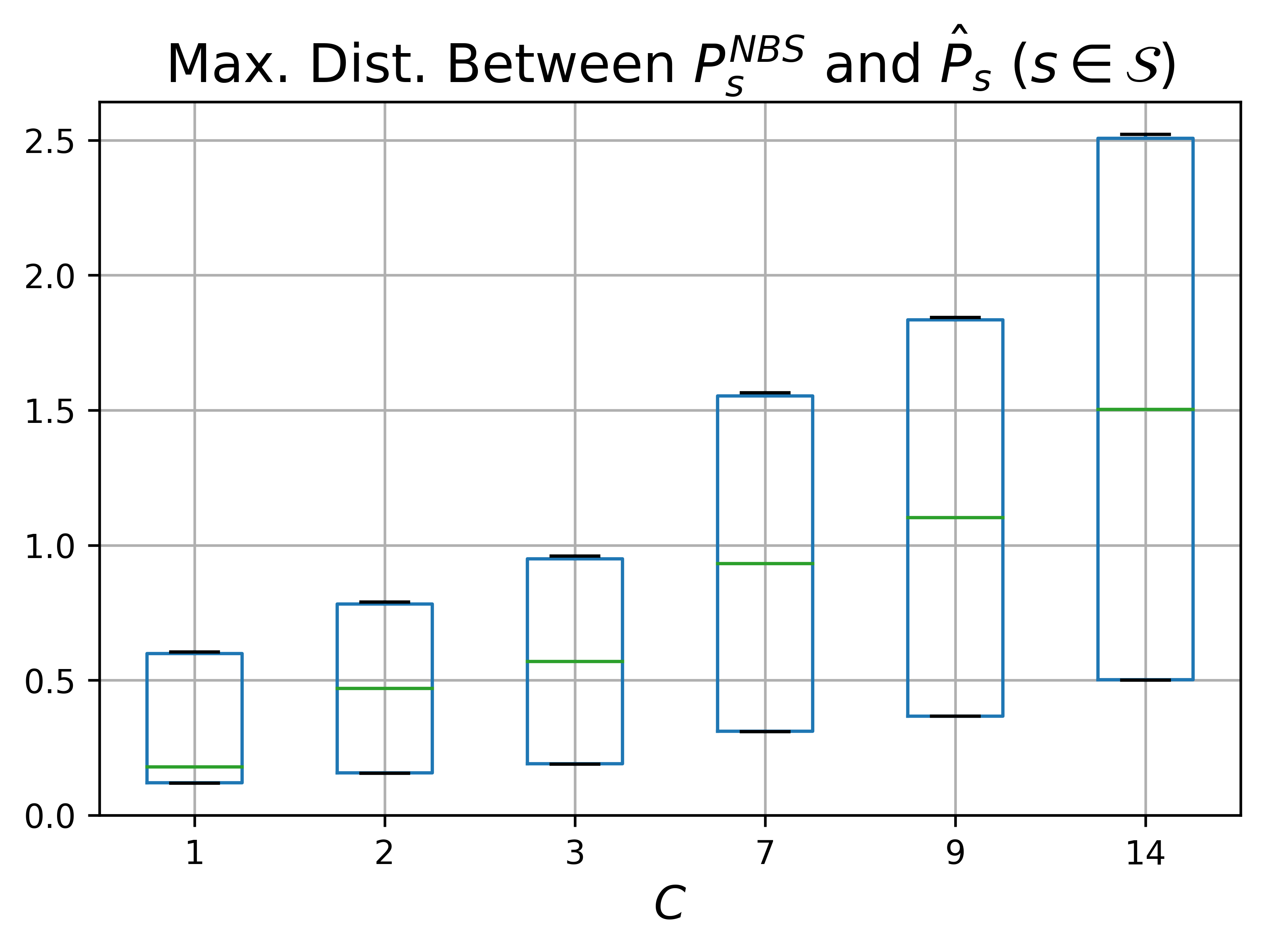}
        \caption{}
        \label{fig:NP_vs_Phat}
    \end{subfigure}
    \caption{Boxplots of $\max_{s \in \m{S}}\sum_{a \in \m{A}} d_{\phi}\l(\bm{P}^y_{s, a}, \hat{\bm{P}}_{s, a}\r)$ for (a) $y=\text{PBS}$ and (b) $y=\text{NBS}$ (binomial).}
    \label{fig:P_hat_dists}
\end{figure}

We also compare the values from NBS with those from PBS in Figure~\ref{fig:v_diffs_NBS}. The quantities plotted here are the mean differences between $\bm{v}^{\text{PBS}}$ and $\bm{v}^{\text{NBS}}$, calculated using:
\begin{equation}\label{eq:v_diffs_NBS}
    \frac{1}{\lvert S \rvert}\sum_{s \in \m{S}}\l(v^{\text{PBS}}_s - v^{\text{NBS}}_s\r).
\end{equation}
Figure~\ref{fig:v_diffs_NBS} shows that the value functions were generally quite close. The smallest and largest mean difference between $\bm{v}^{\text{PBS}}$ and $\bm{v}^{\text{NBS}}$ were $-2.24$ and $1.56$ respectively. However, there is a clear pattern in the value function differences as $C$ increases. For small $C$, Figure~\ref{fig:v_diffs_NBS} shows that $\bm{v}^{\text{NBS}}$ and $\bm{v}^{\text{PBS}}$ were very close, with $\bm{v}^{\text{PBS}}$'s mean value (taken over $s$) being slightly higher. However, as $C$ increases past 3 we see a clear pattern of $\bm{v}^{\text{NBS}}$'s mean value becoming larger than $\bm{v}^{\text{PBS}}$'s. The magnitude of this difference grows as $C$ increases. This indicates that the worst-case distributions resulting from the parametric ambiguity set can be \textit{worse} than those from the non-parametric set, i.e.\ they can lead to lower worst-case rewards.
This result can be explained by differences between the parametric and non-parametric ambiguity sets. Recall that the non-parametric ambiguity set $\m{P}_s$ is defined as containing all distributions $\bm{P}_s$ that satisfy the inequality $\sum_{a \in \m{A}} d_{\phi}\l(\bm{P}_{s, a}, \hat{\bm{P}}_{s, a}\r) \le \kappa$.
In contrast, the parametric ambiguity set $\Theta^{\alpha}_s$ is defined using an inequality restricting the distance from $\hat{\bm{\theta}}_s$ that $\bm{\theta}_s$ can take. This does not restrict the distance from $\hat{\bm{P}}_s$ that the parametric worst-case can take in the same way as the non-parametric ambiguity set does. To understand this, we evaluate $\max_{s \in \m{S}}\sum_{a \in \m{A}} d_{\phi}\l(\bm{P}^y_{s, a}, \hat{\bm{P}}_{s, a}\r)$ for $y \in \{\text{PBS}, \text{NBS}\}$, for every instance solved. We provide boxplots of these values in Figure~\ref{fig:P_hat_dists}. Figure~\ref{fig:P_hat_dists} shows that, for every value of $C$, the parametric worst-case distribution was allowed to be further from $\hat{\bm{P}}_s$ than the non-parametric worst-case. As $C$ increases, difference between the maximum distances for the parametric and non-parametric worst-case distributions increases, explaining why the value functions are more different for large $C$. This happens since larger $C$ occurs when $A$ is larger, meaning the LHS of the inequalities defining the ambiguity sets are sums of more terms, and also $\chi^2_{oA, 1-\alpha}$ is increasing in $A$.

In general, this means that a distribution being binomial may lead to its inclusion in the parametric confidence set even though it is further from $\hat{\bm{P}}$ than any distribution in the non-parametric confidence set. Similarly, distributions that are not binomial need to be much closer to $\hat{\bm{P}}$ in order to be considered as candidates for the true distribution. The result in Figure~\ref{fig:P_hat_dists} suggests that, for this problem, the parametric ambiguity sets are more risk-averse.  It is worth noting, however, that NBS's values being slightly higher on average does not necessarily mean that more long-run reward would be obtained under the non-parametric model. This depends on the initial distribution $\bm{Q}$. For example, studying the value functions we see that $v^{\text{PBS}}_0 \ge v^{\text{NBS}}_0$ was true in 66\% of instances, $v^{\text{PBS}}_1 \ge v^{\text{NBS}}_1$ was true in 65\% of instances, and $v^{\text{PBS}}_2 \ge v^{\text{NBS}}_2$ was true in 50.4\% of the instances with $C \ge 2$. If the initial distribution satisfied, for example, $Q_s > 0$ only for $s \le 2$, then the non-parametric approach would typically not achieve more long-term expected reward. 

We now compare the policies from PBS and NBS. The first characteristic we study is determinism. In these experiments, we find that PBS's optimal policy was deterministic in 42\% of instances. NBS's policy was deterministic in 16\% of instances. This indicates that the parametric model is more likely to yield deterministic policies than the non-parametric model. In addition, we studied the expected actions for each state in order to determine how conservative each model is. Specifically, we calculated $\sum_{a \in \m{A}} a\pi^y_{s, a}$ for each $s$, for $y \in \{\text{PBS}, \text{NBS}\}$ for every instance that we solved. Generally speaking, the expected actions under PBS and NBS were quite similar. However, we found two main results. Firstly, when the inventory level $s$ is below 13, PBS will typically purchase slightly more stock. For the very smallest states, PBS only ordered between 5\% and 6\% more stock than NBS. However, for $s=10$ and $s=11$ PBS ordered 31\% and 24\% more respectively. Secondly, when the stock level is high, i.e.\ 13 or 14, NBS will typically purchase more stock. This was most noticeable for $s=14$, where NBS ordered 38\% more than PBS on average.  From this, we can conclude that PBS's policies are typically less conservative for the majority of states, but that NBSs's are slightly less conservative for the largest 2 states. 

It may seem odd that either algorithm makes positive orders when $s=C=14$ since any stock above $C$ is lost. This occurs due to differences in the worst-case distributions for different actions. For example, when in state $C$ it may be the case that ordering a non-zero amount of stock leads to a much higher worst-case probability of then transitioning to state zero (and hence selling all stock). For example, this may happen when $\hat{p}_{C, a}$ is much larger for $a > 0$ than for $a=0$, which can occur simply due to sampling variation. In some cases, spending some additional multiple of $w$ to purchase wasted stock results in a higher expected reward due to the fact that the newsvendor is then much more likely to sell all of their stock. Since NBS ordered more in the higher states, clearly this was of more benefit under NBS's worst-case distributions than PBS's. If this is not something that the newsvendor would like to allow, the policy can always be constrained to enforce that $\pi_{s, a} = 0$ for all $a \in \m{A}$ such that $s+a >C$. 


\subsection{Numerical experiments with Poisson demands}\label{sec:poisson}

We now carry out the same experiments as in Section~\ref{sec:binom}, but where $X_{s, a} \sim \text{Pois}(\lambda^0_{s, a})$ for $(s, a) \in \m{S} \times \m{A}$. In Section~\ref{sec:pois_AS} we formulate the Poisson ambiguity sets. Following this, we describe the results of the experiments in Sections~\ref{sec:times_pois} and~\ref{sec:P_vs_NP_pois}.

\subsubsection{Ambiguity sets}\label{sec:pois_AS}

Suppose that $X_{s, a} \sim \text{Pois}(\lambda^0_{s, a})$ for $(s, a) \in \m{S} \times \m{A}$, and therefore:
\begin{equation}
    f_{X_{s, a}}(x|\lambda^0_{s,a}) = \frac{(\lambda^0_{s, a})^x \exp(-\lambda^0_{s, a})}{x!} \quad (x \in \mathbb{N}_0).    
\end{equation}
The distribution of $\bm{X}$ is uniquely parameterised by $\bm{\lambda}^0 = (\lambda^0_{s, a})_{(s, a) \in \m{S} \times \m{A}}$. Similarly to in Section~\ref{sec:binom_AS}, we have $o=1$ and suppose that we take the sample $\bm{x}_{s, a} = (x^1_{s, a}, \dots, x^N_{s, a})$ from $X_{s, a}$ for each $(s, a) \in \m{S} \times \m{A}$. Then, the MLE $\hat{\lambda}_{s, a}$ of $\lambda^0_{s, a}$ is given by:
\begin{equation}
    \hat{\lambda}_{s, a} = \frac{\sum_{j=1}^N x^j_{s, a}}{N} \fa (s, a) \in \m{S} \times \m{A}.
\end{equation}
The Fisher information is now given by:
\begin{equation}
    I_{\E}(\hat{\lambda}_{s,a}) = \frac{N}{\hat{\lambda}_{s, a}}.
\end{equation}
Hence, an approximate $100(1-\alpha)\%$ confidence set for $\bm{\lambda}^0_s$ (for large $N$) is given by:
\begin{equation}\label{eq:pois_conf}
    \Theta^{\alpha}_s = \l\{\bm{\lambda}_s \in \R^A_+: \sum_{a \in \m{A}} \frac{N(\lambda_{s, a} - \hat{\lambda}_{s, a})^2}{\hat{\lambda}_{s, a}} \le \chi^2_{A, 1-\alpha} \r\}.
\end{equation}
As in Section~\ref{sec:binom_AS}, we will construct a discretisation of $\Theta^{\alpha}_s$ by creating a set $\Theta^{\text{base}}_s$ with $\Theta^{\alpha}_s \subseteq \Theta^{\text{base}}_s$ and discretising this set. Then, we extract elements of this discrete set that also lie in $\Theta^{\alpha}_s$. The definition of $\Theta^{\alpha}_s$ implies that every $\bm{\lambda}_s \in \Theta^{\alpha}_s$ satisfies:
\begin{equation}
    \lambda_{s, a} \in \lambda^{\text{I}}_{s, a} = \l[\max\l\{0, \hat{\lambda}_{s, a} - \sqrt{\frac{\chi^2_{A, 1-\alpha}\hat{\lambda}_{s, a}}{N}}\r\}, \hat{\lambda}_{s, a} + \sqrt{\frac{\chi^2_{A, 1-\alpha}\hat{\lambda}_{s, a}}{N}}\r]
\end{equation}
for all $a \in \m{A}$. Hence, we define $\Theta^{\text{base}}_s = \lambda^{\text{I}}_{s, 1} \times \dots \times \lambda^{\text{I}}_{s, A}$ and
we have $\Theta^{\alpha}_s \subseteq \Theta^{\text{base}}_s$. Furthermore, define $\lambda^l_{s, a}$ and $\lambda^u_{s, a}$ as the lower and upper bounds of $\lambda^{\text{I}}_{s, a}$ for each $(s, a) \in \m{S} \times \m{A}$. We calculate the following discretisations of each $\lambda^{\text{I}}_{s, a}$, containing $M$ points, as follows:
\begin{equation}
    \tilde{\lambda}^{\text{I}}_{s, a} = \l\{\lambda^l_{s, a} + m \frac{\lambda^u_{s, a} - \lambda^l_{s, a}}{M - 1}\r\}.
\end{equation}
Then, a discretisation of $\Theta^{\text{base}}_s$ is given by $(\Theta^{\text{base}}_s)' = \tilde{\lambda}^{\text{I}}_{s, 1} \times \dots \times \tilde{\lambda}^{\text{I}}_{s, A}$. Finally, a discretisation of $\Theta^{\alpha}_s$ is given by $(\Theta^{\alpha}_s)' = (\Theta^{\text{base}}_s)' \cap \Theta^{\alpha}_s$. 

\subsubsection{Experimental design}\label{sec:exp_des_pois}

All of the main parameters for these experiments are the same as in Section~\ref{sec:exp_des}. The only differences are with respect to the parameters used in the parametric bisection search algorithm of Section~\ref{sec:param_proj}. We again use $\theta^{\min}_{s, a} = 0$ for all $(s, a) \in \m{S} \times \m{A}$, but since $\lambda_{s, a}$ is technically not bounded from above, there is no obvious value for $\theta^{\max}_{s, a}$. However, since any root of $\sum_{s' \in \m{S}} P^{\theta}_{s, a, s'} b_{s'} = \beta$ that has $\lambda_{s, a} > \lambda^u_{s, a}$ cannot be an element of a $\bm{\lambda}$ that lies in the ambiguity set, we set $\theta^{\max}_{s, a} = \lambda^u_{s, a}$ for all $(s, a) \in \m{S} \times \m{A}$. Since this creates a wider range for the potential roots than for the binomial case, so we set $\tilde{\epsilon} = \frac{\theta^{\max}_{s, a} - \theta^{\min}_{s, a}}{100}$. This ensures that the same number of intervals are used here as in the binomial case, where we used $\tilde{\epsilon} = 0.01 \l(= \frac{1-0}{100}\r)$.

\subsubsection{Times taken}\label{sec:times_pois}

We now present the results of our experiments for Poisson demands. Of the 810 instances ran, we found that LP and CS timed out in 54. PBS timed out in 2, although this was again during the final step of finding the optimal policy with CS. PBS never timed out during value iteration. QP and NBS did not time out in any instance, but QP again resulted in convergence issues. When using QP to solve the Bellman updates, value iteration failed to converge in 384 instances. Table~\ref{tab:VI_times_pois} summarises the times taken to run value iteration. Similar results to the binomial case can be found here, with PBS being faster than NBS, CS being slightly faster than LP, and QP being fast per iteration.
\begin{table}[htbp!]
    \centering
    \begin{tabular}{llll}
\toprule
Algorithm &   Mean Time &    Max Time & Mean Time Per Iteration \\
\midrule
      PBS &  0:00:24.44 &  0:03:59.44 &              0:00:00.85 \\
       CS &  0:17:26.57 &     4:00:00 &              0:04:20.96 \\
       LP &  0:23:45.70 &     4:00:00 &              0:16:16.30 \\
       QP &  0:05:37.69 &  0:58:13.20 &              0:00:00.38 \\
      NBS &  0:00:42.01 &  0:10:52.69 &              0:00:01.45 \\
\bottomrule
\end{tabular}

    \caption{Summary of times taken to run value iteration (Poisson)}
    \label{tab:VI_times_pois}
\end{table}

\begin{figure}[htbp!]
    \begin{subfigure}[t]{0.45\textwidth}
        \centering
        \includegraphics[width=0.9\textwidth]{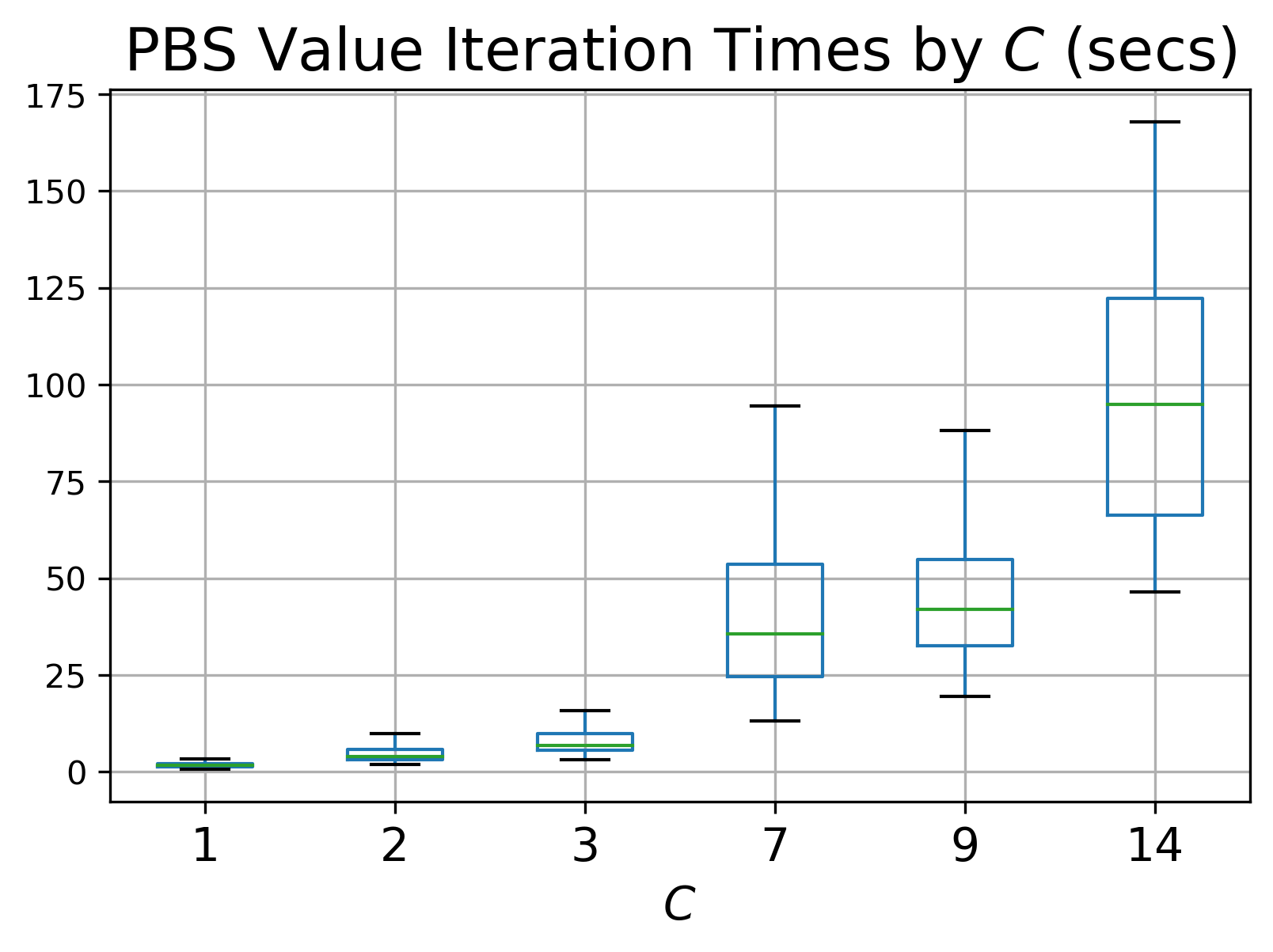}
        \caption{}
        \label{fig:PBS_by_C_pois}
    \end{subfigure}
    \begin{subfigure}[t]{0.45\textwidth}
        \centering
        \includegraphics[width=0.90\textwidth]{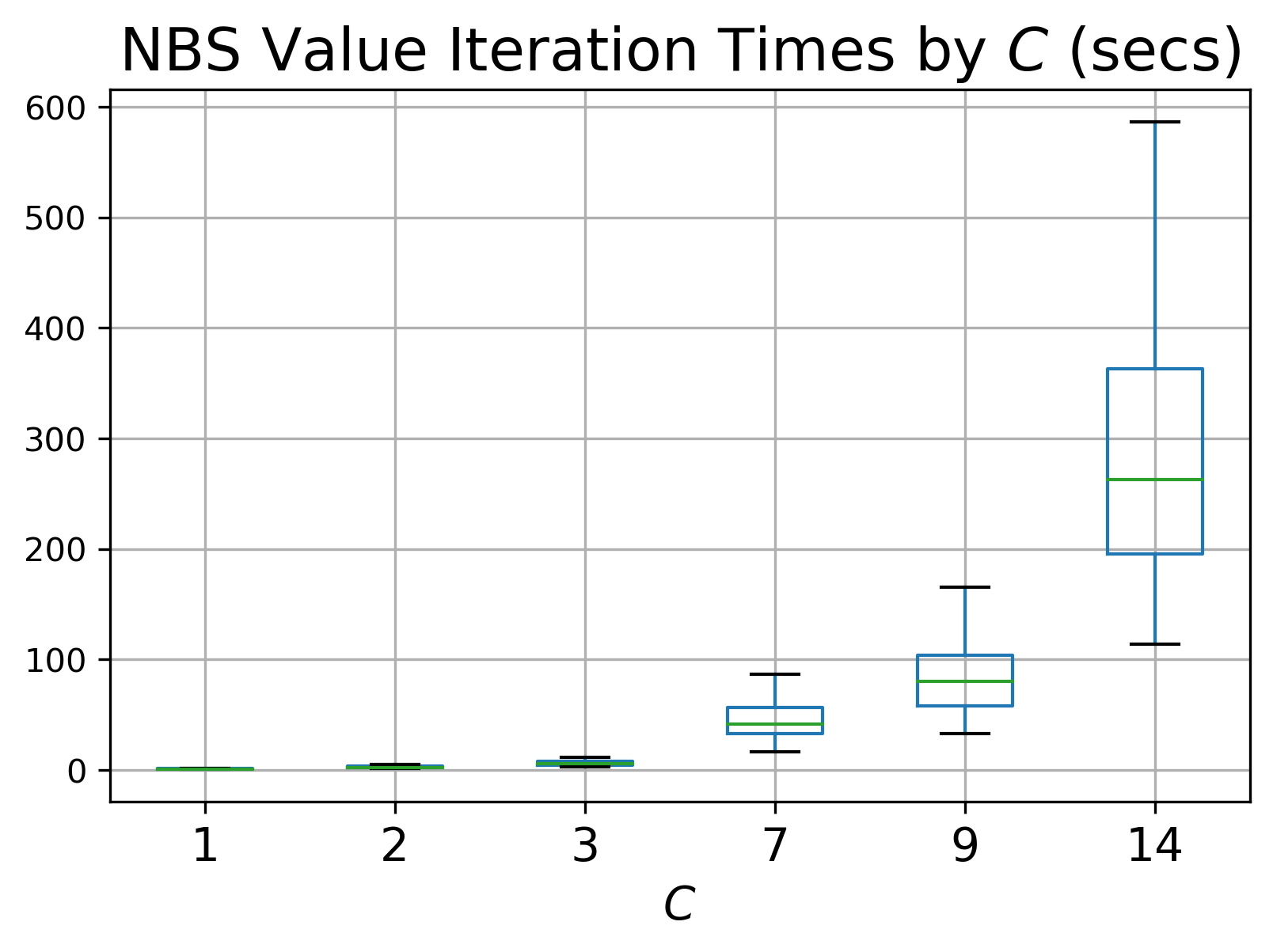}
        \caption{}
        \label{fig:NBS_by_C_pois}
    \end{subfigure}
    \caption{Boxplots of value iteration run times of (a) PBS and  (b) NBS, by $C$ (Poisson).}
    \label{fig:times_C_pois}
\end{figure}

Figure~\ref{fig:times_C_pois} compares the value iteration run times of PBS and NBS more closely. It shows that, while NBS was slightly faster for small $C$, PBS scales much better with large $C$. For $C = 14$, PBS typically took no more than 3 minutes to finish value iteration. However, NBS took up to 10 minutes. 

\begin{figure}[htbp!]
    \begin{subfigure}[t]{0.45\textwidth}
        \centering
        \includegraphics[width=0.85\textwidth]{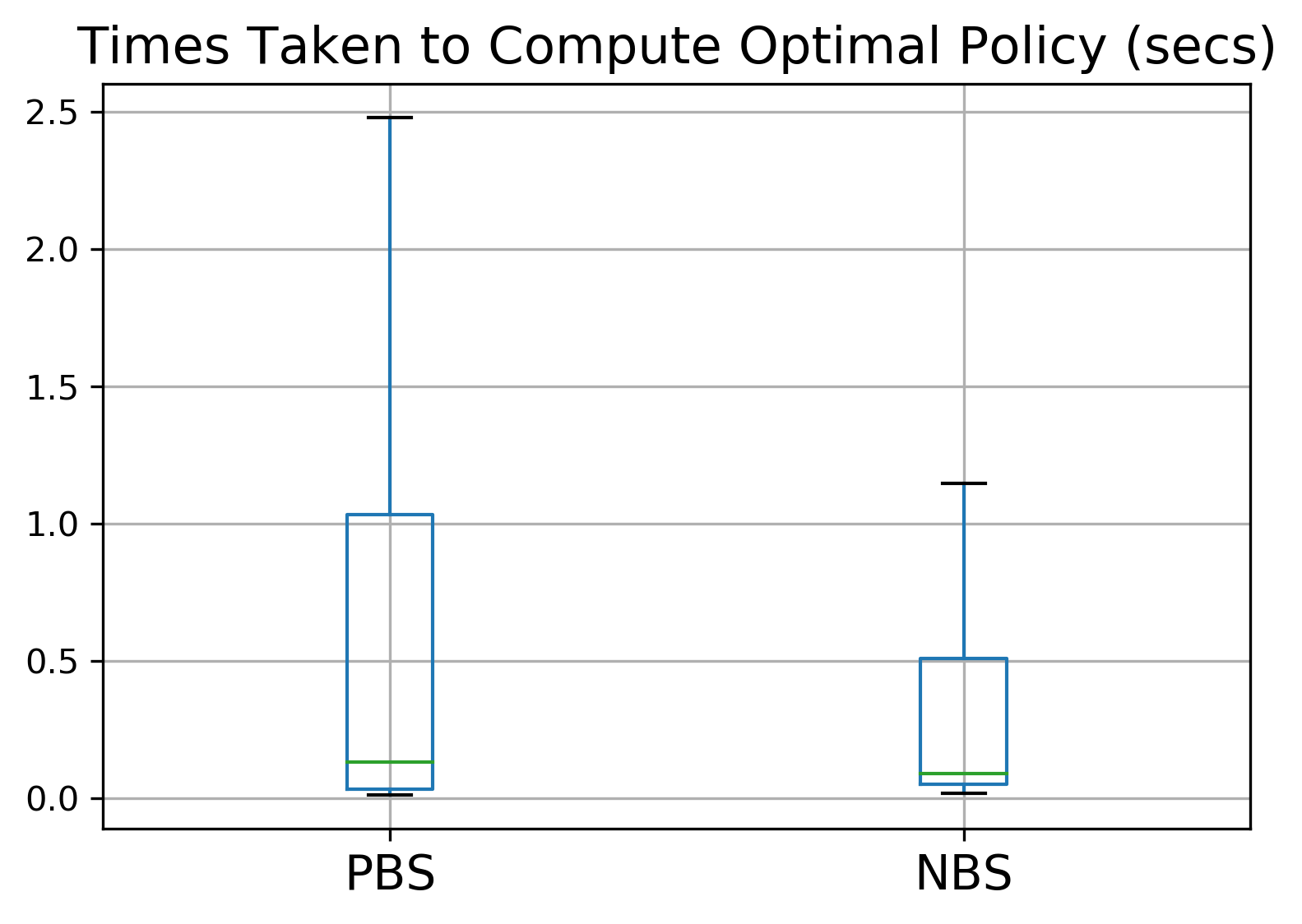}
        \caption{}
        \label{fig:policy_times_pois}
    \end{subfigure}
    \begin{subfigure}[t]{0.45\textwidth}
        \centering
        \includegraphics[width=0.85\textwidth]{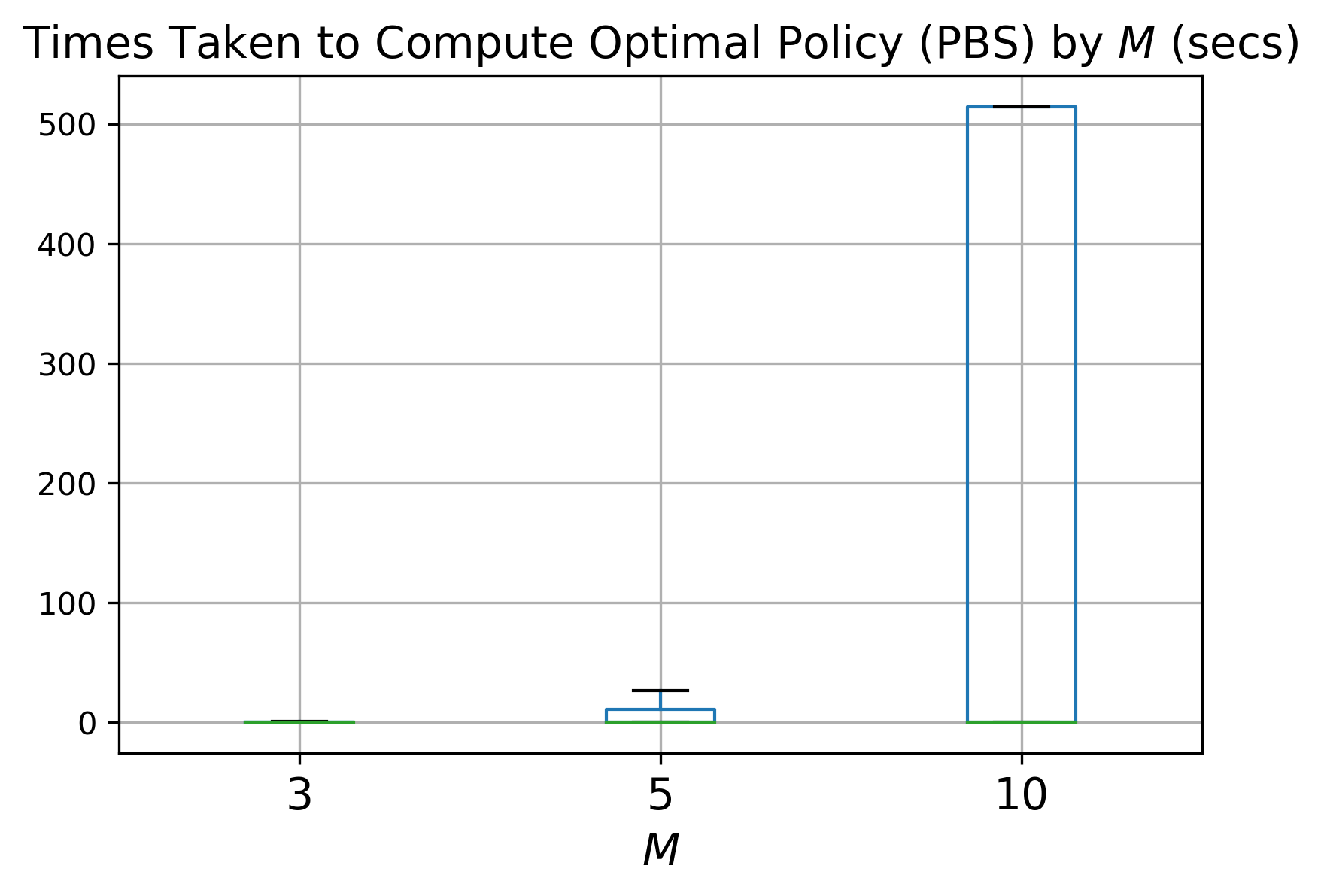}
        \caption{}
        \label{fig:policy_times_M_pois}
    \end{subfigure}
    \caption{Boxplots of times taken to compute optimal policy for (a) both BS algorithms and (b) PBS by $M$ (Poisson).}
    \label{fig:policy_times_pois_}
\end{figure}

While PBS is fast at finding the optimal values, it also took much longer to find a policy after running value iteration with PBS than with NBS. However, as before, this is due to CS being slow for large instances, and has nothing to do with PBS itself. On average, it took 8 minutes and 13 seconds for CS to find a policy for PBS's values, as opposed to 0.38 seconds for QP to compute the policy for NBS. However, the parametric times were skewed by large instances; the median time to compute the policy for PBS was 0.14 seconds. Figure~\ref{fig:policy_times_pois_} shows boxplots of the times taken to find the policy for PBS's and NBS's values. While Figure~\ref{fig:policy_times_pois} suggests that the speeds were similar for PBS and NBS's values for the majority of instances, Figure~\ref{fig:policy_times_M_pois} show the drastic times taken under the parametric model for $M = 10$. As before, this is due to CS's slowness, not PBS's. It is likely that if PBS were to be used in practice, a heuristic algorithm could be used in CS's place that would drastically speed up these times.

\subsubsection{Comparison of value functions, distributions and policies}\label{sec:P_vs_NP_pois}

\begin{figure}[htbp!]
        \begin{subfigure}[t]{0.45\textwidth}
        \centering
        \includegraphics[width=0.9\textwidth]{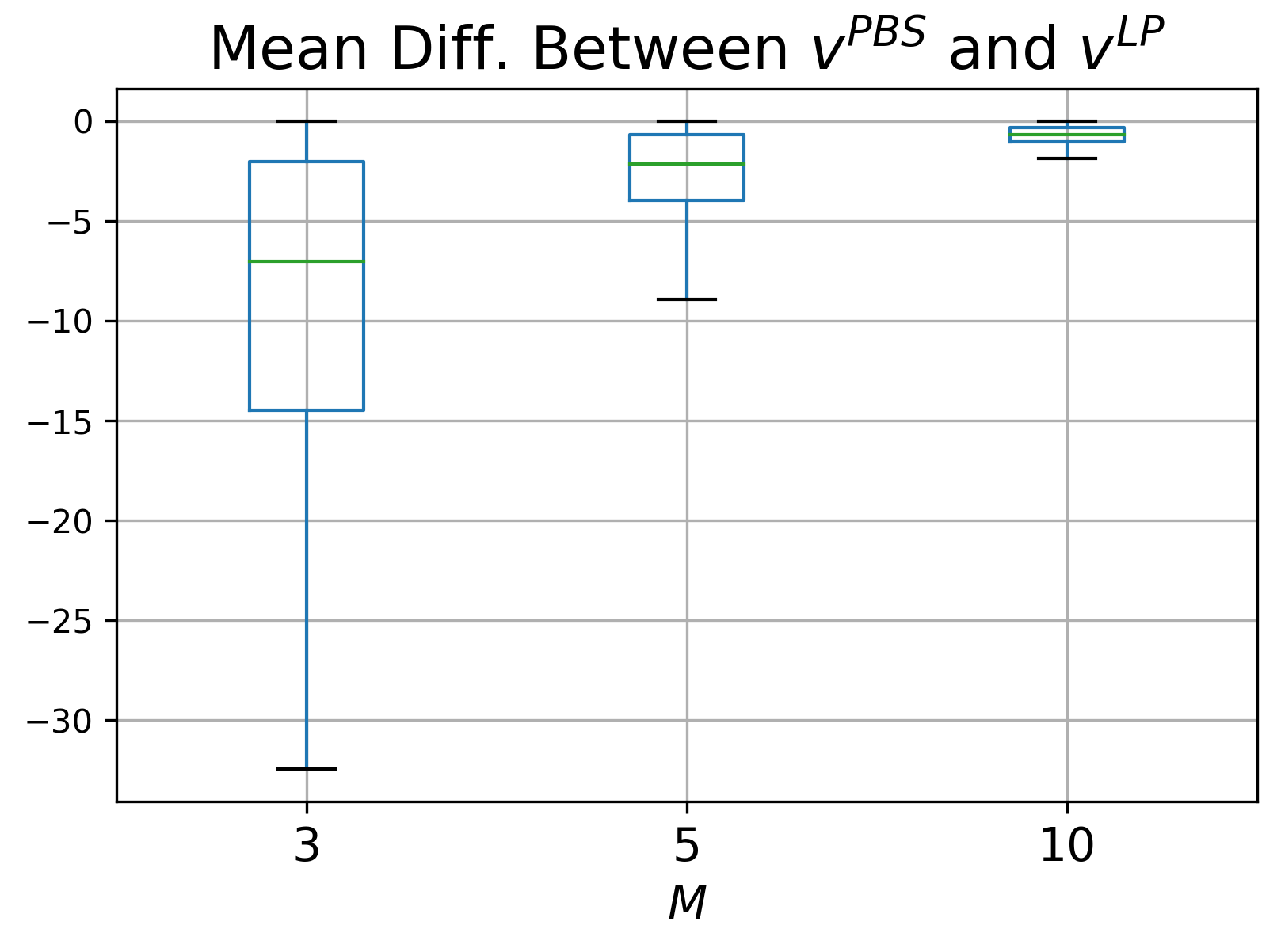}
        \caption{}
        \label{fig:LP_PBS_M_pois}
    \end{subfigure}
        \begin{subfigure}[t]{0.45\textwidth}
        \centering
        \includegraphics[width=0.90\textwidth]{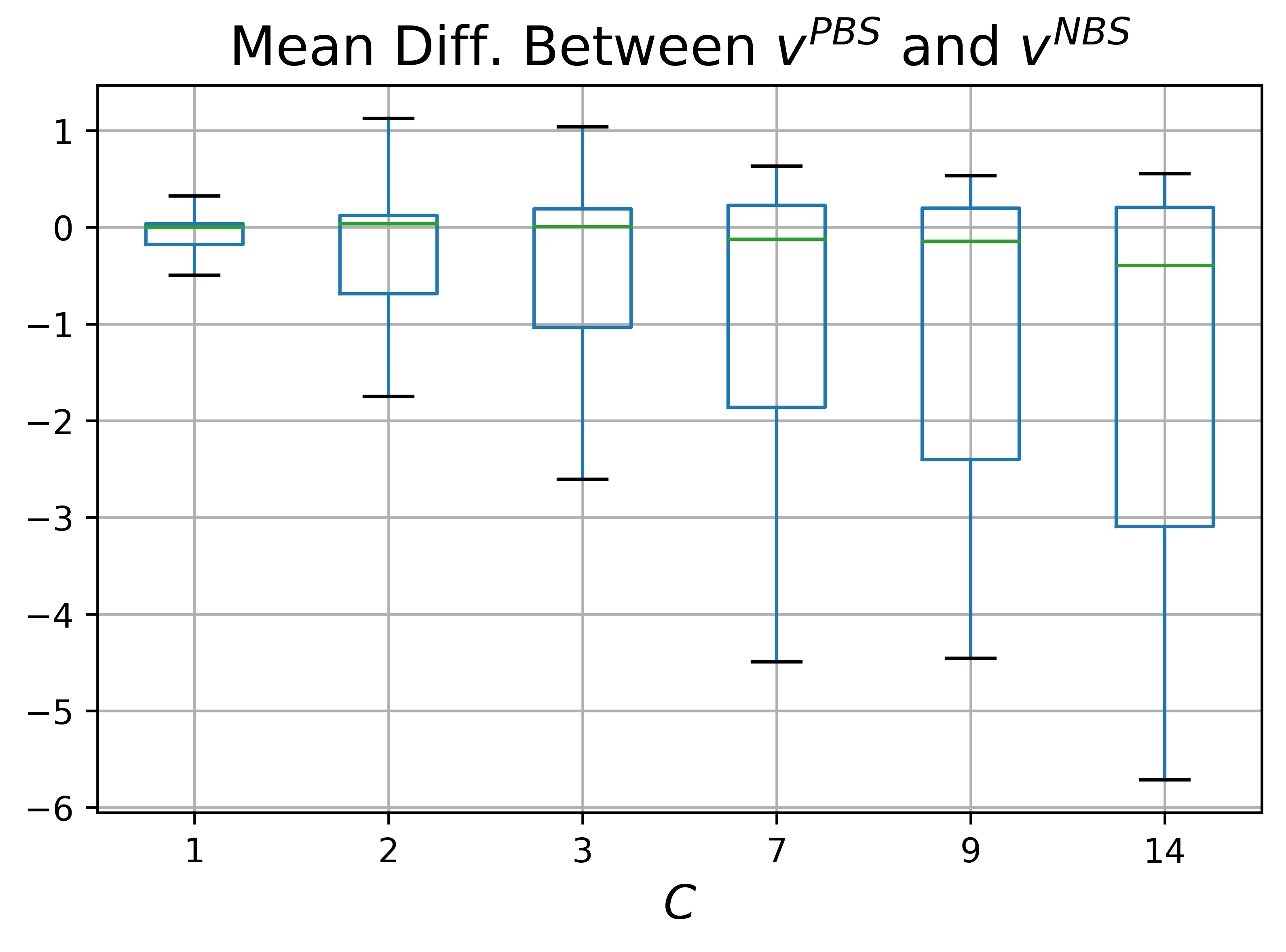}
        \caption{}
        \label{fig:PBS_NBS_pois}
    \end{subfigure}
    \caption{Boxplots of mean difference between $\bm{v}^{\text{PBS}}$ and (a) $\bm{v}^{\text{LP}}$ by $M$ and (b) $\bm{v}^{\text{NBS}}$ by $C$ (Poisson).}
    \label{fig:v_pois}
\end{figure}

We now compare the outputs from the parametric and non-parametric models. Firstly, we compare the value functions from PBS with those from LP and NBS. Boxplots comparing these values are shown in Figure~\ref{fig:v_pois}. As a reminder, these plots show $\frac{1}{\lvert S \rvert}\sum_{s \in \m{S}}\l(v^{\text{PBS}}_s - v^{y}_s\r)$ for $y \in \{\text{LP}, \text{NBS}\}$. Figure~\ref{fig:LP_PBS_M_pois} shows the convergence of LP's values to PBS's as $M$ increases. As for the binomial model, we see that LP's values are always higher, and they grow closer to PBS's on average as $M$ increases. In addition, Figure~\ref{fig:PBS_NBS_pois} shows that PBS and NBS's values are similar for small $C$, but NBS's values are typically higher than PBS's for large $C$.

\begin{figure}[htbp!]
    \begin{subfigure}[t]{0.45\textwidth}
        \centering
        \includegraphics[width=0.9\textwidth]{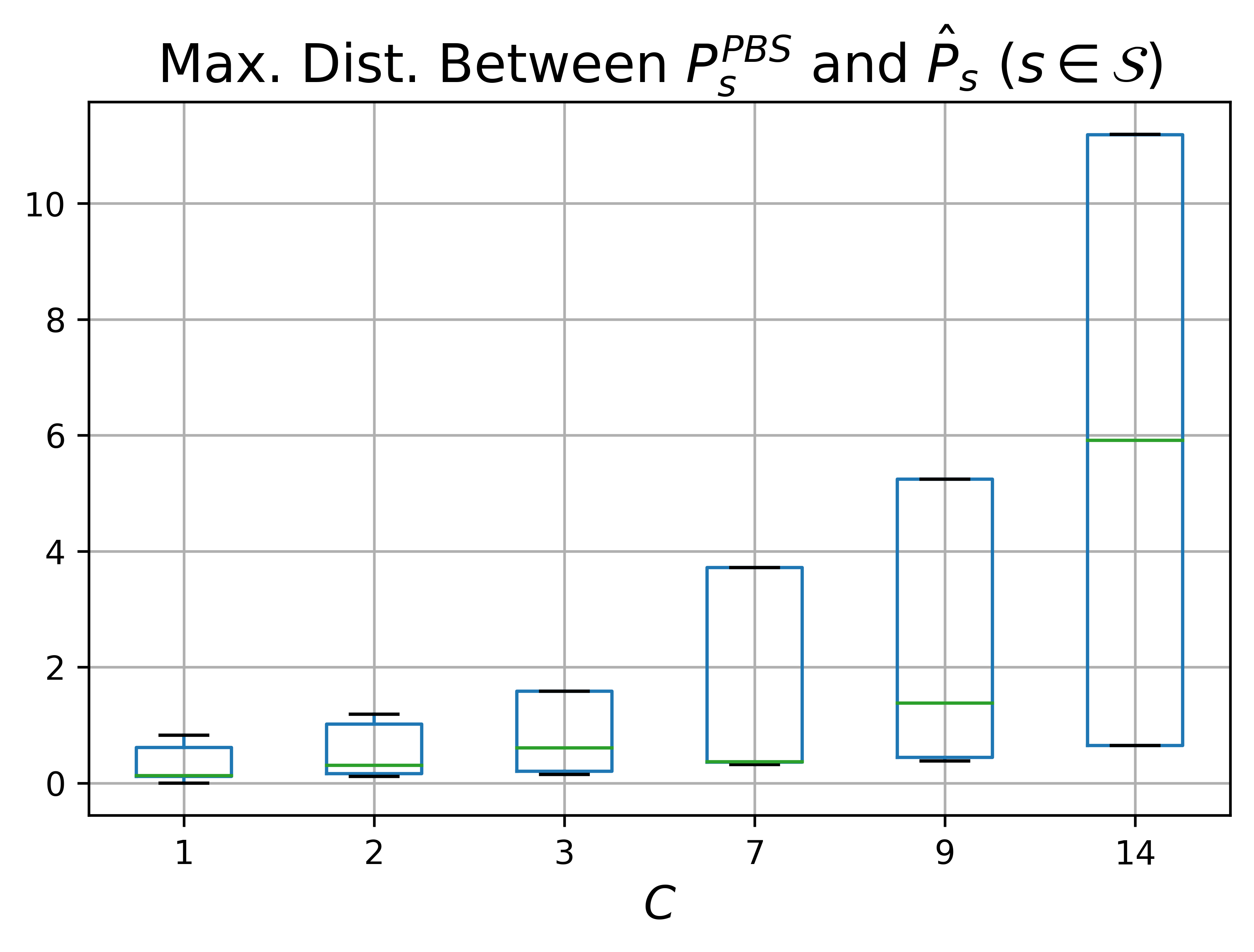}
        \caption{}
        \label{fig:P_vs_Phat_pois}
    \end{subfigure}
    \begin{subfigure}[t]{0.45\textwidth}
        \centering
        \includegraphics[width=0.90\textwidth]{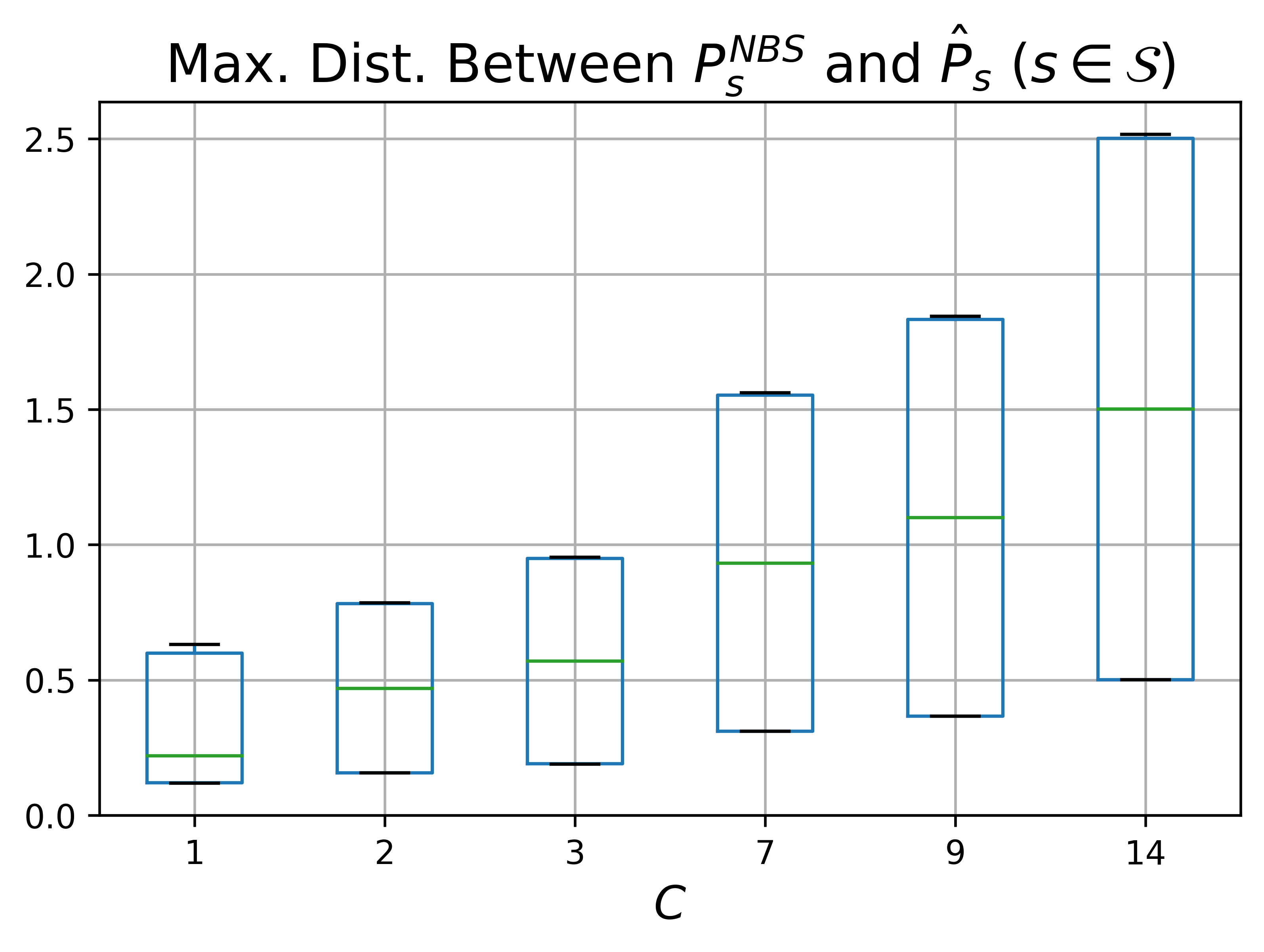}
        \caption{}
        \label{fig:NP_vs_Phat_pois}
    \end{subfigure}
    \caption{Boxplots of $\max_{s \in \m{S}}\sum_{a \in \m{A}} d_{\phi}\l(\bm{P}^y_{s, a}, \hat{\bm{P}}_{s, a}\r)$ for (a) $y=\text{PBS}$ and (b) $y=\text{NBS}$ (Poisson).}
    \label{fig:P_hat_dists_pois}
\end{figure}

In order to confirm that the same property of the ambiguity sets is responsible for this as for the binomial model, we plot the maximum distances from $\hat{\bm{P}}_s$ attained by $\bm{P}^{\text{PBS}}_s$ and $\bm{P}^{\text{NBS}}_s$ in Figure~\ref{fig:P_hat_dists_pois}. This plot again suggests that the parametric worst-case distributions are much further from $\hat{\bm{P}}_s$, especially for large $C$. 

Finally, we compare the policies from the PBS and NBS. Similarly to the binomial model, we find that PBS's policies were more often deterministic. Specifically, 32\% of PBS's policies were deterministic where only 14\% of NBS's were. In addition, we also find that the parametric policies were slightly less conservative in terms of purchasing, for most states. However, for the largest 2 states, i.e.\ $s=13, 14$, the non-parametric policies were less conservative. NBS's polices ordered 30\% more stock in these final two states than PBS's, whereas PBS's policies ordered between 5\% and 25\% more in the lower states.

\section{Conclusions and further research}

In this paper, we studied robust Markov decision processes under parametric transition distributions. We focused on robust value iteration for $s$-rectangular ambiguity sets in particular. Based on a fast projection-based bisection search algorithm found in the literature for robust MDPs with $\phi$-divergence ambiguity sets, we created a projection-based bisection algorithm for the parametric model in the case where the transition distribution is parametrised by one parameter. We also presented two other algorithms for solving the robust Bellman update, a linear programming algorithm and a cutting surface algorithm. These algorithms both discretise the ambiguity set for the true parameter in order to create a linear programming reformulation of the robust Bellman update. In addition, we showed how to use maximum likelihood estimation to create confidence sets for use as ambiguity sets in the parametric model.  

In order to test our algorithms, we presented a dynamic multi-period newsvendor model and applied them to it. In particular, we carried out numerical experiments on the case where the demands in the newsvendor problem are binomial and Poisson. In both cases, we solved the non-parametric model in addition to the parametric model, in order to compare run times and solutions. We found a number of main results. Firstly, our parametric bisection search algorithm was very fast at finishing value iteration. In fact, it was faster than its non-parametric equivalent and offered significant time savings in comparison with the linear programming and cutting surface algorithms. This is due to the fact that the bisection algorithm does not rely on any discretisation of the ambiguity set, and hence does not need to carry out the large amount of pre-computation required by the solver-based algorithms.  However, since our bisection algorithm does not return an optimal policy, one of the two solver-based algorithms had to be used to generate the policy after value iteration ended. This meant that using the parametric model sometimes resulted in large overall run times. Comparing the solutions from the parametric and non-parametric models, we found that the non-parametric value functions were typically higher than the parametric ones on average. This was due to the result that the parametric confidence sets allowed the corresponding worst-case distribution to be much further from the nominal distribution than was allowed by the non-parametric set. However, we also found that the parametric policies were less conservative than the non-parametric, for all states apart from the largest two.

There are two main directions for future research arising from this paper. The most obvious direction is with regards to computing the optimal policy. Since both of our solver-based algorithms were very slow at this, it would be beneficial to find a faster way to extract the policy once our bisection algorithm finishes value iteration. Another potential area for further research is  with regards to discretisation. As is the case in parametric distributionally robust optimisation, discretisation of the ambiguity set is required to create a tractable linear programming model that approximates the true problem. This means that the transition matrices for every parameter in the ambiguity set must be computed prior to building the model, and also that the resulting model becomes very large for large ambiguity sets. Hence, we would like to study ways in which to circumvent the need for discretisation and reduce overall solution times. 

\Large{\textbf{Acknowledgements}}

\normalsize
We would like to acknowledge the support of the Engineering and Physical
Sciences Research Council funded (EP/L015692/1) STOR-i Centre for Doctoral
Training. We would like to thank BT for their funding, and Russell Ainslie, Mathias Kern and
Gilbert Owusu from BT for their support. 

\bibliographystyle{apalike}
\bibliography{bib}

\newpage
\begin{appendices}

\section{Derivation of reformulation of robust Bellman update}\label{sec:reformulation}

\subsection{General reformulation}\label{sec:gen_reform}

The inner problem of~\eqref{eq:robust_bell} can be written as:
\begin{align}
    \min_{\bm{P}_s} &\sum_{a \in \m{A}}\pi_{s,a} \sum_{s' \in \m{S}} P_{s, a, s'}B_{s, a, s'}\\
    \text{s.t. }& \sum_{a \in \m{A}} d_a(\bm{P}_{s, a}, \hat{\bm{P}}_{s, a})\le \kappa \\
    & \sum_{s' \in \m{S}} P_{s, a,s'} = 1 \fa a \in \m{A},
\end{align}
with $\bm{B}_{s, a} = \bm{r}_{s,a} + \gamma \bm{v}$ and $\bm{v} = \bm{v}^n$. The Lagrangian of this problem is given by:
\begin{align}
    L(\bm{\pi}, \bm{\nu}, \eta) &= \sum_{a \in \m{A}}\pi_{s,a} \sum_{s' \in \m{S}} P_{s, a, s'}B_{s, a, s'} + \eta\l(\sum_{a \in \m{A}} d_a(\bm{P}_{s, a}, \hat{\bm{P}}_{s, a})- \kappa\r) + \sum_{a \in \m{A}} \nu_a \l(1- \sum_{s' \in \m{S}}P_{s, a, s'}\r) \\
    &= -\kappa \eta + \nu + \sum_{a \in \m{A}} \sum_{s'\in\m{S}}\l[ \pi_{s, a} P_{s, a, s'}B_{s, a,s'} + \eta \hat{P}_{s, a,s'}\phi\l(\frac{P_{s, a, s'}}{\hat{P}_{s, a,s'}}\r) - \nu_a P_{s, a,s'}\r],
\end{align}
with $\nu = \sum_{a \in \m{A}} \nu_a$. Therefore, the objective of the dual of the inner problem is given by:
\begin{align}
    g(\bm{\nu}, \eta) &= \inf_{\bm{P}_s \ge 0}\l\{
        -\kappa \eta + \nu + \sum_{a \in \m{A}} \sum_{s'\in\m{S}}\l[ \pi_{s, a} P_{s, a, s'}B_{s, a,s'} + \eta \hat{P}_{s, a,s'}\phi\l(\frac{P_{s, a, s'}}{\hat{P}_{s, a,s'}}\r) - \nu_a P_{s, a,s'}\r]
    \r\}\\
    &=  -\kappa \eta + \nu + \sum_{a \in \m{A}} \sum_{s'\in\m{S}}\inf_{P_{s, a, s'} \ge 0}\l\{
     \pi_{s, a} P_{s, a, s'}B_{s, a,s'} + \eta \hat{P}_{s, a,s'}\phi\l(\frac{P_{s, a, s'}}{\hat{P}_{s, a,s'}}\r) - \nu_a P_{s, a,s'}
    \r\}\\
    &= -\kappa \eta + \nu + \sum_{a \in \m{A}} \sum_{s'\in\m{S}}\eta \hat{P}_{s, a, s'}\inf_{P_{s, a, s'} \ge 0}\l\{
     \frac{P_{s, a, s'}}{\hat{P}_{s, a, s'}}\frac{\pi_{s, a} B_{s, a,s'}}{\eta} + \phi\l(\frac{P_{s, a, s'}}{\hat{P}_{s, a, s'}}\r) - \frac{\nu_a}{\eta} \frac{P_{s, a,s'}}{\hat{P}_{s, a, s'}}
    \r\}\\
    &= -\kappa \eta + \nu + \sum_{a \in \m{A}} \sum_{s'\in\m{S}}\eta \hat{P}_{s, a, s'}\inf_{P_{s, a, s'} \ge 0}\l\{
     \frac{P_{s, a, s'}}{\hat{P}_{s, a, s'}}\frac{\pi_{s, a} B_{s, a,s'} - \nu_a}{\eta} + \phi\l(\frac{P_{s, a, s'}}{\hat{P}_{s, a, s'}}\r)\r\}\\
     &= -\kappa \eta + \nu - \sum_{a \in \m{A}} \sum_{s'\in\m{S}}\eta \hat{P}_{s, a, s'}\sup_{P_{s, a, s'} \ge 0}\l\{
     \frac{P_{s, a, s'}}{\hat{P}_{s, a, s'}}\frac{\nu_a - \pi_{s, a} B_{s, a,s'}}{\eta} - \phi\l(\frac{P_{s, a, s'}}{\hat{P}_{s, a, s'}}\r)\r\}\label{eq:inf_sup}\\
     &=-\kappa \eta + \nu - \sum_{a \in \m{A}} \sum_{s'\in\m{S}}\eta \hat{P}_{s, a, s'}\phi^*\l(\frac{\nu_a - \pi_{s, a} B_{s, a,s'}}{\eta}\r).
\end{align}
Here, we used $\inf(A) = -\sup(-A)$ for any set $A$ in order to replace inf with sup. Therefore, the dual of the inner problem is given by:
\begin{align}
    \max_{\eta, \bm{\nu}}\l\{-\kappa \eta + \nu - \sum_{a \in \m{A}} \sum_{s'\in\m{S}}\eta \hat{P}_{s, a, s'}\phi^*\l(\frac{\nu_a - \pi_{s, a} B_{s, a,s'}}{\eta}\r): \eta \in \R_+, \bm{\nu} \in \R^A\r\}.
\end{align}
Combining with the outer problem, our reformulation is given by:
\begin{align}
    \max_{\bm{\pi}_s, \eta, \bm{\nu}}&\l\{-\kappa \eta + \nu - \sum_{a \in \m{A}} \sum_{s'\in\m{S}}\eta \hat{P}_{s, a, s'}\phi^*\l(\frac{\nu_a - \pi_{s, a} B_{s, a,s'}}{\eta}\r)\r\}\\
    \text{s.t. }& \sum_{a \in \m{A}} \pi_{s, a} = 1 \\
    & \eta \in \R_+, \bm{\nu} \in \R^A.
\end{align}

\subsection{\texorpdfstring{Reformulation for modified $\chi^2$ divergence}{}}\label{sec:mchisq_reform}

For the modified $\chi^2$ distance, we have:
\begin{equation}
    \phi^*(z) = \max\l\{1 + \frac{z}{2}, 0\r\}^2 - 1.
\end{equation}
Hence, we have:
\begin{align}
    \eta\phi^*\l(\frac{\nu_a - \pi_{s, a}B_{s, a, s'}}{\eta}\r) &= \eta\max\l\{1 + \frac{\nu_a - \pi_{s, a}B_{s, a, s'}}{2\eta}, 0\r\}^2 - \eta \\
    &= \frac{1}{4\eta}\max\l\{2\eta + \nu_a - \pi_{s, a}B_{s, a, s'}, 0\r\}^2 - \eta.
\end{align}
We can reformulate this using conic quadratic constraints as follows. Firstly, define the dummy variables $\zeta_{s, a, s'}$ for $a \in \m{A}$ and $s' \in \m{S}$ using the following constraints:
\begin{align}
    \zeta_{s, a, s'} &\ge 2\eta + \nu_a - \pi_{s, a}B_{s, a, s'} \fa a \in \m{A} \fa s' \in \m{S}\\
    \zeta_{s, a, s'} &\ge 0 \fa a \in \m{A} \fa s' \in \m{S}.
\end{align}
Now define dummy variables $u_{s, a, s'} \fa a \in \m{A} \fa s' \in \m{S}$ using:
\begin{equation}
    u_{s, a, s'} \ge \frac{\zeta_{s, a, s'}^2}{\eta}
\end{equation}
which is equivalently represented by:
\begin{equation}
    \sqrt{4\zeta_{s,a,s'}^2 + (\eta - u_{s,a,s'})^2} \le (\eta + u_{s, a, s'}).
\end{equation}
Then, at optimality we will have $\eta\phi^*\l(\frac{\nu_a - \pi_{s, a}B_{s, a, s'}}{\eta}\r) = \frac{1}{4} u_{s, a, s'} - \eta$. Therefore, the CQP reformulation of~\eqref{eq:robust_bell} is given by:
\begin{align}
    \max_{\bm{\pi}_s}&\l\{\nu + \eta (A-\kappa) - \frac{1}{4}\sum_{a \in \m{A}}\sum_{s' \in \m{S}} \hat{P}_{s, a, s'} u_{s, a, s'}\r\}\\
    \text{s.t. } & \sqrt{4\zeta_{s,a,s'}^2 + (\eta - u_{s,a,s'})^2} \le (\eta + u_{s, a, s'}) \fa a \in \m{A} \fa s' \in \m{S}\\
    &\zeta_{s, a, s'} \ge 2\eta + \nu_a - \pi_{s, a}B_{s, a, s'} \fa a \in \m{A} \fa s' \in \m{S}\\
    &\zeta_{s, a, s'} \ge 0 \fa a \in \m{A} \fa s' \in \m{S}\\
    &u_{s, a, s'} \ge 0 \fa a \in \m{A} \fa s' \in \m{S}\\
    &\sum_{a \in \m{A}} \pi_{s, a} = 1 \\
    &\pi_{s, a} \ge 0 \fa a \in \m{A} \\
    &\eta \ge 0\\
    &\bm{\nu} \in \R^A.
\end{align}
Note that the term $A\eta$ comes from $\sum_{a \in \m{A}}\sum_{s' \in \m{S}}\hat{P}_{s, a, s'} \eta = A\eta$.

\section{\texorpdfstring{Solving modified $\chi^2$ distance projection problems}{}}

Since we focus on the modified $\chi^2$  divergence in this paper, we will describe the method for this distance only. 

\subsection{Solution by sorting and subproblems}\label{sec:sort_alg}

The method used by~\cite{ho2022robust} consists of the following steps. Firstly, use Lagrangian duality to reformulate the projection problem~\eqref{eq:projection} as:
\begin{align}
    \max_{\xi, \psi} & -\beta \xi + \psi - \sum_{s' \in \m{S}} \hat{P}_{s, a, s'} \phi^*(-\xi b_{s'} +  \psi) \\
    \text{s.t. } & \xi \in \R^+, \psi \in \R.
\end{align}
Then, recalling that $\phi^*(z) = \max\l\{1+\frac{z}{2}, 0\r\}^2 - 1$, we wish to eliminate the max operator in order to make the model tractable. In order to do so, we observe that at optimality, we necessarily have that $\phi^*(-\xi b_{s'} + \psi) = -1$ holds for exactly $\hat{S}$ values of $s'$, for some $\hat{S} \in \{0,\dots,S\}$. In order to find the optimal solution, we can therefore solve the model resulting from enforcing each value of $\hat{S}$ explicitly, and select the solution with the best objective value. In order to do so, w.l.o.g.\ we sort the elements of $\bm{b}$ so that they are non-increasing. Then, for each $\hat{S} \in \{0,\dots,S\}$ we create a subproblem of the reformulated projection problem by constraining $\xi, \psi$ to enforce that $\phi^*(-\xi b_{s'} + \psi) = -1$ for each $s' \in \{1,\dots,\hat{S}\}$. The final term in the objective function becomes:
\begin{align*}
    - \sum_{s' \in \m{S}} \hat{P}_{s, a, s'} \phi^*(-\xi b_{s'} +  \psi) &= \sum_{s' = 1}^{\hat{S}} \hat{P}_{s, a, s'} - \sum_{s' = \hat{S} + 1}^S  \hat{P}_{s, a, s'}\l(\l(1+\frac{-\xi b_{s'} + \psi}{2}\r)^2 - 1\r)\\
    &= \sum_{s' = 1}^{\hat{S}} \hat{P}_{s, a, s'} - \sum_{s' = \hat{S} + 1}^S  \hat{P}_{s, a, s'}\l((-\xi b_{s'} + \psi) + \frac{(-\xi b_{s'} + \psi)^2}{4}\r).
\end{align*}
Therefore, the subproblem is given by~\eqref{eq:sp1}-\eqref{eq:spn}.
\begin{align}
    \max_{\xi, \psi} &- \beta \xi + \psi + \sum_{s'=1}^{\hat{S}} \hat{P}_{s,a,s'} - \sum_{s'=\hat{S}+1}^{S} \hat{P}_{s,a,s'}\l((-\xi b_{s'} + \psi) + \frac{(-\xi b_{s'} + \psi)^2}{4}\r)\label{eq:sp1}\\
    \text{s.t. }& -\xi b_{s'} + \psi \le -2 \fa s' \in \{1, \dots, \hat{S}\}\label{eq:zeta_c1}\\
    & -\xi b_{s'} + \psi \ge -2\fa s' \in \{\hat{S}+1, \dots, S\}\label{eq:zeta_c2}\\
    &\xi \in \R_+, \psi \in \R.\label{eq:spn} 
\end{align}
Note that for $\hat{S} = 0$, constraint~\eqref{eq:zeta_c1} is redundant and can be removed. Similarly, for $\hat{S} = S$, constraint~\eqref{eq:zeta_c2} can be removed. Given this formulation, the solution of the subproblem is obtained from solving at most 3 problems, each with an analytical solution. By~\cite{ho2022robust}, for a fixed $\hat{S}$ and $\xi$, the solution of this subproblem in $\psi$ is given by:
\begin{equation}\label{eq:zeta_star}
    \psi^* = \begin{cases}
        -2 + \xi b_{\hat{S}+1} & \text{ if } H(\xi) \le -2 + \xi b_{\hat{S}+1} \\
        -2 + \xi b_{\hat{S}}  & \text{ if } H(\xi) \ge -2 + \xi b_{\hat{S}} \\
        H(\xi) & \text{ otherwise},
    \end{cases}
\end{equation}
where
\begin{equation}
    H(\xi) = \frac{2 \sum_{s'=1}^{\hat{S}}\hat{P}_{s, a,s'} + \xi \sum_{s' = \hat{S}+1}^S b_{s'} \hat{P}_{s, a, s'}}{\sum_{s'=\hat{S}+1}^S \hat{P}_{s, a, s'}}.
\end{equation}
For some border cases, we do not need to solve the problem in all 3 of these cases. In particular, we have the following special cases:
\begin{enumerate}
    \item $\hat{S} = 0$. In this case, the second case is not defined as $b_{\hat{S}}$ does not exist. 
    \item $\hat{S} = S - 1$. In this case we have:
    \begin{align}
        H(\xi) &= \frac{2 \sum_{s'=1}^{S-1}\hat{P}_{s, a,s'} + \xi  b_{S} \hat{P}_{s, a, S}}{\hat{P}_{s, a, S}}\\
        &= \frac{2(1-\hat{P}_{s, a, S})}{\hat{P}_{s, a, S}} + b_S\xi\\
        &= \frac{2}{\hat{P}_{s, a, S}} + (-2 + b_{S}\xi)\\
        & > -2 + \xi b_{S}\\
        &= -2 + \xi b_{\hat{S}+1}.
    \end{align}
    Hence, the first case in~\eqref{eq:zeta_star} is impossible. In addition,  for $\hat{S} = S - 1$ the problem becomes:
    \begin{align}
    \max_{\xi, \psi} &- \beta \xi + \psi + \sum_{s'=1}^{S-1} \hat{P}_{s,a,s'} - \hat{P}_{s,a,S}\l((-\xi b_{S} + \psi) + \frac{(-\xi b_{S} + \psi)^2}{4}\r)\\
    \text{s.t. }& -\xi b_{s'} + \psi \le -2 \fa s' \in \{1, \dots, S - 1\}\\
    & -\xi b_{S} + \psi \ge -2\\
    &\xi \in \R_+, \psi \in \R. 
    \end{align}
    In the third case of~\eqref{eq:zeta_star}, we have $\psi = \frac{2}{\hat{P}_{s, a, S}} + (-2 + b_{S}\xi)$ and so the objective function is given by:
    \begin{equation}
        -\beta\xi + \frac{2}{\hat{P}_{s, a, S}} + (-2 + b_{S}\xi) - \hat{P}_{s, a, S}\l(\l(\frac{2}{\hat{P}_{s, a, S}} - 2\r) + \frac{1}{4}\l(\frac{2}{\hat{P}_{s, a, S}} - 2\r)^2\r).
    \end{equation}
    Therefore, the derivative of the objective function is $(b_S-\beta)\xi \le 0 \fa \xi \ge 0$, since $\beta \ge \min{\bm{b}} = b_S$. Hence, $\xi$ should be set at zero if it is unconstrained.
    \item $\hat{S} = S$. In this case, the problem becomes~\eqref{eq:dont_need}-\eqref{eq:something}:
    \begin{align}
    \max_{\xi, \psi} &- \beta \xi + \psi + 1 \label{eq:dont_need} \\
    \text{s.t. }& -\xi b_{s'} + \psi \le -2 \fa s' \in \{1, \dots, S\}\label{eq:zeta_c4}\\
    &\xi \in \R_+, \psi \in \R.\label{eq:something} 
    \end{align}
    Constraint~\eqref{eq:zeta_c4} implies that $\psi \le \xi b_{S} - 2$. Since the objective is increasing in $\psi$, this means $\psi^* = -2 + \xi b_{S}$. Hence, the second case in~\eqref{eq:zeta_star} is guaranteed. Furthermore, the objective is given by $\max \xi(b_S - \beta) - 1$. Since the assumption made by~\cite{ho2022robust} is that $\min{\bm{b}} \le \beta$ and $b_S = \min{\bm{b}}$, the objective is decreasing in $\xi$ and so the optimal solution is $(\xi^*, \psi^*) = (0, -2)$. The optimal objective value is $-1$. 
\end{enumerate}
Now, in each case defined by~\eqref{eq:zeta_star}, the problem can be reformulated as a univariate program with one constraint. In the first case, \cite{ho2022robust} show that the model becomes:
\begin{align}
    \max_{\xi} &\l\{- \beta \xi + \xi b_{\hat{S}+1} - 2 + \sum_{s'=1}^{\hat{S}} \hat{P}_{s,a,s'} - \sum_{s'=\hat{S}+1}^{S} \hat{P}_{s,a,s'}\l((-\xi b_{s'} + \xi b_{\hat{S}+1} - 2) + \frac{(-\xi b_{s'} + \xi b_{\hat{S}+1} - 2)^2}{4}\r)\r\}\\
    \text{s.t. }& \xi \ge 2 \l(\sum_{s'=\hat{S}+1}^S(b_{\hat{S}+1} - b_{s'})\hat{P}_{s, a, s}\r)^{-1}.
\end{align}
Differentiating the objective function, we find that it's derivative is given by:
\begin{equation}
    -\beta + b_{\hat{S}+1} - \sum_{s'=\hat{S}+1}^S\hat{P}_{s,a,s'}\l(b_{\hat{S}+1} - b_{s'} + \frac{1}{2}(b_{\hat{S}+1} - b_{s'})(-\xi b_{s'} + \xi b_{\hat{S}+1} - 2)\r).
\end{equation}
which can be written as:
\begin{equation}
    -\beta + b_{\hat{S}+1} - \sum_{s'=\hat{S}+1}^S\hat{P}_{s,a,s'}\xi\l(b_{\hat{S}+1} - b_{s'}\r)^2.
\end{equation}
which means the globally optimal $\xi$ is given by:
\begin{equation}
    \xi^*_1 = \frac{-\beta + b_{\hat{S}+1}}{\sum_{s'=\hat{S}+1}^S\hat{P}_{s,a,s'}\l(b_{\hat{S}+1} - b_{s'}\r)^2}.
\end{equation}
In the second case, it is easy to see that $\xi^*_2$ is obtained by replacing $b_{\hat{S}+1}$ with $b_{\hat{S}}$. The model is therefore: 
\begin{align}
    \max_{\xi} &\l\{- \beta \xi + \xi b_{\hat{S}} - 2 + \sum_{s'=1}^{\hat{S}} \hat{P}_{s,a,s'} - \sum_{s'=\hat{S}+1}^{S} \hat{P}_{s,a,s'}\l((-\xi b_{s'} + \xi b_{\hat{S}} - 2) + \frac{(-\xi b_{s'} + \xi b_{\hat{S}} - 2)^2}{4}\r)\r\}\\
    \text{s.t. }& \xi \le 2 \l(\sum_{s'=\hat{S}+1}^S(b_{\hat{S}} - b_{s'})\hat{P}_{s, a, s}\r)^{-1}.
\end{align}
The corresponding globally optimal solution is given by:
\begin{equation}
    \xi^*_2 = \frac{-\beta + b_{\hat{S}}}{\sum_{s'=\hat{S}+1}^S\hat{P}_{s,a,s'}\l(b_{\hat{S}} - b_{s'}\r)^2}.
\end{equation}
In the final case, we note that:
\begin{equation}
    H'(\xi) = \frac{ \sum_{s' = \hat{S}+1}^S b_{s'} \hat{P}_{s, a, s'}}{\sum_{s'=\hat{S}+1}^S \hat{P}_{s, a, s'}}.
\end{equation}
The model then becomes:
\begin{align}
    \max_{\xi} &\l\{- \beta \xi + H(\xi) + \sum_{s'=1}^{\hat{S}} \hat{P}_{s,a,s'} - \sum_{s'=\hat{S}+1}^{S} \hat{P}_{s,a,s'}\l((-\xi b_{s'} + H(\xi)) + \frac{(-\xi b_{s'} + H(\xi))^2}{4}\r)\r\}\\
    \text{s.t. }& \xi \le 2 \l(\sum_{s'=\hat{S}+1}^S(b_{\hat{S}+1} - b_{s'})\hat{P}_{s, a, s}\r)^{-1}\\
    & \xi \ge 2 \l(\sum_{s'=\hat{S}+1}^S(b_{\hat{S}} - b_{s'})\hat{P}_{s, a, s}\r)^{-1} \label{eq:case3_LB}
\end{align}
The derivative of the objective is given by:
\begin{equation}
    -\beta + H'(\xi) - \sum_{s'=\hat{S}+1}^S\hat{P}_{s,a,s'}\l(H'(\xi) - b_{s'} + \frac{1}{2}(H'(\xi) - b_{s'})(-\xi b_{s'} + H(\xi))\r).
\end{equation}
From the same steps as for the first case, this leads to:
\begin{equation}
    \xi^*_3 = \frac{-\beta + H'(\xi)}{\sum_{s'=\hat{S}+1}^S\hat{P}_{s,a,s'}\l(H'(\xi) - b_{s'}\r)^2}.
\end{equation}
Then solving the problem in each case corresponds to checking if the optimal $\xi$ lies within the allowed range, and selecting one of the bounds if it does not.

\subsection{Reformulation of projection problem}

As shown by~\cite{ho2022robust}, a general projection problem can be reformulated as:
\begin{align}
    \max_{\psi, \xi} &- \beta\xi + \psi - \sum_{s' \in \m{S}} \hat{P}_{s, a, s'}\phi^*(-\xi b_{s'} + \psi)\\
    \text{s.t. } & \xi \in \R_+, \psi \in \R.
\end{align}
For the modified $\chi^2$ distance, we have $\phi^*(z) = \max\l\{1+\frac{z}{2}, 0\r\}^2 - 1$, or equivalently $\phi^*(z) = \frac{1}{4}\max\l\{z + 2, 0\r\}^2 - 1$. Hence, we can represent $\phi^*(-\xi b_{s'} + \psi)$ via:
\begin{align}
    \zeta_{s'} &\ge -\xi b_{s'} + \psi + 2 \fa s' \in \m{S} \\
    \zeta_{s'} &\ge 0 \fa s' \in \m{S}\\
    u_{s'} &\ge \frac{1}{4}\zeta^2_{s'} \fa s' \in \m{S}.
\end{align}
Then, the model becomes:
\begin{align}
    \max_{\xi, \psi, \bm{\zeta}, \bm{u}} &- \beta\xi + \psi - \sum_{s' \in \m{S}} \hat{P}_{s, a, s'}(u_{s'} - 1)\\
    \text{s.t. } &\zeta_{s'} \ge -\xi b_{s'} + \psi + 2 \fa s' \in \m{S} \\
    &u_{s'} \ge \frac{1}{4}\zeta^2_{s'} \fa s' \in \m{S}\\
    &\zeta_{s'} \ge 0 \fa s' \in \m{S}\\
    & \xi \in \R_+, \psi \in \R.
\end{align}

\section{A newsvendor model incorporating backorder costs}\label{sec:X_reformulation}

Suppose that action $a$ is taken when in state $s$ and assume that $b'$ now represents a backorder cost per unit of unmet demand. For a given realisation $x$ of the demand random variable $X_{s, a}$, we define the one-period reward incorporating backorder costs as:
\begin{equation}
    r'_{s, a, x} = c\min\{x, \bar{s}\} - wa - h(\bar{s} - \min\{x, \bar{s}\}) - b'\max\{x - \bar{s}, 0\}.
\end{equation}
In addition, let $\bm{P}'_{s, a} = \l(P'_{s, a, x}\r)_{x \in \m{X}_{s, a}}$ represent a (non-parametric) candidate for the distribution of $X_{s, a}$. We can then formulate the non-parametric robust Bellman update as:
\begin{equation}\label{eq:backorder_bell}
    v^{n+1}_s = \max_{\bm{\pi}_s \in \Delta_A} \min_{\bm{P}'_s \in \m{P}'_s} \sum_{a \in \m{A}} \pi_{s, a} \sum_{x \in \m{X}_{s, a}} P'_{s, a, x}\l(r'_{s, a, x} + \gamma v^n_{g(x|s, a)}\r) \fa s \in \m{S},
\end{equation}
where $\m{P}'_s$ is an ambiguity set for the true distribution of $X_{s, a}$ (not the true transition distribution). This set can be defined using $\phi$-divergences as follows:
\begin{equation}
    \m{P}'_s = \l\{\bm{P}'_s \in \Delta_{|\m{X}_{s, 1}|} \times \dots \times \Delta_{|\m{X}_{s, A}|}: \sum_{a \in \m{A}} d_a(\bm{P}'_{s, a}, \hat{\bm{P}'}_{s, a}) \le \kappa \r\},
\end{equation}
where $\hat{\bm{P}'}_{s, a} = \l(f_{X_{s, a}}(x| \bm{\hat{\theta}})\r)_{x \in \m{X}_{s, a}}$, for example. Similarly, we can formulate the parametric update problem as:
\begin{equation}
    v^{n+1}_s = \max_{\bm{\pi}_s \in \Delta_A} \min_{\bm{\theta}_s \in \Theta_s} \sum_{a \in \m{A}} \pi_{s, a} \sum_{x \in \m{X}_{s, a}} f_{X_{s,a}}(x|\bm{\theta}_{s, a})\l(r'_{s, a, x} + \gamma v^n_{g(x|s, a)}\r) \fa s \in \m{S}.
\end{equation}
In these formulations, we could simplify the terms relating to backorder costs as follows:
\begin{align}
    \sum_{x \in \m{X}_{s, a}} P'_{s, a, x} \max\l\{x-\bar{s}, 0\r\} &= \sum_{x = \bar{s} + 1}^{|\m{X}_{s,a}|} P'_{s, a, x}(x-\bar{s}) \label{eq:nonpara_backorder}\\
    \sum_{x \in \m{X}_{s, a}} f_{X_{s,a}}(x|\bm{\theta}_{s, a}) \max\l\{x-\bar{s}, 0\r\} &= \sum_{x = \bar{s} + 1}^{|\m{X}_{s,a}|} f_{X_{s,a}}(x|\bm{\theta}_{s, a})(x-\bar{s}).\label{eq:para_backorder}
\end{align}
If we have infinite support demands, i.e.\ $|\m{X}_{s, a}| \ = \infty$, then this implies that an infinite number of decision variables are required for the non-parametric model. This means that a completely different treatment is required. In many cases, however, the parametric expression can be further simplified. For example, if the demand random variable is $X_{s,a} \sim \text{Pois}(\lambda_{s, a})$, then we have:
\begin{align}
    \sum_{x = \bar{s}}^{|\m{X}_{s,a}|} f_{X_{s,a}}(x|\bm{\theta}_{s, a})(x-\bar{s}) &= \sum_{x = \bar{s} + 1}^{\infty}\frac{\lambda_{s, a}^x \exp(-\lambda_{s, a})}{x!}(x - \bar{s})\\
    &= \lambda_{s, a}\sum_{x=\bar{s}+1}^{\infty}\frac{\lambda^{x-1}_{s, a} \exp(-\lambda_{s, a})}{(x-1)!} - \bar{s}\l(1-\sum_{x=0}^{\bar{s}} \frac{\lambda_{s, a}^x \exp(-\lambda_{s, a})}{x!}\r)\\
    &= \lambda_{s, a}\sum_{x=\bar{s}}^{\infty}\frac{\lambda^{x}_{s, a} \exp(-\lambda_{s, a})}{x!} - \bar{s} \l(1-F_{X_{s,a}}(\bar{s} |\lambda_{s, a})\r)\\
    &= \lambda_{s,a}\l(1-F_{X_{s, a}}(\bar{s} - 1 |\lambda_{s, a})\r) - \bar{s} \l(1-F_{X_{s,a}}(\bar{s} |\lambda_{s, a})\r),
\end{align}
which only involves finite sums. Without incorporating further information on the true distribution of $X_{s, a}$ such as its moments, the expression in~\eqref{eq:nonpara_backorder} cannot be simplified further. The infinite number of variables required means that the algorithms in this paper are not applicable to the robust Bellman update problem in~\eqref{eq:backorder_bell}.

\end{appendices}

\end{document}